\documentclass[12pt,reqno]{amsart}
\usepackage{amsmath, amssymb, amsthm}
\usepackage{epsfig, times, stackrel, bm}
\usepackage{color}
\usepackage{esint}
\usepackage[margin=0.95in]{geometry}
\usepackage{hyperref}

	\numberwithin{equation}{section}
	
\newcommand{\field}[1]{\mathbb{#1}}

\newcommand{\R}{\field{R}}    % real numbers
\newcommand{\C}{\field{C}}    % complex numbers
\def\be{\begin{equation}}
  \def\ee{\end{equation}}

\newtheorem{theorem}{Theorem}
\newtheorem{proposition}{Proposition}

\title[Nonlinear stability of periodic waves in the Benjamin--Ono equation]{Nonlinear stability of periodic waves \\ in the Benjamin--Ono equation}

\date{\today}

\author[J. Chen]{Jinbing Chen}
\address[J. Chen]{School of Mathematics, Southeast University, Nanjing, Jiangsu 210096, People’s Republic of China}
\email{cjb@seu.edu.cn}

\author[C. Grant]{Colin Grant}
\address[C. Grant]{Department of Mathematics and Statistics, McMaster University, Hamilton, Ontario, Canada, L8S 4K1}
\email{grantc12@mcmaster.ca}

\author[D. E. Pelinovsky]{Dmitry E. Pelinovsky}
\address[D. E. Pelinovsky]{Department of Mathematics and Statistics, McMaster University, Hamilton, Ontario, Canada, L8S 4K1}
\email{pelinod@mcmaster.ca}

\begin{document}
\maketitle

\begin{abstract}
We construct the Lyapunov functional for the periodic waves of the Benjamin--Ono equation based on the conservation of momentum, energy, and a higher-order energy, which exists due to integrability. Positivity of the second variation of the Lyapunov functional for perturbations 
in Sobolev space $H^1(\mathbb{R})$ is established by 
utilizing the basis properties of eigenfunctions of the linearized 
Benjamin--Ono equation in $L^2(\mathbb{R})$. 
Coercivity of the Lyapunov functional for co-periodic and multiple-periodic perturbations yields the nonlinear orbital stability of the periodic waves in the periodic domain.
\end{abstract}

\section{Introduction}

The Benjamin--Ono (BO) equation is a popular model to describe dynamics 
of long internal waves in a deep fluid \cite{Ben67,Ono75}. Although 
it has been known to be formally integrable with the inverse scattering method for a long time \cite{AF83,KM98}, it was only recently that many groundbreaking results were obtained in the mathematical analysis of this model. 

We are particularly concerned with the stability of traveling wave solutions with spatially decaying or spatially periodic profiles. These are called {\em solitary waves} and {\em periodic waves}, respectively. 
We take the BO equation in the standard form: 
\begin{equation}
\label{BO}
u_t + 2 u u_x + H(u_{xx}) = 0, 
\end{equation}
where $u(x,t) : \mathbb{R}\times \mathbb{R}\to \mathbb{R}$ and 
$H$ is the Hilbert transform defined in the Fourier basis by 
$H(e^{i \xi x}) = i {\rm sgn}(\xi) e^{i \xi x}$, $\xi \in \mathbb{R}$. 
The solitary wave is characterized by algebraically decaying tails:
\begin{equation}
\label{sol-wave}
u(x,t) = \frac{2c}{1 + c^2 \xi^2},
\end{equation}
where $\xi = x - c t - \xi_0$ is the traveling wave coordinate with 
an arbitrary wave speed $c > 0$ and the phase parameter $\xi_0 \in \mathbb{R}$. The periodic wave can be expressed in terms of elementary functions:
\begin{equation}
\label{per-wave}
u(x,t) = \frac{k \sinh \phi}{\cos(k\xi) + \cosh \phi},
\end{equation}
where $\xi$ and $\xi_0$ are the same as in (\ref{sol-wave}), the 
two parameters $k, \phi > 0$ are arbitrary, and the wave speed $c$ is uniquely found from $c = k \coth \phi$. 
Without loss of generality, we can set $\xi_0 = 0$ due to the translational invariance of the BO equation (\ref{BO}). Due to the scaling invariance of the BO equation (\ref{BO}) under the transformation $u(x,t) \to k u(kx,k^2t)$, one can set $k=1$ and restrict 
the periodic solutions to the period $2\pi$, for which the requirement $\phi > 0$ implies $c > 1$. For computational transparency, we will leave the scaling parameter $k > 0$ arbitrary. 

\subsection{Motivations and background references}

Regarding the solitary waves (\ref{sol-wave}), uniqueness of spatially 
decaying traveling wave solutions of the BO equation for each wave speed $c > 0$ was shown in \cite{AT91}. Linear stability of a single solitary 
wave was shown in \cite{ChenKaup1980} by constructing a complete basis of eigenfunctions in the linearized BO equation. The complete basis of eigenfunctions of  the linearized BO equation was extended to the multi-soliton solutions \cite{MatsunoKaup,MatsunoKaup2}, which implies their linear stability.

For the nonlinear stability analysis, we recall the first three conserved quantities for the BO equation (\ref{BO}):
\begin{align}
M(u) &= \frac{1}{2} \int u^2 dx, \label{mass} \\
E(u) &= -\frac{1}{2} \int u H(u_x) dx - \frac{1}{3} \int u^3 dx, \label{energy} \\
R(u) &= \frac{1}{2} \int u_x^2 dx + \frac{3}{4} \int u^2 H(u_x) dx + \frac{1}{4} \int u^4 dx. \label{higher-energy}
\end{align}
These quantities are well-defined in the Sobolev spaces $H^0(\R)$, $H^{1/2}(\R)$, and $H^1(\mathbb{R})$ respectively. The list of conserved quantities can be continued to Sobolev spaces of higher regularity $H^{N/2}(\R)$ with $N \in \mathbb{N}$ due to integrability of the BO equation.

Nonlinear orbital stability of the solitary wave (\ref{sol-wave}) in $H^{1/2}(\R)$ was shown in \cite{Albert87,Bona83} by minimizing the energy
$E(u)$ subject to fixed momentum $M(u)$. The method of using more conserved 
quantities for the Lyapunov functional has been employed to show the nonlinear stability of multi-soliton solutions. Nonlinear orbital stability of two-soliton solutions in $H^1(\R)$ was shown in \cite{NL2006} by minimizing the higher-order energy $R(u)$ subject to fixed momentum $M(u)$ and energy $E(u)$. The same strategy was extended to show the nonlinear stability of $N$-solitons in Sobolev space $H^{N/2}(\R)$ by proving their local constrained minimization properties in a linear combination of higher-order conserved quantities \cite{Matsuno2006}. The recent work \cite{LanWang2025} connected the proof of nonlinear stability of multi-solitons with the completeness of eigenfunctions of the linearized BO equation.
 
The nonlinear stability of multi-solitons in Sobolev spaces of low regularity $H^s(\R)$ with $s > -\frac{1}{2}$ was obtained in  \cite{BKV2025} by using a new variational characterization of multi-solitons (avoiding the energy conserved quantities) and by extending the spectral analysis of Lax operators to spaces of low regularity. A more recent work \cite{GassotMiller2026} extends the nonlinear stability of multi-solitons to the level of their global asymptotic stability and gives a proof of the so-called soliton resolution conjecture. The long-time scattering of time-dependent solutions of the BO equation to the superposition of several solitary waves requires sufficient regularity and spatial decay of the initial data and is based on an explicit representation of solutions obtained in \cite{Gerard2023}. The latest work \cite{BKV2026} gives the proof of the asymptotic stability of multi-solitons in $L^2(\R)$ without imposing the spatial decay condition on the initial data. 

Regarding the periodic waves (\ref{per-wave}), their linear stability was shown in \cite{MS94}  by using the method of squared eigenfunctions. Modulation theory for the multi-periodic waves was developed in \cite{DK-91} and  \cite{Matsuno1998} based on Whitham's modulation theory. Rigorous characterization of the Lax spectrum in spaces of co-periodic and decaying perturbations was developed recently in \cite{GK-21,GKT-20,GKT-22}. This construction was used in \cite{ChenPel24} to identify the bright and dark breathers (solitary waves on the spatially periodic background which exhibit temporally periodic dynamics) according to the eigenvalues in the semi-infinite and finite spectral gaps of the Lax spectrum. The unique property of the breathers dynamics in the BO equation is the zero phase shift in the soliton-wave interaction \cite{ChenPel24} similarly to the zero phase shift in multi-soliton solutions \cite{Matsuno1979}.

{\em The purpose of this work is to establish the nonlinear stability of the 
periodic waves (\ref{per-wave}) by using the higher-order variational characterization of the periodic wave profiles.} The characterization is based on a linear combination of the momentum $M(u)$, the energy $E(u)$, and the higher-order energy $R(u)$ given by (\ref{mass}), (\ref{energy}), and (\ref{higher-energy}), which yields the Lyapunov functional 
in Sobolev space $H^1(\mathbb{R})$. 

The idea of finding a Lyapunov functional for the periodic waves 
from the higher-order energies has been previously used for 
the Korteweg--de Vries (KdV) equation \cite{KD,LHP,ND2010} as well as 
for the defocusing versions of the modified KdV equation \cite{Decon} 
and the nonlinear Schr\"{o}dinger equation \cite{DN,GP}. An additional challenge for the BO equation appears from the nonlocal terms in the energy 
functionals (\ref{energy}) and (\ref{higher-energy}). Nevertheless, we show how to adapt the squared eigenfunctions of the linearized BO equation from \cite{MS94} to prove the Lyapunov property for periodic waves in the BO equation. 

The method can be extended for the $N$-periodic solutions at the price of using a higher-order variational formulation in Sobolev spaces of a higher regularity $H^{(N+1)/2}(\mathbb{R})$ compared to the standard variational formulation in Sobolev space $H^{N/2}(\mathbb{R})$. It is not clear if any of the very recent treatments of multi-solitons of the BO equation can be used to obtain similar nonlinear stability results for periodic waves, 
e.g. working with the commuting flows \cite{BKV2025,BKV2026} or utilizing the explicit solution formula \cite{GassotMiller2026,Gerard2023}. Given the sharp well-posedness of the BO equation in $H^s(\mathbb{R})$ for $s > -\frac{1}{2}$ \cite{GKT-23,KLV}, it is natural to look for the nonlinear stability of periodic waves in Sobolev spaces of low regularity, which is an open problem.

\subsection{Main results and applications}

In what follows, we refer to $U$ as the spatial profile of the periodic waves $u(x,t) = U(\xi)$ with $\xi = x - ct$. It is defined by the explicit 
formula (\ref{per-wave}) as the $\frac{2\pi}{k}$-periodic function. 
Regarding perturbations to the periodic profile $U$, they can be 
considered in Sobolev spaces of co-periodic, multiple-periodic, or decaying functions. These spaces are usually denoted 
as 
$$
H^s_{\rm per} = H^s\left(\mathbb{R} / \frac{2\pi}{k} \mathbb{Z}\right), \quad 
H^s_{\rm Nper} = H^s\left(\mathbb{R} / \frac{2\pi N}{k} \mathbb{Z}\right), \quad H^s(\mathbb{R})
$$ 
for $N \in \mathbb{N}$ and $s \in \mathbb{R}$, and we restrict these spaces to real-valued functions. 

In what follows, we use the first and second derivatives 
of the functionals (\ref{mass}), (\ref{energy}), and (\ref{higher-energy}):
\begin{align*}
\begin{cases}
M'(u) = u, \\
E'(u) = -H(u_x) - u^2, \\
R'(u) = -u_{xx} + \frac{3}{2} u H(u_x) + \frac{3}{2} H(u u_x) + u^3, 
\end{cases}
\end{align*}
and
\begin{align*}
\begin{cases}
M''(u) = 1, \\
E''(u) = -H \partial_x - 2u, \\
R''(u) = -\partial_x^2 + \frac{3}{2} H(u_x) 
+ \frac{3}{2} u H \partial_x + \frac{3}{2} H \partial_x (u \cdot ) + 3 u^2.
\end{cases}
\end{align*}
Since $U$ is the periodic profile, these derivatives are well defined for perturbations in Sobolev spaces of co-periodic and multiple-periodic functions. By subtracting the densities of $M(u)$, $E(u)$, and $R(u)$ computed at the periodic profile $U$ from the corresponding densities computed at $U + v$, these derivatives are also well-defined for the perturbation $v$ in Sobolev spaces of decaying functions provided that $U$ is a critical point of the corresponding energy functionals.

The following theorem presents the main result of this work.

\begin{theorem}
	\label{th-main}
The Lyapunov functional for the periodic wave with the spatial 
profile $U$ is given by  
	\begin{equation}
	\label{Lyapunov}
	S(u) = R(u) + \frac{1}{2} (3c-k) E(u) + \frac{1}{4} (3c + k) (c-k) M(u),
	\end{equation}
in the sense that $S'(U) = 0$ and $\langle S''(U) v, v \rangle_{L^2(\R)} \geq 0$ for every $v \in H^1(\mathbb{R})$. 
The periodic profile $U$ is a local minimizer of $S(u)$ 
in $H^1_{\rm Nper}$ for any fixed $N \in \mathbb{N}$, which is only degenerate due to the translational symmetry in the sense that $\langle S''(U) v, v \rangle_{L^2_{\rm Nper}} = 0$ if and only if $v = \alpha U'$ in $H^1_{\rm Nper}$ for arbitrary $\alpha \in \mathbb{R}$.
\end{theorem}

An immediate consequence of Theorem \ref{th-main} is 
the orbital stability of periodic waves with respect to co-periodic or multiple-periodic perturbations. 

\begin{theorem}
	\label{th-orbit}
For every fixed $N \in \mathbb{N}$, the periodic wave with the spatial profile $U$ is orbitally stable in the time evolution of the BO equation in $H^1_{\rm Nper}$. To be precise, for every $\varepsilon > 0$, there exists $\delta > 0$ such that for every $u_0 \in H^1_{\rm Nper}$ such that $\| u_0 - U \|_{H^1_{\rm Nper}} < \delta$, the unique solution 
	$u \in C(\mathbb{R},H^1_{\rm Nper})$ satisfies 
\begin{equation}
\label{bound-on-sol}
	\sup_{t \in \mathbb{R}} \inf_{\theta \in \mathbb{R}} \| u(\cdot,t) - U(\cdot - \theta) \|_{H^1_{\rm Nper}} < \varepsilon.
\end{equation}
\end{theorem}

The result of Theorem \ref{th-orbit} follows from the coercivity of 
the quadratic form $\langle S''(U) v, v \rangle_{L^2_{\rm Nper}}$ 
for $v \in H^1_{\rm Nper} |_{\{ U'\}^{\perp}}$, where the constraint 
is satisfied by the choice of the translational parameter $\theta$ in 
(\ref{bound-on-sol}). The coercivity is obtained using the Lyapunov property of the functional $S(u)$ due to the spectral gap between the 
zero eigenvalue of $S''(U)$ in $L^2_{\rm Nper}$ with ${\rm Ker}_{L^2_{{\rm Nper}}}(S''(U)) = {\rm span}(U')$ and the rest of the spectrum of $S''(U)$ in $L^2_{\rm Nper}$. The proof of Theorem \ref{th-orbit} is standard and follows the previous works \cite{KD,GP}, hence it is omitted.

As is well-known, the problem of nonlinear stability of periodic waves with respect to spatially decaying perturbations has been open in Hamiltonian systems 
due to the lack of the spectral gap of $S''(U)$ in $L^2(\mathbb{R})$. 
This problem has been recently solved in \cite{BukiedaRijk2026} 
for the KdV equation, where variational arguments, Floquet–Bloch theory, and Duhamel-based estimates were combined with
the spatiotemporal modulations of the periodic waves. The translational modulations were shown in \cite[Theorem 1.4]{BukiedaRijk2026} 
to be growing in time in the $L^2(\mathbb{R})$ norm due to the 
phase shift accumulating after multiple soliton-wave interactions. See also \cite{BGR} for a related work.

Compared to the KdV equation, where the phase shift was analyzed in \cite{HMP23}, the BO equation admits no 
phase shift in the soliton-wave interactions \cite{ChenPel24}. 
Therefore, it is an interesting open question to extend the proof 
of nonlinear stability of periodic waves under the spatially decaying perturbations to the BO equation and to study if the translational 
modulations of the periodic waves can be bounded in space and time.

\subsection{Methods and organization of the paper}

To prove Theorem \ref{th-main}, we build a complete basis of eigenfunctions of the linearized BO equation 
\begin{equation}
\label{BO-lin}
q_t + 2 (uq)_x + H(q_{xx}) = 0, 
\end{equation}
where $u(x,t) = U(x-ct)$ is the traveling wave solution with the spatial profile (\ref{per-wave}). The construction of eigenfunctions 
was developed in \cite{MS94} based on the formal dressing method used previously in \cite{KSF}. Orthogonality and completeness of eigenfunctions were shown in Appendices B and C of \cite{MS94} with a brief discussion. For transparency of presentation, we are going to recover the main 
conclusions of \cite{MS94} with a complete study, where the construction, orthogonality, and completeness of eigenfunctions are rigorously proven 
with explicit details. These results are included in Sections 
\ref{sec-2}, \ref{sec-3}, and \ref{sec-4}, where we use the following generalized identities which follow from the theory of distributions and Fourier transform:
\begin{equation}
\label{distrib-1}
\int_{\mathbb{R}} e^{i \alpha x} dx = 2\pi \delta(\alpha), \quad \alpha \in \mathbb{R}
\end{equation}
and 
\begin{equation}
\label{distrib-2}
\int_{\mathbb{R}} \frac{d}{dx} e^{i \alpha x} dx = 0, \quad \alpha \in \mathbb{R},
\end{equation}
where $\delta$ is the Dirac delta distribution.

Once the complete basis of eigenfunctions is verified, we identify the spatial profile $U$ as the relative equilibrium in two commuting flows:
\begin{align}
u_t = \partial_x E'(u)
\label{BO-flow} 
\end{align}
and
\begin{align}
u_{\tau} = \partial_x R'(u).
\label{higher-BO-flow}
\end{align}
The flow (\ref{BO-flow}) corresponds to the BO equation (\ref{BO}) 
with the time variable $t$. The flow (\ref{higher-BO-flow}) corresponds 
to the higher-order BO equation with the time variable $\tau$. Similar to the standard variational characterization, we show that $U$ is a critical point of the augmented energy $E(u) + c M(u)$ in $H^{1/2}_{\rm per}$. Due to the commuting flow, we also show that $U$ is a critical point of the augmented higher-order energy $R(u) + \gamma M(u)$ in $H^1_{\rm per}$ with $\gamma = -\frac{1}{4} (k^2 + 3 c^2)$. 

For the convexity analysis, we diagonalize the quadratic forms associated 
to the Hessian operators $L : H^1(\R) \subset L^2(\R) \to L^2(\R)$ and 
$\tilde{L} : H^2(\R) \subset L^2(\R) \to L^2(\R)$ given by 
\begin{align}
\begin{cases}
L = E''(U) + c M''(U), \\
\tilde{L} = R''(U) + \gamma M''(U).
\end{cases}
\label{Hessian}
\end{align} 
Due to the duality of $H^{-1/2}(\R)$ and $H^{1/2}(R)$ as well as 
$H^{-1}(\R)$ and $H^1(\R)$ under the $L^2(\R)$ inner product, the quadratic 
forms $\langle L q, q \rangle_{L^2(\R)}$ and $\langle \tilde{L} q, q \rangle_{L^2(\R)}$ are well-defined for every $q \in H^{1/2}(\R)$ and $q \in H^1(\R)$ respectively.

To diagonalize the linearized operators (\ref{Hessian}), we use the complete basis of eigenfunctions of the linearized 
BO equation (\ref{BO-lin}) with $u(x,t) = U(x-ct)$, as well as 
the complete basis of eigenfunctions of the adjoint linearized 
BO equation  
\begin{equation}
\label{BO-adjoint}
\psi_t + 2 u \psi_x + H(\psi_{xx}) = 0. 
\end{equation}
As a result of explicit computations, we show that the quadratic forms  $\langle L q, q \rangle_{L^2(\R)}$ and $\langle \tilde{L} q, q \rangle_{L^2(\R)}$ are not  positive semi-definite, but their linear combination  $\langle (\tilde{L} + \delta L) q, q \rangle_{L^2(\R)}$ with $\delta = \frac{1}{2} (3c - k)$ is positive semi-definite. The linear combination yields the Lyapunov functional (\ref{Lyapunov}) and thus the proof of the first assertion of Theorem \ref{th-main}. This result is described in Section \ref{sec-5}.

Finally, we extend the proof of orthogonality and completeness of squared eigenfunctions in $L^2_{\rm per}$. Using the complete basis of eigenfunctions in $L^2_{\rm per}$, we prove that ${\rm Ker}_{L^2_{{\rm per}}}(S''(U)) = {\rm span}(U')$ and the rest of the spectrum of $S''(U)$ in $L^2_{\rm per}$ is bounded away from the simple zero eigenvalue 
of $S''(U)$ in $L^2_{\rm per}$. Moreover, we extend the proof to $L^2_{\rm Nper}$ for every $N \in \mathbb{N}$ with the Floquet--Bloch theory and this yields the proof of the second assertion of Theorem \ref{th-main}.  This result is described in Section \ref{sec-6}.

\section{Eigenfunctions of the linearized BO equation} 
\label{sec-2}

The BO equation (\ref{BO}) arises as a compatibility condition 
of the following system of Lax equations \cite{AF83,KM98}:
\begin{equation}
\label{Lax-pair}
\begin{cases} 
i \varphi^+_x + \lambda( \varphi^+ - \varphi^-) + u \varphi^+ = 0, \\
i \varphi^{\pm}_t - 2 i \lambda \varphi^{\pm}_x + \varphi^{\pm}_{xx} - [\pm i u_x + H(u_x)] \varphi^{\pm} = 0,
\end{cases}
\end{equation}
where $\varphi^{\pm}(x,t) : \mathbb{R}\times \mathbb{R}\to \mathbb{C}$ are continued off the real line for $x$ as analytic functions in $\mathbb{C}^{\pm}$ and $\lambda \in \mathbb{R}$ is the spectral parameter.

It was formally found in \cite{DK-91}  and rigorously justified in \cite{GK-21,GKT-20,GKT-22} that if $u(x,t) = U(x-ct)$ is the traveling wave with the spatial profile (\ref{per-wave}), then the Lax equations (\ref{Lax-pair}) admit the purely continuous spectrum in $L^2(\R)$ located at
\begin{equation}
\label{Lax-spectrum-continuous}
\sigma_L |_{L^2(\R)} = [\lambda_0,\mu_0] \cup [0,\infty), \quad \lambda_0 = -\frac{c+k}{2}, \quad \mu_0 = -\frac{c-k}{2},
\end{equation}
see \cite[Proposition 2.2]{GK-21} and \cite[Proposition C.2]{GK-21}. We note that $\mu_0 = \lambda_0 + k$. In addition, we can write down some relations between $\lambda_0$ and $\mu_0$ obtained from (\ref{Lax-spectrum-continuous}) with $c = k \coth \phi$:
\begin{equation}
\label{lambda-mu}
\lambda_0 = -\frac{k e^{\phi}}{e^{\phi} - e^{-\phi}}, \quad \mu_0 = -\frac{k e^{-\phi}}{e^{\phi} - e^{-\phi}}, \quad \Rightarrow \quad \lambda_0 e^{-\phi} = \mu_0 e^{\phi}.
\end{equation}

To introduce the eigenfunctions, we are going to use the Hirota bilinear form for the periodic traveling wave (\ref{per-wave}) with the spatial profile $U$ given by 
\begin{equation}
\label{Hirota}
U = i \partial_x \log \frac{f'}{f} = i \left( \frac{f'_x}{f'} - \frac{f_x}{f} \right),  \quad 
\left\{ \begin{array}{l} f = 1 + e^{i k \xi + \phi}, \\
f' = 1 + e^{i k \xi - \phi} \end{array}.\right.
\end{equation}
We note that $f$ has no zeros in $\C^-$ and $f'$ has no zeros in $\C^+$ 
in the complex $x$-plane due to the constraints $k, \phi > 0$. As a result, 
$\frac{1}{f}$ is analytic in $\mathbb{C}^-$ and $\frac{1}{f'}$ is analytic in $\mathbb{C}^+$. If $\psi^{\pm}$ is analytically continued off $\R$ to $\mathbb{C}^{\pm}$ and if $\psi^{\pm}$ decays to $0$ sufficiently fast (e.g., exponentially fast) as ${\rm Im}(x) \to \pm \infty$, then $H \psi^{\pm} = \pm i \psi^{\pm}$. As a result, we can write
\begin{equation}
\label{Hirota-H}
H U = - \left( \frac{f'_x}{f'} + \frac{f_x}{f} \right) + i k,
\end{equation}
where we have represented 
$$
-\frac{i f_x}{f} = \frac{k e^{i k \xi+ \phi}}{1 + e^{ik \xi + \phi}} 
= k - \frac{k}{f} \quad \Rightarrow \quad H\left(-\frac{i f_x}{f}\right) = \frac{ik}{f} = ik -\frac{f_x}{f},
$$
in order to ensure the exponential decay of $-\frac{k}{f}$ to $0$ as ${\rm Im}(x) \to -\infty$, whereas $H(k) = 0$. No modification is needed for $\frac{i f'_x}{f'}$ since this function decays exponentially to $0$ as ${\rm Im}(x) \to +\infty$. 
 
For the two parts of the spectrum in (\ref{Lax-spectrum-continuous}), 
the two sets of eigenfunctions of the Lax equations (\ref{Lax-pair}) have been derived in \cite{ChenPel24} in the explicit form:
\begin{equation}
\varphi_I^+(\lambda;x,t) = e^{i \lambda x + i \lambda^2 t} \frac{f(x-ct)}{f'(x-ct)}, \quad \varphi_I^-(\lambda;x,t) = 0, \quad \lambda \in [0,\infty)
\label{eigenfunction-1}
\end{equation}
and for any $\lambda \in \mathbb{R}\backslash \{ \lambda_0\}$, 
\begin{equation}
\varphi_{II}^+(\lambda;x,t) = 1 - \frac{k}{(\lambda - \lambda_0) f'(x-ct)}, \quad 
\varphi_{II}^-(\lambda;x,t) = 1 - \frac{k}{(\lambda - \lambda_0) f(x-ct)}. 
\label{eigenfunction-2}
\end{equation}

Based on the computations in \cite{MS94}, we have recovered the set of four particular solutions to the adjoint linearized BO equation (\ref{BO-adjoint}) from the sets of eigenfunctions (\ref{eigenfunction-1}) and (\ref{eigenfunction-2}). We give a self-contained proof of the following proposition for ease of reading.

\begin{proposition}
	\label{prop-eigenfunctions}
	The adjoint linearized BO equation (\ref{BO-adjoint}) admits the following four solutions:
\begin{align}
\psi_{I}(\lambda;x,t) &= e^{i \lambda x + i \lambda^2 t} \frac{f(x-ct)}{f'(x-ct)} \left[ 1 + \frac{k}{(\lambda - \mu_0) f(x-ct)} \right], \label{solution-1} \\
\psi_{II}(\lambda;x,t) &= e^{-i \lambda x - i \lambda^2 t} \frac{f'(x-ct)}{f(x-ct)} \left[ 1 - \frac{k}{(\lambda - \lambda_0) f'(x-ct)} \right] \label{solution-2}
\end{align}
where $\lambda \in (0,\infty)$, and 
\begin{align}
\psi_{III}(\lambda;x,t) &= \frac{e^{-i (\lambda_0 x + \lambda_0^2 t)}}{f(x-ct)} \psi_I(\lambda;x,t), 
\label{solution-3} \\
\psi_{IV}(\lambda;x,t) &= \frac{e^{i (\mu_0 x + \mu_0^2 t)}}{f'(x-ct)} \psi_{II}(\lambda;x,t),
\label{solution-4}
\end{align}
where $\lambda \in (\lambda_0,\mu_0)$.
\end{proposition}

\begin{proof}
	By using (\ref{eigenfunction-1}) and (\ref{eigenfunction-2}), we observe that 
	\begin{align*}
	\bar{\varphi}^+_I(\lambda;x,t) &= e^{-i\lambda x - i \lambda^2 t} \frac{f'(x-ct)}{f(x-ct)} e^{2 \phi}, \\
	\bar{\varphi}^+_{II}(\lambda;x,t) &= \frac{\lambda - \mu_0}{\lambda - \lambda_0} \left[1 + \frac{k}{(\lambda - \mu_0) f(x-ct)} \right]. 
	\end{align*}
	As a result, eigenfunctions $\psi_I$ and $\psi_{II}$ given by (\ref{solution-1}) and (\ref{solution-2}) are proportional to the squared eigenfunctions of the Lax equations (\ref{Lax-pair}), namely,
	\begin{equation*}
	\psi_I(\lambda;x,t) = \frac{\lambda - \lambda_0}{\lambda - \mu_0} \varphi^+_I(\lambda;x,t) \bar{\varphi}^+_{II}(\lambda;x,t), \quad 
	\psi_{II}(\lambda;x,t) = e^{-2 \phi} \bar{\varphi}^+_I(\lambda;x,t) \varphi^+_{II}(\lambda;x,t).
	\end{equation*}
Since we have 
	\begin{equation}
	\label{conjugation-1}
\psi_{II} = e^{-2 \phi} \frac{\lambda - \mu_0}{\lambda - \lambda_0} \bar{\psi}_I,
	\end{equation}
and $u$ is real, it follows that $\psi_{II}$ is a solution of the adjoint linearized equation (\ref{BO-adjoint}) if $\psi_I$ is a solution of the same equation. 

On the other hand, eigenfunctions  $\psi_{III}$ and $\psi_{IV}$ given by
(\ref{solution-3}) and (\ref{solution-4}) are constructed 
nontrivially in terms of the squared eigenfunctions $\psi_I$ and $\psi_{II}$, namely
\begin{align*}
\psi_{III}(\lambda;x,t) &= -\frac{1}{k} \psi_I(\lambda;x,t) \lim_{\lambda \to \lambda_0} (\lambda - \lambda_0) \psi_{II}(\lambda;x,t), \\
\psi_{IV}(\lambda;x,t) &= \frac{1}{k} \psi_{II}(\lambda;x,t) \lim_{\lambda \to \mu_0} (\lambda - \mu_0) \psi_{I}(\lambda;x,t).
\end{align*}
By using $\mu_0 = \lambda_0 + k$ and $2 \lambda_0 + k = -c$, we obtain 
\begin{align*}
\lim_{\lambda \to \mu_0} (\lambda - \mu_0) \psi_I(\lambda;x,t) &= \frac{k e^{i \mu_0 x + i \mu_0^2 t}}{f'(x-ct)} = \frac{k e^{i \lambda_0 x + i \lambda_0^2 t + i k x + i k (2\lambda_0 + k)t}}{1 + e^{i k (x-ct) - \phi}} \\
&= \frac{k e^{i \lambda_0 x + i \lambda_0^2 t + \phi}}{1 + e^{-i k (x-ct) + \phi}} = \frac{k e^{i \lambda_0 x + i \lambda_0^2 t + \phi}}{\bar{f}(x-ct)} \\
&= -e^{\phi} \lim_{\lambda \to \lambda_0} (\lambda - \lambda_0) \bar{\psi}_{II}(\lambda;x,t),
\end{align*}
which implies that 
	\begin{equation}
\label{conjugation-2}
\psi_{IV} = e^{-\phi} \frac{\lambda - \mu_0}{\lambda - \lambda_0} \bar{\psi}_{III}.
	\end{equation}
Since $u$ is real, it follows that $\psi_{IV}$ is a solution of the adjoint linearized equation (\ref{BO-adjoint}) if $\psi_{III}$ is a solution of the same equation. 

Hence, it suffices to prove that $\psi_I$ and $\psi_{III}$ satisfy the adjoint linearized equation (\ref{BO-adjoint}).

\medskip

\underline{Solution $\psi_I$.} Since $f'$ has no zeros in $\mathbb{C}^+$ in the complex $x$-plane, it follows from (\ref{solution-1}) that $\psi_I$ is analytic in $\C^+$. Furthermore, because $\lambda > 0$, $e^{i \lambda x}$ and hence $\psi_I$ decays to $0$ exponentially as ${\rm Im}(x) \to +\infty$. As a result, we have $H(\psi_I) = i \psi_{I}$ and the adjoint linearized equation (\ref{BO-adjoint}) for $\psi_I$ becomes
$$
\psi_{I,t} + 2 u \psi_{I,x} + i \psi_{I,xx} = 0.
$$
Substituting $\psi_I \eqsim \varphi^+_I \bar{\varphi}^+_{II}$ into this equation and dropping the multiplicative constant, we obtain the equivalent equation in the form:
\begin{equation}
\label{tech-eq-1}
0=\varphi_{I,t}^+\bar\varphi^+_{II}+\varphi_{I}^+\bar\varphi_{II,t}^+ +2u(\varphi_{I,x}^+\bar\varphi^+_{II}+\varphi_{I}^+\bar\varphi_{II,x}^+)+i(\varphi_{I,xx}^+\bar\varphi_{II}^+ + 2\varphi_{I,x}^+\bar\varphi_{II,x}^+ +\varphi_I^+\bar\varphi_{II,xx}^+).
\end{equation}
By using the second equation of the system (\ref{Lax-pair}) for $\varphi_{II}^+$ after complex conjugation,
\begin{align*}
\bar\varphi_{II,t}^+ &= -i\bar\varphi_{II,xx}^+ + \bar\varphi_{II}^+(u_x+iH(u_x))+2\lambda\bar\varphi_{II,x}^+,
\end{align*}
we rewrite (\ref{tech-eq-1}) as 
\begin{align*}
0 &= \varphi_{I,t}^+ \bar\varphi_{II}^+ + 2 u (\varphi_{I,x}^+ \bar\varphi_{II}^+ + \varphi_I^+ \bar\varphi_{II,x}^+) + i (\varphi_{I,xx}^+ \bar\varphi_{II}^+ + 2 \varphi_{I,x}^+ \bar\varphi_{II,x}^+) \\
& \qquad + \varphi_I^+ (2\lambda\bar\varphi_{II,x}^+ + (u_x +iH(u_x)) \bar\varphi_{II}^+) \\
&= \bar\varphi_{II,x}^+ \left( 2u \varphi_I^+ + 2i \varphi_{I,x}^+ + 2 \lambda \varphi_I^+ \right)  \\
& \qquad + \bar\varphi_{II}^+ \left( \varphi_{I,t}^+ 
+ 2 u \varphi_{I,x}^+ + i \varphi_{I,xx}^+ +  (u_x +iH(u_x))\varphi_I^+
\right).
\end{align*}
By using the system (\ref{Lax-pair}) for $\varphi_I^+$ 
with $\varphi_I^- = 0$,
\begin{align*}
\varphi_{I,x}^+ &= i (\lambda+u) \varphi_{I}^+ , \\
\varphi^{+}_{I,t} &= 2 \lambda \varphi^{+}_{I,x} + i \varphi^{+}_{I,xx} + [u_x - i H(u_x)] \varphi^{+}_I,
\end{align*}
we transform the remaining equation to the form:
\begin{align*}
0&=2\lambda\varphi_{I,x}^+ + 2 u_x \varphi_I^+ + 2 u \varphi_{I,x}^+
+ 2 i \varphi_{I,xx}^+\\
&= 2 \partial_x \left( i \varphi_{I,x}^+ + \lambda \varphi_{I}^+ + u \varphi_I^+ \right),
\end{align*}
which is satisfied because $\varphi_{I,x}^+ = i (\lambda+u) \varphi_{I}^+$.

\medskip

\underline{Solution $\psi_{III}$.} Since $\mu_0 = \lambda_0 + k$ and $\mu_0 = \lambda_0 e^{-2\phi}$, we can rewrite $\psi_{III}$ in the form:
\begin{align}
\psi_{III} & = e^{i (\lambda - \lambda_0) x + i (\lambda^2 - \lambda^2_0) t} \left[
\frac{1}{f'} + \frac{k}{(\lambda - \mu_0) f f'} \right] \notag \\
& = \frac{e^{i (\lambda - \lambda_0) x + i (\lambda^2 - \lambda^2_0) t} }{\lambda - \mu_0} \left[ 
\frac{\lambda}{f'} - \frac{\lambda_0}{f} \right].
\label{temp-psi}
\end{align}
We recall that $\frac{1}{f'}$ and $\frac{1}{f}$ are analytic in $\mathbb{C}^\pm$. Since $\lambda \in (\lambda_0,\mu_0)$, the function $\frac{1}{f'} e^{ i (\lambda - \lambda_0)x}$ decays to $0$ as ${\rm Im}(x) \rightarrow +\infty$ and the function $\frac{1}{f} e^{ i (\lambda - \lambda_0)x} \sim e^{i(\lambda - \mu_0) x - \phi}$ decays to $0$ as ${\rm Im}(x) \rightarrow -\infty$. Therefore, we obtain
\begin{equation} \label{c-1}
H(\psi_{III}) = \frac{e^{i (\lambda - \lambda_0) x + i (\lambda^2 - \lambda^2_0) t}}{\lambda - \mu_0}  \left[
\frac{i \lambda}{f'} + \frac{i \lambda_0}{f} \right].
\end{equation}
Additionally, it follows from  \eqref{Hirota} that $u$ can be rewritten in the form
\begin{equation}\label{c-2}
u = k \left(\frac{1}{f'} - \frac{1}{f}\right).
\end{equation}
We now substitute (\ref{temp-psi}), (\ref{c-1}), and (\ref{c-2}) into 
each term of the adjoint linearized equation (\ref{BO-adjoint}) and drop the multiplicative constant. For brevity of writing, we denote the common factor $e^{i (\lambda - \lambda_0) x + i (\lambda^2 - \lambda^2_0) t}$ by $\mathcal{E}$. Using the quotient rule, we obtain with the aid of $\lambda_0 e^{-\phi} = \mu_0 e^{\phi}$, $\mu_0 = \lambda_0 + k$, and $c = -2 \lambda_0 - k$ that 
\begin{align*}
\frac{\partial}{\partial t} \left[ \frac{\mathcal{E}}{f'} \right]
& = i \mathcal{E} \left[ \frac{\lambda^2 - \mu_0^2}{f'}
- \frac{ k c }{(f')^2}\right], \\
\frac{\partial}{\partial t} \left[ \frac{\mathcal{E}}{f} \right]
& = i \mathcal{E} \left[ \frac{\lambda^2 - \mu_0^2}{f}
- \frac{k c }{f^2}\right], \\
\frac{\partial}{\partial x}  \left[ \frac{\mathcal{E}}{f'} \right]
& = i \mathcal{E} \left[ \frac{\lambda - \mu_0}{f'} + \frac{k }{(f')^2} \right],\\
\frac{\partial}{\partial x}  \left[ \frac{\mathcal{E}}{f} \right]
& = i \mathcal{E} \left[ \frac{\lambda - \mu_0}{f} + \frac{k }{f^2} \right],\\
\frac{\partial^2}{\partial x^2}  \left[ \frac{\mathcal{E}}{f'} \right]
& = - \mathcal{E} \left[\frac{(\lambda - \mu_0)^2}{f'} + \frac{k (2 \lambda - 2 \mu_0 - k)}{(f')^2}+ 
\frac{2 k^2}{(f')^3}\right],\\
\frac{\partial^2}{\partial x^2}  \left[ \frac{\mathcal{E}}{f} \right]
& = - \mathcal{E} \left[\frac{(\lambda - \mu_0)^2}{f} + \frac{k (2 \lambda - 2 \mu_0 - k)}{f^2}+ 
\frac{2 k^2}{f^3}\right]. 
\end{align*}
Substituting these formulas into the adjoint linearized equation \eqref{BO-adjoint} we obtain that 
\begin{align*}
&\psi_t + 2 u \psi_x + H(\psi_{xx}) \\
& = - \frac{2 i \mathcal{E}}{f^2 (f')^2} \left[k (1 -f) + \lambda_0 (f' - f) \right] \left[
\lambda^2 f f' - \lambda f ( \mu_0 f' - k) + k \lambda_0 f' \right].
\end{align*}
Since $k (1 -f) + \lambda_0 (f' - f) = (1-f) (k + \lambda_0 - \lambda_0 e^{-2\phi}) = 0$, it follows that $\psi_{III}$ satisfies the adjoint linearized equation \eqref{BO-adjoint}. 
\end{proof}

If $\psi$ is a solution of the adjoint linearized BO equation (\ref{BO-adjoint}), then $q = \psi_x$ is a solution of  the linearized BO equation (\ref{BO-lin}). From the adjoint eigenfunctions of Proposition \ref{prop-eigenfunctions}, we define four eigenfunctions by using the following normalizations:
\begin{align}
q_{I}(\lambda;x,t) = -\frac{i}{\lambda} \partial_x \psi_I(\lambda;x,t), &\quad 
q_{II}(\lambda;x,t) = \frac{i}{\lambda} \partial_x \psi_{II}(\lambda;x,t), 
\quad \lambda \in (0,\infty),
\label{solution-adjoint-12} 
\end{align}
and
\begin{equation}
\label{solution-adjoint-34}
q_{III}(\lambda;x,t) = -\frac{i}{\lambda} \partial_x \psi_{III}(\lambda;x,t), \quad 
q_{IV}(\lambda;x,t) = \frac{i}{\lambda} \partial_x \psi_{IV}(\lambda;x,t), \quad \lambda \in (\lambda_0,\mu_0).
\end{equation}
It follows from (\ref{conjugation-1}) and (\ref{conjugation-2}) that the eigenfunctions in (\ref{solution-adjoint-12}) and (\ref{solution-adjoint-34}) are related by the same complex conjugate symmetry:
	\begin{equation}
\label{conjugation-1-q}
q_{II} = e^{-2 \phi} \frac{\lambda - \mu_0}{\lambda - \lambda_0} \bar{q}_I
\end{equation}
and
	\begin{equation}
\label{conjugation-2-q}
q_{IV} = e^{-\phi} \frac{\lambda - \mu_0}{\lambda - \lambda_0} \bar{q}_{III}.
\end{equation}
The following proposition gives a simpler form of the eigenfunctions (\ref{solution-adjoint-12}) and  (\ref{solution-adjoint-34}). 

\begin{proposition}
	\label{prop-eigenfunctions-adjoint}
Eigenfunctions (\ref{solution-adjoint-12}) and  (\ref{solution-adjoint-34}) can be rewritten in the explicit form:
		\begin{align}
		q_{I}(\lambda;x,t) &= e^{i \lambda x + i \lambda^2 t} \frac{f(x-ct)}{f'(x-ct)} \left[ 1 + \frac{k}{(\lambda - \mu_0) f'(x-ct)} \right], \label{solution-1-adjont} \\
		q_{II}(\lambda;x,t) &= e^{-i \lambda x - i \lambda^2 t} \frac{f'(x-ct)}{f(x-ct)} \left[ 1 - \frac{k}{(\lambda - \lambda_0) f(x-ct)} \right], \label{solution-2-adjoint}
		\end{align}
		and
	\begin{align}
q_{III}(\lambda;x,t) &= e^{i (\lambda-\lambda_0) x + i (\lambda^2 - \lambda_0^2) t} \frac{1}{f'(x-ct)} \left[ 1 + \frac{k}{(\lambda - \mu_0) f'(x-ct)} \right] \notag \\
& \qquad -\frac{\lambda_0}{\lambda} e^{i (\lambda-\lambda_0) x + i (\lambda^2 - \lambda_0^2) t} \frac{1}{f(x-ct)} \left[ 1 + \frac{k}{(\lambda - \mu_0) f(x-ct)} \right],  \label{solution-3-adjont} \\
q_{IV}(\lambda;x,t) &= e^{-i (\lambda - \mu_0) x - i (\lambda^2 - \mu_0^2) t} \frac{1}{f(x-ct)} \left[ 1 - \frac{k}{(\lambda - \lambda_0) f(x-ct)} \right] \notag \\
& \qquad - \frac{\mu_0}{\lambda} e^{-i (\lambda - \mu_0) x - i (\lambda^2 - \mu_0^2) t} \frac{1}{f'(x-ct)} \left[ 1 - \frac{k}{(\lambda - \lambda_0) f'(x-ct)} \right]. \label{solution-4-adjoint}
\end{align}		
\end{proposition}

\begin{proof}
It follows from \eqref{lambda-mu} and \eqref{Hirota} that
	\begin{equation}\label{j-1}
\left(\frac{1}{f}\right)_x = \frac{ik (1-f)}{f^2},\quad  \left(\frac{1}{f'}\right)_x = \frac{ik (1-f')}{(f')^2}, \quad 
\frac{1}{f f'} = \frac{1}{k}\left(\frac{\mu_0}{f'} - \frac{\lambda_0}{f}\right).
	\end{equation}
By using the product and quotient rules, we obtain 
from (\ref{solution-1}), (\ref{solution-2}), and (\ref{solution-adjoint-12}):
	\begin{align*}
	q_I &=  e^{i \lambda x + i \lambda^2 t} \left[ 
	\frac{f}{f'} + \frac{k}{(\lambda - \mu_0) f'} + \frac{k}{\lambda} \left( \frac{f - f'}{(f')^2} - \frac{k(f'-1)}{(\lambda - \mu_0) (f')^2} \right) 
	\right] \\
	&=  e^{i \lambda x + i \lambda^2 t} \left[ 
	\frac{f}{f'} + 
	\frac{k f}{(\lambda - \mu_0) (f')^2} + 
	\frac{k}{\lambda (\lambda - \mu_0) (f')^2} \left( - \mu_0 f + \mu_0 f' - k f' + k \right) 
	\right] \\
	&=  e^{i \lambda x + i \lambda^2 t} \left[ 
	\frac{f}{f'} + 
	\frac{k f}{(\lambda - \mu_0) (f')^2} 	\right]
	\end{align*}
	and
	\begin{align*}
q_{II} &=  e^{-i \lambda x - i \lambda^2 t} \left[ 
\frac{f'}{f} - \frac{k}{(\lambda - \lambda_0) f} - \frac{k}{\lambda} \left( \frac{f' - f}{f^2} + \frac{k(f-1)}{(\lambda - \lambda_0) f^2} \right) 
\right] \\
&=  e^{-i \lambda x- i \lambda^2 t} \left[ 
\frac{f'}{f} - 
\frac{k f'}{(\lambda - \lambda_0) f^2} -
\frac{k}{\lambda (\lambda - \lambda_0) f^2} \left( - \lambda_0 f' + \lambda_0 f + k f - k \right) 
\right] \\
&=  e^{-i \lambda x - i \lambda^2 t} \left[ 
\frac{f'}{f} -
\frac{k f'}{(\lambda - \lambda_0) f^2} 	\right],
\end{align*}
where we have used (\ref{lambda-mu}) and (\ref{Hirota}) for simplifications of the resulting expressions.
After factorization, these formulas are written in the form (\ref{solution-1-adjont}) and (\ref{solution-2-adjoint}).	

Similarly, by using the product and quotient rules as well as 
(\ref{j-1}), we obtain 
from (\ref{solution-3}), (\ref{solution-4}), and (\ref{solution-adjoint-34}):
\begin{align*}
q_{III} &= - \frac{i}{\lambda (\lambda - \mu_0)} \partial_x \left[ e^{i (\lambda - \lambda_0) x + i (\lambda^2 - \lambda^2_0) t}  \left( \frac{\lambda}{f'} - \frac{\lambda_0}{f} \right) \right] \\
&= \frac{e^{i (\lambda - \lambda_0) x + i (\lambda^2 - \lambda^2_0) t} }{\lambda (\lambda - \mu_0)}  \left[ {(\lambda - \lambda_0)}\left( \frac{\lambda}{ f'} - \frac{\lambda_0}{f} \right) + {k \lambda}\left( \frac{1}{(f')^2} - \frac{1}{f'} \right)
- {k \lambda_0}\left( \frac{1}{f^2} - \frac{1}{f} \right)\right]\\
& = e^{i (\lambda - \lambda_0) x + i (\lambda^2 - \lambda^2_0) t} \left[
\frac{1}{f'} \left(1 + \frac{k}{(\lambda - \mu_0) f'}\right) - \frac{\lambda_0}{\lambda f } \left(1 + \frac{k}{(\lambda - \mu_0) f} \right) \right]
\end{align*}
and
\begin{align*}
q_{IV} & = \frac{i}{\lambda (\lambda - \lambda_0)} \partial_x \left[ e^{- i (\lambda - \mu_0) x - i (\lambda^2 - \mu^2_0) t}  \left( \frac{\lambda}{f} - \frac{\mu_0}{f'} \right) \right] \\
&= \frac{e^{- i (\lambda - \mu_0) x - i (\lambda^2 - \mu^2_0) t}}{\lambda (\lambda - \lambda_0)} \left[
(\lambda - \mu_0) \left(\frac{\lambda}{f} - \frac{\mu_0}{f'}\right) - k \lambda \left(\frac{1}{f^2} - \frac{1}{f}\right) + k \mu_0
\left(\frac{1}{(f')^2} - \frac{1}{f'}\right) \right] \\
& = e^{- i (\lambda - \mu_0) x - i (\lambda^2 - \mu^2_0) t } \left[\frac{1}{f} \left(1 - \frac{k}{(\lambda - \lambda_0) f}\right)
- \frac{\mu_0}{\lambda f'} \left(1 - \frac{k}{(\lambda - \lambda_0) f'} \right) \right].
\end{align*}
These expressions coincide with those in \eqref{solution-3-adjont} and \eqref{solution-4-adjoint}. 
\end{proof}

\section{Orthogonality of eigenfunctions in $L^2(\mathbb{R})$} 
\label{sec-3}

The eigenfunctions $q_{I}$, $q_{II}$, $q_{III}$, and $q_{IV}$ in Proposition \ref{prop-eigenfunctions-adjoint} are particular solutions of the linearized BO equation (\ref{BO-lin}), whereas the eigenfunctions $\psi_{I}$, $\psi_{II}$, $\psi_{III}$, and $\psi_{IV}$ in Proposition \ref{prop-eigenfunctions} are particular solutions of the adjoint linearized BO equation (\ref{BO-adjoint}). To simplify notations, we do not write arguments $(x,t)$ in the eigenfunctions  $q_{I,II,III,IV}$ and $\psi_{I,II,III,IV}$.

Here, we use the standard inner product
$$
f,g \in L^2(\mathbb{R}) : \quad \langle f, g \rangle_{L^2(\R)} = \int_{\mathbb{R}} f(x) \bar{g}(x) dx
$$
to establish the orthogonality of eigenfunctions in $L^2(\R)$. The 
orthogonality relations for the bounded and oscillatory eigenfunctions, which are non-decaying to zero at infinity, are understood in the sense of generalized distributional identities (\ref{distrib-1}) and (\ref{distrib-2}) for the Fourier modes. 

The following three propositions establish orthogonality between all the eigenfunctions. 

\begin{proposition}
	\label{prop-orthog-1}
	The set $\{ \bar{\psi}_I, \bar{\psi}_{II}\}$ is orthogonal to the set $\{ q_I, q_{II} \}$ in the sense:
	\begin{equation}
	\label{orth-1}
\begin{cases}
\langle q_I(\lambda), \bar{\psi}_I(\lambda') \rangle_{L^2(\R)} =  
\langle q_{II}(\lambda), \bar{\psi}_{II}(\lambda') \rangle_{L^2(\R)} = 0, \\	
\langle q_I(\lambda), \bar{\psi}_{II}(\lambda') \rangle_{L^2(\R)} =  
\langle q_{II}(\lambda), \bar{\psi}_{I}(\lambda') \rangle_{L^2(\R)} = 2\pi \delta(\lambda - \lambda'), 
\end{cases}
	\end{equation}
where $\lambda, \lambda' \in (0,\infty)$ and $\delta(\lambda - \lambda')$ is the Dirac delta-distribution. 
\end{proposition}

\begin{proof}
To prove the first line in (\ref{orth-1}), we use (\ref{solution-1}) and (\ref{solution-1-adjont}) to write 
\begin{equation}
\label{representation-1}
q_{I}(\lambda) \psi_{I}(\lambda') = e^{i (\lambda + \lambda') x + i (\lambda^2 + (\lambda')^2) t} \frac{f}{f'} \left[ 1 + \frac{k}{(\lambda - \mu_0) f'} \right] 
\left[ \frac{f}{f'} + \frac{k}{(\lambda' - \mu_0) f'} \right].
\end{equation}
Since the function $f'$ has no zeros in $\mathbb{C}^+$ of the complex $x$-plane, the function 
$q_I(\lambda) \psi_I(\lambda')$ is analytic in $\mathbb{C}^+$. In addition, the exponential factor $e^{i(\lambda + \lambda') x}$ decays to $0$ exponentially fast as ${\rm Im}(x) \to +\infty$ if $\lambda, \lambda' \in (0,\infty)$, whereas the multiplicative factor is bounded. Integrating  $q_{I}(\lambda) \psi_{I}(\lambda')$ over the closed contour
$$
[-R,R] \cup \{ R e^{i \theta}, \; 0 \leq \theta \leq \pi \}, \quad R > 0,
$$
yields the zero integral by Cauchy's integral theorem for every $R > 0$. Therefore, 
\begin{align*}
\langle q_I(\lambda), \bar{\psi}_I(\lambda') \rangle_{L^2(\R)} &= \lim_{R\to \infty} 
\int_{-R}^R  q_I(\lambda;x,t) \psi_I(\lambda';x,t) dx \\
&= - \lim_{R\to \infty} \int_{0}^{\pi}  q_I(\lambda;R e^{i \theta},t) \psi_I(\lambda';R e^{i \theta},t) i R e^{i \theta} d \theta = 0,
\end{align*}	
where the last equality is obtained by Lebesgue's dominated convergence theorem due
to the exponential convergence of $e^{i(\lambda + \lambda') x} \to 0$ for $x = R e^{i \theta}$, $0 < \theta < \pi$ as $R \to \infty$. This yields $\langle q_I(\lambda), \bar{\psi}_I(\lambda') \rangle_{L^2(\R)} = 0$ in (\ref{orth-1}). The proof for $\langle q_{II}(\lambda), \bar{\psi}_{II}(\lambda') \rangle_{L^2(\R)} = 0$ is similar but the contour of integration is closed in $\mathbb{C}^-$.

To prove the second line in (\ref{orth-1}), we use (\ref{solution-2}) and (\ref{solution-1-adjont}) to write
\begin{align}
\label{representation-2}
q_I(\lambda) \psi_{II}(\lambda') &= e^{i (\lambda - \lambda')x + i (\lambda^2 - (\lambda')^2)t} \left[  1 + \frac{k}{(\lambda - \mu_0) f'} \right] \left[ 1 - \frac{k}{(\lambda' - \lambda_0) f'} \right].
\end{align}
Integrating $e^{i (\lambda - \lambda')x + i (\lambda^2 - (\lambda')^2)t}$ in $x$ over $\R$ yields $2\pi \delta(\lambda - \lambda')$ for $\lambda, \lambda' \in (0,\infty)$ due to (\ref{distrib-1}). Hence, we get $\langle q_I(\lambda), \bar{\psi}_{II}(\lambda') \rangle_{L^2(\R)} = 2\pi \delta(\lambda - \lambda')$ if and only if the integral of all remaining terms after expanding the product vanishes:
\begin{align*}
0 &= \frac{k}{(\lambda - \mu_0) (\lambda' - \lambda_0)} \int_{\R} e^{i (\lambda - \lambda')x + i (\lambda^2 - (\lambda')^2)t}  \left( \frac{(\lambda' - \lambda)}{f'}  + \frac{k}{f'} - \frac{k}{(f')^2}\right) dx \\
&= \frac{ik}{(\lambda - \mu_0) (\lambda' - \lambda_0)} \int_{\R} \frac{d}{dx} \left( \frac{e^{i (\lambda - \lambda')x + i (\lambda^2 - (\lambda')^2)t} }{f'} \right) dx \\
&= \frac{ik}{(\lambda - \mu_0) (\lambda' - \lambda_0)} \int_{\R} \sum_{n=0}^{\infty} (-1)^n e^{-n \phi} \frac{d}{dx} 
e^{i (\lambda - \lambda' + nk) x + i (\lambda^2 - (\lambda')^2 - nck) t} dx 
\end{align*}
which is true for $\lambda, \lambda' \in (0,\infty)$ due to (\ref{distrib-2}) and the absolute convergence of the geometric series for $\phi > 0$. The result for 
$$
\langle q_{II}(\lambda), \bar{\psi}_{I}(\lambda') \rangle_{L^2(\R)} 
= \overline{\langle q_{I}(\lambda), \bar{\psi}_{II}(\lambda') \rangle_{L^2(\R)}} = 2\pi \delta(\lambda - \lambda')
$$ 
holds by the complex conjugate symmetry (\ref{conjugation-1}) and (\ref{conjugation-1-q}).
\end{proof}

\begin{proposition} 
		\label{prop-orthog-2}
The set $\{\bar{\psi}_I,\bar{\psi}_{II}\}$ is orthogonal to $\{ q_{III},q_{IV}\}$ and the set $\{\bar{\psi}_{III},\bar{\psi}_{IV}\}$ is orthogonal to $\{ q_{I},q_{II}\}$ in the sense:
\begin{equation}
\label{orth-3-b}
\begin{cases}
\langle q_{I}(\lambda), \bar{\psi}_{III}(\lambda') \rangle_{L^2(\R)} =  
\langle q_{II}(\lambda), \bar{\psi}_{III}(\lambda') \rangle_{L^2(\R)} = 0, \\	
\langle q_{I}(\lambda), \bar{\psi}_{IV}(\lambda') \rangle_{L^2(\R)} =  
\langle q_{II}(\lambda), \bar{\psi}_{IV}(\lambda') \rangle_{L^2(\R)} = 0, 
\end{cases} 
\end{equation} 
and 
\begin{equation}
\label{orth-3-c}
\begin{cases}
\langle q_{III}(\lambda'), \bar{\psi}_{I}(\lambda) \rangle_{L^2(\R)} =  
\langle q_{IV}(\lambda'), \bar{\psi}_{I}(\lambda) \rangle_{L^2(\R)} = 0, \\	
\langle q_{III}(\lambda'), \bar{\psi}_{II}(\lambda) \rangle_{L^2(\R)} =  
\langle q_{IV}(\lambda'), \bar{\psi}_{II}(\lambda) \rangle_{L^2(\R)} = 0, 
\end{cases}
\end{equation} 
where $\lambda \in (0,\infty)$ and $\lambda' \in (\lambda_0,\mu_0)$.
\end{proposition}

\begin{proof}
To simplify the presentation, we set $t = 0$ in all functions. 
	
It follows from (\ref{solution-3}) and (\ref{solution-1-adjont}) that 
\begin{equation*}
q_{I}(\lambda) \psi_{III}(\lambda') = e^{i (\lambda + \lambda' - \lambda_0) x} \frac{1}{f'} \left[ 1 + \frac{k}{(\lambda - \mu_0) f'} \right] 
\left[ \frac{f}{f'} + \frac{k}{(\lambda' - \mu_0) f'} \right].
\end{equation*}
Since the function $f'$ has no zeros in $\mathbb{C}^+$ of the complex $x$-plane, the function $q_{I}(\lambda) \psi_{III}(\lambda')$ is analytic in $\mathbb{C}^+$. 
In addition, $\lambda + \lambda' - \lambda_0 > 0$ if $\lambda \in (0,\infty)$ and $\lambda' \in (\lambda_0,\mu_0)$ so that the exponential factor $e^{i (\lambda + \lambda' - \lambda_0) x}$ decays to zero as ${\rm Im}(x) \to +\infty$, whereas the multiplicative factor is bounded. By using Cauchy's integral theorem as in the proof of Proposition \ref{prop-orthog-1}, we have $\langle q_{I}(\lambda), \bar{\psi}_{III}(\lambda') \rangle_{L^2(\R)} = 0$. 

Computations for $\langle q_{I}(\lambda), \bar{\psi}_{IV}(\lambda') \rangle_{L^2(\R)} =  0$ are similar since 
\begin{equation*}
q_{I}(\lambda) \psi_{IV}(\lambda') = e^{i (\lambda + \mu_0 - \lambda') x} \frac{1}{f'} \left[ 1 + \frac{k}{(\lambda - \mu_0) f'} \right] 
\left[ 1 - \frac{k}{(\lambda' - \lambda_0) f'} \right].
\end{equation*}
and $\lambda + \mu_0 - \lambda' > 0$ if $\lambda \in (0,\infty)$ and $\lambda' \in (\lambda_0,\mu_0)$.

It follows from (\ref{solution-3}) and (\ref{solution-2-adjoint}) that 
\begin{align*}
q_{II}(\lambda) \psi_{III}(\lambda')  = e^{- i (\lambda - \lambda' + \lambda_0) x}
\frac{1}{f} \left[1 - \frac{k}{(\lambda - \lambda_0)f}\right] \left[1 + \frac{k}{ (\lambda' - \mu_0) f}\right].
\end{align*}
Since the function $f$ has no zeros in $\mathbb{C}^-$ of the complex $x$-plane, the function $q_{II}(\lambda) \psi_{III}(\lambda')$ is analytic in $\mathbb{C}^-$. If $\lambda > \lambda' - \lambda_0$, the exponential factor $e^{- i (\lambda - \lambda' + \lambda_0) x}$ decays to zero as ${\rm Im}(x) \to -\infty$, and Cauchy's integral theorem implies that $\langle q_{II}(\lambda), \bar{\psi}_{III}(\lambda') \rangle_{L^2(\R)} = 0$. However, this argument is not applicable for the case $\lambda < \lambda' - \lambda_0$. Nevertheless, it follows 
from \eqref{Hirota} that 
\begin{equation}
\label{reflected-functions}
	f(-x) = e^{- i k x + \phi} f'(x), \qquad f'(-x) = e^{- i k x - \phi} f(x),
\end{equation}
and 
\begin{align*}
1 - \frac{k}{(\lambda - \lambda_0)f(-x)} &= 1 - \frac{k(f'(x)-1)}{(\lambda - \lambda_0) f'(x)} = \frac{\lambda - \mu_0}{\lambda - \lambda_0} \left[ 1 + \frac{k}{(\lambda - \mu_0) f'(x)} \right], \\
1 + \frac{k}{ (\lambda' - \mu_0) f(-x)} &= 1 + \frac{k(f'(x)-1)}{(\lambda' - \mu_0) f'(x)} = \frac{\lambda' - \lambda_0}{\lambda' - \mu_0} \left[ 1 - \frac{k}{(\lambda' - \lambda_0) f'(x)} \right],
\end{align*}
This allows us to write 
\begin{align*}
q_{II}(\lambda;-x) \psi_{III}(\lambda';-x)  &= \frac{(\lambda - \mu_0) (\lambda' - \lambda_0)}{(\lambda - \lambda_0) (\lambda' - \mu_0)} e^{i (\lambda + \mu_0 - \lambda') x - \phi} \\
& \quad \times 
\frac{1}{f'} \left[1 + \frac{k}{(\lambda - \mu_0) f'}\right] \left[1 - \frac{k}{ (\lambda' - \lambda_0) f'}\right].
\end{align*}
The function $q_{II}(\lambda;-x) \psi_{III}(\lambda';-x)$ is analytic in $\mathbb{C}^+$ and the exponential factor $e^{i (\lambda + \mu_0 - \lambda') x}$ decays to zero as ${\rm Im}(x) \to +\infty$ since $\lambda + \mu_0 - \lambda' > 0$. By using Cauchy's integral theorem, we get 
\begin{equation*} 
\langle q_{II}(\lambda; x), \bar{\psi}_{III}(\lambda'; x) \rangle_{L^2(\R)} = \langle q_{II}(\lambda; - x), \bar{\psi}_{III}(\lambda'; -x) \rangle_{L^2(\R)} = 0.
\end{equation*}

Similarly, we write the representation:
\begin{align*}
q_{II}(\lambda) \psi_{IV}(\lambda') = e^{- i (\lambda + \lambda' - \mu_0) x}
\frac{1}{f} \left[1 - \frac{k}{(\lambda - \lambda_0)f}\right] \left[ \frac{f'}{f} - \frac{k}{ (\lambda' - \lambda_0) f}\right].
\end{align*}
By using (\ref{reflected-functions}) and 
\begin{align*}
\frac{f'(-x)}{f(-x)} - \frac{k}{ (\lambda' - \lambda_0) f(-x)} &= 
e^{-2\phi}  \left[ \frac{f(x)}{f'(x)} - \frac{k (f(x) - 1)}{(\lambda' - \lambda_0) f'(x)} \right] \\
&= e^{-2\phi} \frac{\lambda' - \mu_0}{\lambda' - \lambda_0}  \left[ \frac{f(x)}{f'(x)} + \frac{k}{(\lambda' - \mu_0) f'(x)} \right], 
\end{align*}
we reflect $x$ in the previous expression and rewrite it in the equivalent form:
\begin{align*}
q_{II}(\lambda;-x) \psi_{IV}(\lambda';-x) &= \frac{(\lambda - \mu_0) (\lambda' - \mu_0)}{(\lambda - \lambda_0) (\lambda' - \lambda_0)} e^{i (\lambda + \lambda' - \lambda_0) x - 3 \phi} \\
& \quad \times 
\frac{1}{f'} \left[1 + \frac{k}{(\lambda - \mu_0)f'}\right] \left[ \frac{f}{f'} + \frac{k}{ (\lambda' - \mu_0) f'}\right].
\end{align*}
The function $q_{II}(\lambda;-x) \psi_{IV}(\lambda';-x)$ is analytic in $\mathbb{C}^+$ and the exponential factor $e^{i (\lambda + \lambda' - \lambda_0) x}$ decays to zero as ${\rm Im}(x) \to +\infty$ since $\lambda + \lambda' - \lambda_0 > 0$. By using Cauchy's integral theorem, we get 
\begin{equation*} 
\langle q_{II}(\lambda; x), \bar{\psi}_{IV}(\lambda'; x) \rangle_{L^2(\R)} = \langle q_{II}(\lambda; - x), \bar{\psi}_{IV}(\lambda'; -x) \rangle_{L^2(\R)} = 0.
\end{equation*}
Thus, the orthogonality relations \eqref{orth-3-b} have been proven. 

The orthogonality relations (\ref{orth-3-c}) are equivalent to (\ref{orth-3-b}) after integration by parts with the help of (\ref{distrib-2}). Hence, (\ref{orth-3-b}) implies (\ref{orth-3-c}).
\end{proof}

\begin{proposition}
		\label{prop-orthog-3}
	The set $\{ \bar{\psi}_{III}, \bar{\psi}_{IV}\}$ is orthogonal to the set $\{ q_{III},q_{IV}\}$ in the sense:
	\begin{equation}
	\label{orth-3-a}
\begin{cases}
	\langle q_{III}(\lambda), \bar{\psi}_{III}(\lambda') \rangle_{L^2(\R)} =  
	\langle q_{IV}(\lambda), \bar{\psi}_{IV}(\lambda') \rangle_{L^2(\R)} = 0, \\	
	\langle q_{III}(\lambda), \bar{\psi}_{IV}(\lambda') \rangle_{L^2(\R)} =  
	\langle q_{IV}(\lambda), \bar{\psi}_{III}(\lambda') \rangle_{L^2(\R)} = \frac{2\pi}{(e^{\phi} - e^{-\phi})} \delta(\lambda - \lambda'), 
\end{cases}
	\end{equation}
	where $\lambda, \lambda' \in (\lambda_0,\mu_0)$ and $\delta(\lambda - \lambda')$ is the Dirac delta-distribution. 
\end{proposition}

\begin{proof}
	To prove the second line in (\ref{orth-3-a}), we use \eqref{solution-4} and (\ref{solution-3-adjont}) to write
\begin{align*}
q_{III}(\lambda) \psi_{IV}(\lambda') &= e^{i (\lambda - \lambda' + k) x + i  (\lambda^2 - (\lambda')^2 - \lambda_0^2 + \mu_0^2) t} \left[1 - \frac{k}{(\lambda' - \lambda_0) f'} \right] \\
& \quad \times \frac{1}{f}   \left( \frac{1}{f'} \left[ 1 + \frac{k}{(\lambda - \mu_0) f'} \right] 
- \frac{\lambda_0}{\lambda f} \left[ 1 + \frac{k}{(\lambda - \mu_0) f} \right] \right).
\end{align*}
By using 
\begin{equation}
\label{c-3.14}
1 + \frac{k}{(\lambda - \mu_0) f} = \frac{1}{\lambda - \mu_0}
\left[ \lambda - \frac{\lambda_0 f'}{f} \right], \quad 
1 - \frac{k}{(\lambda' - \lambda_0) f'} = \frac{1}{\lambda' - \lambda_0}
\left[ \lambda' - \frac{\mu_0 f}{f'} \right],
\end{equation}
we rewrite the factor behind the exponential function as follows: 
\begin{align*}
& \frac{1}{f}   \left( \frac{1}{f'} \left[ 1 + \frac{k}{(\lambda - \mu_0) f'} \right] 
- \frac{\lambda_0}{\lambda f} \left[ 1 + \frac{k}{(\lambda - \mu_0) f} \right] \right)
\left[1 - \frac{k}{(\lambda' - \lambda_0) f'} \right]\\
&= -\frac{\mu_0}{(\lambda' - \lambda_0) (f')^2} \left[ 1 + \frac{k}{(\lambda - \mu_0) f'} \right]  - \frac{\lambda' \lambda_0}{\lambda (\lambda - \mu_0) (\lambda' - \lambda_0) f^2} \left[ \lambda - \frac{\lambda_0 f'}{f} \right]  \\
& \quad + \frac{\lambda'}{(\lambda' - \lambda_0) f f'} \left[ 1 + \frac{k}{(\lambda - \mu_0) f'} \right] + \frac{\mu_0 \lambda_0}{\lambda ( \lambda - \mu_0) (\lambda' - \lambda_0) f f'} \left[ \lambda - \frac{\lambda_0 f'}{f} \right] \\
&= -\frac{\mu_0}{(\lambda - \mu_0) (\lambda' - \lambda_0) (f')^2} \left[ \lambda - \lambda' -  \mu_0 + \frac{k}{f'} \right] \\
& \quad - \frac{\lambda_0}{\lambda (\lambda - \mu_0) (\lambda' - \lambda_0) f^2} \left[ \lambda \lambda'  + \lambda_0 \mu_0 - \frac{\lambda' \lambda_0 f'}{f} \right]  + \frac{\lambda' (\lambda  - \lambda_0) + \mu_0 (\lambda_0 - \lambda')}{(\lambda - \mu_0) (\lambda' - \lambda_0) f f'}.
\end{align*}
Since $\lambda - \lambda' + k = (\lambda - \lambda_0) + (\mu_0 - \lambda') > 0$ due to $\lambda, \lambda' \in (\lambda_0, \mu_0)$, the function 
$$
-\frac{\mu_0}{(\lambda - \mu_0) (\lambda' - \lambda_0) (f')^2} \left[ \lambda - \lambda' -  \mu_0 + \frac{k}{f'} \right] e^{i (\lambda - \lambda' + k) x + i  (\lambda^2 - (\lambda')^2 - \lambda_0^2 + \mu_0^2) t} 
$$
is analytic in $\mathbb{C}^+$ and the exponential factor 
$e^{i (\lambda - \lambda' + k) x}$ decays to zero as ${\rm Im}(x) \to +\infty$. By using Cauchy's integral theorem in $\mathbb{C}^+$, the corresponding integral vanishes. To show that the next term also vanishes, we set $t = 0$ and rewrite it as 
$$
- \frac{\lambda_0}{\lambda (\lambda - \mu_0) (\lambda' - \lambda_0) f^2} \left[ \lambda \lambda'  + \lambda_0 \mu_0 - \frac{\lambda' \lambda_0 f'}{f} \right]  e^{i (\lambda - \lambda' + k) x}.
$$
Reflecting $x$ and using (\ref{reflected-functions}), we get 
$$
- \frac{\lambda_0}{\lambda (\lambda - \mu_0) (\lambda' - \lambda_0) (f')^2} \left[ \lambda \lambda'  + \lambda_0 \mu_0 - \frac{\lambda' \mu_0  f}{f'} \right] 
e^{-i (\lambda - \lambda' - k) x  - 2\phi}.
$$
Since $-\lambda + \lambda' + k = (\lambda' - \lambda_0) + (\mu_0 - \lambda) > 0$, the exponential factor $e^{-i (\lambda - \lambda' - k) x}$ decays to zero as ${\rm Im}(x) \to +\infty$, whereas the reflected function is analytic in $\mathbb{C}^+$. The corresponding integral is zero again for $t = 0$ and, in fact, for every $t \in \mathbb{R}$. As a result, 
we reduce the inner product between $\bar{\psi}_{IV}(\lambda')$ and $q_{III}(\lambda)$ to the last term of the expansion:
\begin{align*}
\langle q_{III}(\lambda), \bar{\psi}_{IV}(\lambda') \rangle_{L^2(\R)} &= 
\frac{\lambda' (\lambda  - \lambda_0) + \mu_0 (\lambda_0 - \lambda')}{(\lambda - \mu_0) (\lambda' - \lambda_0)} \int_{\R} \frac{e^{i (\lambda - \lambda' + k) x + i  (\lambda^2 - (\lambda')^2 - \lambda_0^2 + \mu_0^2) t}}{f f'} dx.
\end{align*}
By using (\ref{per-wave}) and (\ref{Hirota}), we get 
\begin{equation}
\label{c-3.15}
\frac{e^{i k \xi}}{f f'} = \frac{1}{k (e^{\phi} - e^{-\phi})} U = \frac{i}{k (e^{\phi} - e^{-\phi})} \left(\frac{f'_x}{f'} - \frac{f_x}{f} \right).
\end{equation}
Therefore, we can rewrite the remaining integral in the form:
\begin{align*}
& \quad \langle q_{III}(\lambda), \bar{\psi}_{IV}(\lambda') \rangle_{L^2(\R)} \\
&= 
\frac{\lambda' (\lambda  - \lambda_0) + \mu_0 (\lambda_0 - \lambda')}{(\lambda - \mu_0) (\lambda' - \lambda_0)}  \frac{i}{k (e^{\phi} - e^{-\phi})} \int_{\R} \left(\frac{f'_x}{f'} - \frac{f_x}{f} \right)  e^{i (\lambda - \lambda') x + i  (\lambda^2 - (\lambda')^2) t} dx,
\end{align*}
where we have used the identity $kc = \lambda_0^2 - \mu_0^2$, which follows from (\ref{Lax-spectrum-continuous}). By the same argument with $\lambda - \lambda' + k > 0$, we get 
$$
\int_{\R} \frac{f'_x}{f'} e^{i (\lambda - \lambda') x  + i  (\lambda^2 - (\lambda')^2) t} dx = i k e^{-\phi} \int_{\R} \frac{1}{f'} e^{i (\lambda - \lambda' + k) x  + i  (\lambda^2 - (\lambda')^2 - \lambda_0^2 + \mu_0^2) t} dx = 0,
$$
after the path is closed in $\mathbb{C}^+$. On the other hand, we have 
\begin{equation}
\label{c-3.16}
\frac{f_x}{f} = \frac{ik e^{ik\xi + \phi}}{1 + e^{ik\xi + \phi}} = \frac{ik}{1 + e^{-ik\xi - \phi}} = ik \sum_{n=0}^{\infty} (-1)^n e^{-i k n \xi - n\phi},
\end{equation}
where the geometric series converges since $\phi > 0$. Since 
$$
\lambda - \lambda' - nk = (\lambda - \mu_0) + (\lambda_0 - \lambda') - (n-1)k < 0
$$ 
for every integer $n \geq 1$, we get 
$$
\int_{\R} \sum_{n=1}^{\infty} (-1)^n e^{-i k n \xi - n\phi} e^{i (\lambda - \lambda') x + i  (\lambda^2 - (\lambda')^2) t} dx = 0,
$$
after the path is closed in $\mathbb{C}^-$. This leads us to the only remaining, potentially nonzero term:
\begin{align*}
\langle q_{III}(\lambda), \bar{\psi}_{IV}(\lambda') \rangle_{L^2(\R)} &= 
\frac{\lambda' (\lambda  - \lambda_0) + \mu_0 (\lambda_0 - \lambda')}{(\lambda - \mu_0) (\lambda' - \lambda_0) (e^{\phi} - e^{-\phi})} \int_{\R} e^{i (\lambda - \lambda') x + i  (\lambda^2 - (\lambda')^2) t} dx \\
&= 2\pi \delta(\lambda - \lambda') \frac{\lambda' (\lambda  - \lambda_0) + \mu_0 (\lambda_0 - \lambda')}{(\lambda - \mu_0) (\lambda' - \lambda_0) (e^{\phi} - e^{-\phi})} \biggr|_{\lambda = \lambda'} \\
&= \frac{2\pi}{e^{\phi} - e^{-\phi}} \delta(\lambda - \lambda'),
\end{align*}
where we have used (\ref{distrib-1}). This completes the proof of the second line in (\ref{orth-3-a}), where the equality 
$$
\langle q_{IV}(\lambda), \bar{\psi}_{III}(\lambda') \rangle_{L^2(\R)} = \overline{\langle q_{III}(\lambda), \bar{\psi}_{IV}(\lambda') \rangle_{L^2(\R)}} = \frac{2\pi}{e^{\phi} - e^{-\phi}} \delta(\lambda - \lambda')
$$ 
follows from the complex-conjugate symmetry (\ref{conjugation-2}) and (\ref{conjugation-2-q}).

To prove the first line in (\ref{orth-3-a}), we use \eqref{solution-3} and (\ref{solution-3-adjont}) to write
\begin{align*}
q_{III}(\lambda) \psi_{III}(\lambda') &= e^{i (\lambda + \lambda' - 2 \lambda_0) x + i  (\lambda^2 + (\lambda')^2 - 2\lambda_0^2) t} \left[ 1 + \frac{k}{(\lambda' - \mu_0) f} \right] \\
& \quad \times \frac{1}{f'}   \left( \frac{1}{f'} \left[ 1 + \frac{k}{(\lambda - \mu_0) f'} \right] 
- \frac{\lambda_0}{\lambda f} \left[ 1 + \frac{k}{(\lambda - \mu_0) f} \right] \right).
\end{align*}
By using (\ref{c-3.14}), we rewrite the factor behind the exponential function as follows:
\begin{align*}
&  \quad \frac{1}{f'}   \left( \frac{1}{f'} \left[ 1 + \frac{k}{(\lambda - \mu_0) f'} \right] 
- \frac{\lambda_0}{\lambda f} \left[ 1 + \frac{k}{(\lambda - \mu_0) f} \right] \right)
\left[ 1 + \frac{k}{(\lambda' - \mu_0) f} \right]  \\
&= \frac{\lambda'}{(\lambda' - \mu_0) (f')^2} \left[ 1 + \frac{k}{(\lambda - \mu_0) f'} \right]  + \frac{\lambda_0^2}{\lambda (\lambda' - \mu_0) f^2} \left[ 1 + \frac{k}{(\lambda - \mu_0) f} \right]  \\
& \quad -\frac{\lambda_0 \lambda'}{\lambda (\lambda' - \mu_0) f f'} \left[ 1 + \frac{k}{(\lambda - \mu_0) f} \right] - \frac{\lambda_0}{( \lambda' - \mu_0) f f'} \left[ 1 + \frac{k}{(\lambda - \mu_0) f'} \right] \\
&= \frac{1}{(\lambda' - \mu_0) (f')^2} \left[ \lambda' - \frac{\lambda_0 \mu_0}{\lambda - \mu_0}  + \frac{\lambda' k}{(\lambda - \mu_0) f'} \right]  \\
& \quad + \frac{\lambda_0^2}{\lambda (\lambda' - \mu_0) (\lambda - \mu_0) f^2} \left[ \lambda + \lambda' - \frac{\lambda_0 f'}{f} \right]  
- \frac{\lambda_0 (\lambda + \lambda' - \lambda_0 - \mu_0)}{(\lambda' - \mu_0) (\lambda - \mu_0) f f'}.
\end{align*}
Since $\lambda + \lambda' - 2 \lambda_0 > 0$ due to $\lambda, \lambda' \in (\lambda_0, \mu_0)$, the function 
$$
\frac{1}{(\lambda' - \mu_0) (f')^2} \left[ \lambda' - \frac{\lambda_0 \mu_0}{\lambda - \mu_0}  + \frac{\lambda' k}{(\lambda - \mu_0) f'} \right] e^{i (\lambda + \lambda' - 2 \lambda_0) x + i  (\lambda^2 + (\lambda')^2 - 2\lambda_0^2) t}
$$
is analytic in $\mathbb{C}^+$ and the exponential factor $e^{i (\lambda + \lambda' - 2 \lambda_0) x}$ decays to zero as ${\rm Im}(x) \to +\infty$. By using Cauchy's integral theorem in $\mathbb{C}^+$, the corresponding integral vanishes. To show that the next term also vanishes, we set $t = 0$ and rewrite it as
$$
\frac{\lambda_0^2}{\lambda (\lambda' - \mu_0) (\lambda - \mu_0) f^2} \left[ \lambda + \lambda' - \frac{\lambda_0 f'}{f} \right]  e^{i (\lambda + \lambda' - 2 \lambda_0) x}.
$$
Reflecting $x$ and using (\ref{reflected-functions}), we get 
$$
\frac{\lambda_0^2}{\lambda (\lambda' - \mu_0) (\lambda - \mu_0) (f')^2} \left[ \lambda + \lambda' - \frac{\mu_0 f}{f'} \right]  
e^{i (2 \mu_0 - \lambda - \lambda') x - 2\phi}.
$$
Since $2 \mu_0 - \lambda - \lambda' > 0$, the exponential factor $e^{i (2 \mu_0 - \lambda - \lambda') x}$ decays to zero as ${\rm Im}(x) \to +\infty$, whereas the reflected function is analytic in $\mathbb{C}^+$. The corresponding integral is zero again for $t = 0$ and, in fact, for every $t \in \mathbb{R}$. As a result, we reduce the inner product between $\bar{\psi}_{III}(\lambda')$ and $q_{III}(\lambda)$ to the last term of the expansion:
\begin{align*}
\langle q_{III}(\lambda), \bar{\psi}_{III}(\lambda') \rangle_{L^2(\R)} &= 
- \frac{\lambda_0 (\lambda + \lambda' - \lambda_0 - \mu_0)}{(\lambda' - \mu_0) (\lambda - \mu_0)} \int_{\R} \frac{e^{i (\lambda + \lambda' - 2 \lambda_0) x + i  (\lambda^2 + (\lambda')^2 - 2 \lambda_0^2) t}}{f f'} dx.
\end{align*}
By using (\ref{c-3.15}) and (\ref{c-3.16}), we can write 
$$
\frac{1}{f f'} = \frac{1}{(e^{\phi} - e^{-\phi})} \left( \frac{e^{\phi}}{f} - \frac{e^{-\phi}}{f'} \right) = 
\frac{1}{(e^{\phi} - e^{-\phi})} \left( e^{-ik \xi} \sum_{n=0}^{\infty} (-1)^n e^{-i k n\xi - n \phi} - \frac{e^{-\phi}}{f'} \right),
$$
where the geometric series converges. Since 
$$
\lambda + \lambda' - 2 \lambda_0 - (n+1) k < 0
$$
for every integer $n \geq 1$, we get 
$$
\int_{\R} \sum_{n=1}^{\infty} (-1)^n e^{-i k n\xi - n \phi}  
 e^{-ik \xi} e^{i (\lambda + \lambda' - 2 \lambda_0) x + i  (\lambda^2 + (\lambda')^2 - 2 \lambda_0^2) t} dx = 0
$$
after the path is closed in $\mathbb{C}^-$. Similarly, we get 
$$
\int_{\R} \frac{e^{-\phi}}{f'} e^{i (\lambda + \lambda' - 2 \lambda_0) x + i  (\lambda^2 + (\lambda')^2 - 2 \lambda_0^2) t} dx = 0
$$
after the path is closed in $\mathbb{C}^+$. This leads us to the only remaining, potentially nonzero term:
\begin{align*}
& \qquad \langle q_{III}(\lambda), \bar{\psi}_{III}(\lambda') \rangle_{L^2(\R)} \\
&= 
- \frac{\lambda_0 (\lambda + \lambda' - \lambda_0 - \mu_0)}{(\lambda' - \mu_0) (\lambda - \mu_0) (e^{\phi} - e^{-\phi})} \int_{\R} e^{i (\lambda + \lambda' - \lambda_0 - \mu_0) x + i k c t + i  (\lambda^2 + (\lambda')^2 - 2 \lambda_0^2) t} dx \\
&= \frac{i \lambda_0}{(\lambda' - \mu_0) (\lambda - \mu_0) (e^{\phi} - e^{-\phi})} \int_{\R} \frac{d}{dx} e^{i (\lambda + \lambda' - \lambda_0 - \mu_0) x + i k c t + i  (\lambda^2 + (\lambda')^2 - 2 \lambda_0^2) t} dx,
\end{align*}
which is zero due to (\ref{distrib-2}).  The first line in (\ref{orth-3-a}) is verified since 
$$
\langle q_{IV}(\lambda), \bar{\psi}_{IV}(\lambda') \rangle_{L^2(\R)} 
= e^{-2\phi} \frac{(\lambda - \mu_0) (\lambda' - \mu_0)}{(\lambda - \lambda_0) (\lambda' - \lambda_0)} 
\overline{\langle q_{III}(\lambda), \bar{\psi}_{III}(\lambda') \rangle_{L^2(\R)}} = 0
$$ 
follows by the complex-conjugate symmetry (\ref{conjugation-2}) and (\ref{conjugation-2-q}).
\end{proof}

\section{Completeness of eigenfunctions in $L^2(\mathbb{R})$}
\label{sec-4}

To simplify the proof of completeness, we set $t = 0$ and suppress the time variable from the list of arguments. We claim that solutions (\ref{solution-1-adjont}), (\ref{solution-2-adjoint}), (\ref{solution-3-adjont}), and (\ref{solution-4-adjoint}) provide a basis in $L^2(\R)$ in the sense that every element $q \in L^2(\R)$ can be uniquely represented by the superposition of elements of the basis:
\begin{align}
q(x) &= \int_{0}^{\infty} \left[ a_{I}(\lambda) q_I(\lambda;x) + a_{II}(\lambda) q_{II}(\lambda;x) \right] d \lambda \notag \\
& \qquad + 
\int_{\lambda_0}^{\mu_0} \left[ a_{III}(\lambda) q_{III}(\lambda;x) + a_{IV}(\lambda) q_{IV}(\lambda;x) \right] d \lambda. \label{basis-2}
\end{align}

If $q \in L^2(\R) \cap L^1(\R)$, then the coefficients in the decomposition (\ref{basis-2}) are computed by using the orthogonality relations of Propositions \ref{prop-orthog-1}, \ref{prop-orthog-2}, and \ref{prop-orthog-3}:
\begin{align}
\label{projection-1}
\begin{cases}
a_{I}(\lambda) = \frac{1}{2\pi} \langle q, \bar{\psi}_{II}(\lambda)  \rangle_{L^2(\R)}, \\
a_{II}(\lambda) = \frac{1}{2\pi} \langle q, \bar{\psi}_{I}(\lambda) \rangle_{L^2(\R)}, 
\end{cases}
\quad \lambda \in (0,\infty)
\end{align}
and 
\begin{align}
\label{projection-2}
\begin{cases}
a_{III}(\lambda) = \frac{(e^{\phi} - e^{-\phi})}{2\pi} \langle q, \bar{\psi}_{IV}(\lambda) \rangle_{L^2(\R)}, \\
a_{IV}(\lambda) = \frac{(e^{\phi} - e^{-\phi})}{2\pi} \langle q, \bar{\psi}_{III}(\lambda) \rangle_{L^2(\R)}, 
\end{cases}
\quad  \lambda \in (\lambda_0,\mu_0).
\end{align}
Substituting projections into (\ref{basis-2}) yields the identity
\begin{align*}
q(x) &= \frac{1}{2\pi} \int_{\mathbb{R}} q(y) \left( \int_0^{\infty} \left[ q_{I}(\lambda;x) \psi_{II}(\lambda;y) + q_{II}(\lambda;x) \psi_{I}(\lambda;y) \right] d\lambda \right. \\
& \qquad \left. + (e^{\phi} - e^{-\phi}) \int_{\lambda_0}^{\mu_0} \left[ q_{III}(\lambda;x) \psi_{IV}(\lambda;y) + q_{IV}(\lambda;x) \psi_{III}(\lambda;y) \right] d\lambda \right) dy. 
\end{align*}
Since the function $q \in L^2(\R) \cap L^1(\R)$ is arbitrary, the basis property of the eigenfunctions is equivalent to the completeness relation 
\begin{align}
& \int_0^{\infty} \left[ q_I(\lambda;x) \psi_{II}(\lambda;y) + q_{II}(\lambda;x) \psi_I(\lambda;y) \right] d\lambda  \notag \\
& + (e^{\phi} - e^{-\phi}) \int_{\lambda_0}^{\mu_0} \left[ q_{III}(\lambda;x) \psi_{IV}(\lambda;y) + q_{IV}(\lambda;x) \psi_{III}(\lambda;y) \right] d\lambda = 2\pi \delta(x-y). 
\label{complete-2}
\end{align}
The following proposition confirms the validity of the completeness relation (\ref{complete-2}).

\begin{proposition}
	\label{prop-complete-2}
	The completeness relation (\ref{complete-2}) is satisfied.
\end{proposition}

\begin{proof}
For the first term in (\ref{complete-2}), we obtain from (\ref{solution-2}) and (\ref{solution-1-adjont}) that
\begin{align*}
q_{I}(\lambda;x) \psi_{II}(\lambda;y) &= -i \frac{\partial}{\partial x} \frac{f(x) f'(y)}{f'(x) f(y)} \frac{e^{i \lambda  (x-y)}}{\lambda} \left[ 1 + \frac{k}{(\lambda - \mu_0) f(x)}\right] \left[ 1 - \frac{k}{(\lambda - \lambda_0) f'(y)} \right].
\end{align*}
Expanding the brackets with the help of (\ref{lambda-mu}) yields 
\begin{align*}
& \frac{1}{\lambda} \left[ 1 + \frac{k}{(\lambda - \mu_0) f(x)}\right] \left[ 1 - \frac{k}{(\lambda - \lambda_0) f'(y)} \right] \\
& \quad =  \frac{1}{\lambda} \left[ 1 + \frac{k \lambda (e^{iky - \phi} - e^{i k x + \phi}) + k\mu_0 e^{ikx + \phi} - k \lambda_0 e^{ik y-\phi}}{(\lambda - \lambda_0) (\lambda - \mu_0) f(x) f'(y)}  \right] \\
& \quad = \frac{1}{\lambda} 
+ \frac{k (e^{iky - \phi} - e^{i k x + \phi})}{(\lambda - \lambda_0) (\lambda - \mu_0) f(x) f'(y)} + \frac{k^2 (e^{iky}-e^{ikx})}{\lambda (\lambda - \lambda_0) (\lambda - \mu_0) f(x) f'(y) (e^{\phi} - e^{-\phi})}.
\end{align*}
By using partial fraction decomposition,
\begin{align*}
\frac{1}{\lambda (\lambda - \lambda_0) (\lambda - \mu_0)} = \frac{1}{\lambda_0 \mu_0 \lambda} + \frac{\mu_0 + \lambda_0 - \lambda}{\lambda_0 \mu_0 (\lambda - \lambda_0) (\lambda - \mu_0)}, \quad \lambda_0 \mu_0 = \frac{k^2}{(e^{\phi}-e^{-\phi})^2},
\end{align*}
we get with the help of $\mu_0 + \lambda_0 = -k \coth \phi$,
\begin{align*}
& \frac{1}{\lambda} \left[ 1 + \frac{k}{(\lambda - \mu_0) f(x)}\right] \left[ 1 - \frac{k}{(\lambda - \lambda_0) f'(y)} \right] \\
& \quad = \frac{f'(x) f(y)}{\lambda f(x) f'(y)} 
+ \frac{\lambda (e^{\phi}-e^{-\phi}) (e^{ikx} - e^{iky}) + k (e^{ikx-\phi}-e^{iky+\phi})}{(\lambda - \lambda_0) (\lambda - \mu_0) f(x) f'(y)} \\
& \quad = \frac{f'(x) f(y)}{\lambda f(x) f'(y)} 
+ \frac{(e^{\phi}-e^{-\phi})}{f(x) f'(y)} \left[ \frac{e^{ikx}}{\lambda - \lambda_0} - \frac{e^{iky}}{\lambda - \mu_0} \right],
\end{align*}
which yields 
\begin{align}
\label{term-first}
q_{I}(\lambda;x) \psi_{II}(\lambda;y) &= e^{i \lambda (x-y)} - i 
\frac{\partial}{\partial x} \frac{(e^{\phi}-e^{-\phi})}{f'(x) f(y)} \left[ \frac{e^{i(\lambda + k)x - i\lambda y}}{\lambda - \lambda_0} - \frac{e^{i \lambda x - i (\lambda - k) y}}{\lambda - \mu_0} \right].
\end{align}
Similarly, we compute for the second term in (\ref{complete-2}) that 
\begin{align}
q_{II}(\lambda;x) \psi_{I}(\lambda;y) &= i \frac{\partial}{\partial x} \frac{f'(x) f(y)}{f(x) f'(y)} \frac{e^{-i \lambda  (x-y)}}{\lambda} \left[ 1 - \frac{k}{(\lambda - \lambda_0) f'(x)}\right] \left[ 1 + \frac{k}{(\lambda - \mu_0) f(y)} \right] 
\notag \\
&= e^{-i \lambda (x-y)} + i 
\frac{\partial}{\partial x} \frac{(e^{\phi}-e^{-\phi})}{f(x) f'(y)} \left[ \frac{e^{-i \lambda x + i (\lambda + k) y}}{\lambda - \lambda_0} - \frac{e^{- i (\lambda - k) x + i \lambda y}}{\lambda - \mu_0} \right].
\label{term-second}
\end{align}	
By using (\ref{term-first}) and (\ref{term-second}), we obtain
	\begin{align*}
&	\int_0^{\infty} \left[ q_{I}(\lambda;x) \psi_{II}(\lambda;y) + q_{II}(\lambda;x) \psi_{I}(\lambda;y) \right] d\lambda = \int_0^{\infty} \left[ e^{i \lambda (x-y)} 
+ e^{-i \lambda (x-y)} \right] d\lambda \\ 
& \qquad  - i 
	\frac{\partial}{\partial x} \frac{(e^{\phi}-e^{-\phi})}{f'(x) f(y)} 
	\int_0^{\infty} \left[ \frac{e^{i(\lambda + k)x - i\lambda y}}{\lambda - \lambda_0} - \frac{e^{i \lambda x - i (\lambda - k) y}}{\lambda - \mu_0} \right] d \lambda \\
	& \qquad + i 
	\frac{\partial}{\partial x} \frac{(e^{\phi}-e^{-\phi})}{f(x) f'(y)} \int_0^{\infty} \left[ \frac{e^{-i \lambda x + i (\lambda + k) y}}{\lambda - \lambda_0} - \frac{e^{- i (\lambda - k) x + i \lambda y}}{\lambda - \mu_0} \right] d \lambda.
	\end{align*}
By using (\ref{distrib-1}), we have 
$$
\int_0^{\infty} \left[ e^{i \lambda (x-y)} 
+ e^{-i \lambda (x-y)} \right] d\lambda = \int_{\mathbb{R}} e^{i \lambda (x-y)} d \lambda = 2 \pi \delta(x-y).
$$
On the other hand, we transform the other two integrals as follows:
\begin{align*}
\int_0^{\infty} \left[ \frac{e^{i(\lambda + k)x - i\lambda y}}{\lambda - \lambda_0} - \frac{e^{i \lambda x - i (\lambda - k) y}}{\lambda - \mu_0} \right] d \lambda 
&= e^{i \mu_0 x - i \lambda_0 y} \left[ \int_{-\lambda_0}^{\infty} \frac{e^{i \lambda (x-y)}}{\lambda} d \lambda - \int_{-\mu_0}^{\infty} \frac{e^{i \lambda (x-y)}}{\lambda} d \lambda \right] \\
&= e^{i \mu_0 x - i \lambda_0 y} \int_{\lambda_0}^{\mu_0} \frac{e^{-i \lambda (x-y)}}{\lambda} d \lambda
\end{align*}
and 
\begin{align*}
\int_0^{\infty} \left[ \frac{e^{-i \lambda x + i (\lambda + k) y}}{\lambda - \lambda_0} - \frac{e^{- i (\lambda - k) x + i \lambda y}}{\lambda - \mu_0} \right] d \lambda 
&= e^{i \mu_0 y - i \lambda_0 x} \left[ \int_{-\lambda_0}^{\infty} \frac{e^{-i \lambda (x-y)}}{\lambda} d \lambda - \int_{-\mu_0}^{\infty} \frac{e^{-i \lambda (x-y)}}{\lambda} d \lambda \right] \\
&= e^{i \mu_0 y - i \lambda_0 x} \int_{\lambda_0}^{\mu_0} \frac{e^{i \lambda (x-y)}}{\lambda} d \lambda.
\end{align*}
Thus, we obtain 
\begin{align}
&	\int_0^{\infty} \left[ q_{I}(\lambda;x) \psi_{II}(\lambda;y) + q_{II}(\lambda;x) \psi_{I}(\lambda;y) \right] d\lambda =  2 \pi \delta(x-y) \notag  \\ 
& - i (e^{\phi}-e^{-\phi}) 
\frac{\partial}{\partial x} \frac{e^{i \mu_0 x - i \lambda_0 y}}{f'(x) f(y)} 
\int_{\lambda_0}^{\mu_0} \frac{e^{-i \lambda (x-y)}}{\lambda} d \lambda 
\notag \\
& + i 
(e^{\phi}-e^{-\phi}) \frac{\partial}{\partial x} \frac{e^{i \mu_0 y - i \lambda_0 x}}{f(x) f'(y)} \int_{\lambda_0}^{\mu_0} \frac{e^{i \lambda (x-y)}}{\lambda} d \lambda. \label{term-first-second}
\end{align}

Next, we show that the last two terms in (\ref{term-first-second}) are canceled in the sum of the third and fourth terms in (\ref{complete-2}).
These terms are represented by using \eqref{solution-3}, \eqref{solution-4}, (\ref{solution-adjoint-12}), \eqref{solution-adjoint-34}, (\ref{term-first}), and (\ref{term-second}):
\begin{align}
q_{IV}(\lambda;x) \psi_{III}(\lambda;y)
&= i \frac{\partial}{\partial x} \frac{e^{i \mu_0 x - i \lambda_0 y}}{f'(x) f(y)} \frac{\psi_{II}(\lambda;x) \psi_I(\lambda;y)}{\lambda} \notag \\
&= i\frac{\partial}{\partial x} \frac{e^{i \mu_0 x - i \lambda_0 y}}{f'(x) f(y)} \frac{e^{-i\lambda (x-y)}}{\lambda} \label{term-third} \\
& \quad + i (e^{\phi} - e^{-\phi}) \frac{\partial}{\partial x} \frac{e^{i \mu_0 x - i \lambda_0 y}}{f(x) f'(x) f(y) f'(y)}  \left[ \frac{e^{-i\lambda x + i (\lambda + k) y}}{\lambda - \lambda_0} - \frac{e^{-i (\lambda - k) x +i \lambda y}}{\lambda - \mu_0} \right] \notag
\end{align}
and 
\begin{align}
q_{III}(\lambda;x) \psi_{IV}(\lambda;y) 
&= -i \frac{\partial}{\partial x} \frac{e^{i \mu_0 y - i \lambda_0 x}}{f(x) f'(y)} \frac{\psi_{I}(\lambda;x) \psi_{II}(\lambda;y)}{\lambda} \notag \\
&= -i\frac{\partial}{\partial x} \frac{e^{i \mu_0 y - i \lambda_0 x}}{f(x) f'(y)} \frac{e^{i\lambda (x-y)}}{\lambda} \label{term-fourth}  \\
& \quad - i (e^{\phi} - e^{-\phi}) \frac{\partial}{\partial x} \frac{e^{i \mu_0 y - i \lambda_0 x}}{f(x) f'(x) f(y) f'(y)} \left[ \frac{e^{i (\lambda + k) x - i \lambda y}}{\lambda - \lambda_0} - \frac{e^{i\lambda x - i (\lambda - k) y}}{\lambda - \mu_0} \right]. \notag
\end{align}
Added together and integrated on $[\lambda_0,\mu_0]$, the second terms in (\ref{term-third}) and (\ref{term-fourth}) cancel out because 
\begin{align*}
& \int_{\lambda_0}^{\mu_0} \left( 
\left[ \frac{e^{-i\lambda x + i (\lambda + k) y}}{\lambda - \lambda_0} - \frac{e^{-i(\lambda-k)x + i\lambda y}}{\lambda - \mu_0} \right] 
e^{i \mu_0 x - i \lambda_0 y}
- \left[ \frac{e^{i (\lambda + k) x  - i \lambda y}}{\lambda - \lambda_0} - \frac{e^{i\lambda x - i (\lambda - k) y}}{\lambda - \mu_0} \right] 
e^{i \mu_0 y - i \lambda_0 x} \right) d \lambda \\
&= e^{ik (x+y)} \int_0^k \left[ \frac{e^{-i\lambda (x-y)}}{\lambda} + \frac{e^{i\lambda (x-y)}}{\lambda} -\frac{e^{i\lambda (x-y)}}{\lambda} - \frac{e^{-i\lambda (x-y)}}{\lambda} \right] d \lambda = 0,
\end{align*}
where we have used again that $\mu_0 = \lambda_0 + k$. Therefore, we obtain 
for the sum of the third and fourth terms in (\ref{complete-2}): 
\begin{align*}
& \int_{\lambda_0}^{\mu_0}  \left[ q_{III}(\lambda;x) \psi_{IV}(\lambda;y) + 
q_{IV}(\lambda;x) \psi_{III}(\lambda;y) \right] d\lambda \\
&= i\frac{\partial}{\partial x} \frac{e^{i \mu_0 x - i \lambda_0 y}}{f'(x) f(y)} \int_{\lambda_0}^{\mu_0} \frac{e^{-i\lambda (x-y)}}{\lambda} d \lambda 
-i\frac{\partial}{\partial x} \frac{e^{i \mu_0 y - i \lambda_0 x}}{f(x) f'(y)} \int_{\lambda_0}^{\mu_0} \frac{e^{i\lambda (x-y)}}{\lambda} d \lambda.
\end{align*}
Since these integrals are multiplied by $(e^{\phi} - e^{-\phi})$ in (\ref{complete-2}), they cancel the last two terms in (\ref{term-first-second}). 
This proves the validity of the completeness relation (\ref{complete-2}). 
\end{proof}

\section{Positivity of the second variation $S''(U)$ in Theorem \ref{th-main}}
\label{sec-5}

To prove the first assertion of Theorem \ref{th-main}, we first confirm that the spatial profile $U$ of the periodic wave is a critical point in the two commuting flows given by (\ref{BO-flow}) and (\ref{higher-BO-flow}). Then, we compute 
$\langle L q,q \rangle_{L^2(\R)}$ for every $q \in H^{1/2}(\R)$ and $\langle \tilde{L}q, q \rangle_{L^2(\R)}$ for every $q \in H^1(\R)$, where $L$ and $\tilde{L}$ are the two Hessian operators in the two variational formulations given by (\ref{Hessian}). Finally, we show that the linear combination of the two Hessian operators in the Lyapunov functional (\ref{Lyapunov}) yields positivity of the second variation for any $q \in H^1(\R)$.

\subsection{Variational formulations for the profile $U$ of the periodic wave}

Recall that $U$ is the spatial profile of the periodic waves  
$u(x,t) = U(\xi)$ with $\xi = x - ct$, which is given by the explicit formula (\ref{per-wave}). The following proposition gives two alternative variational characterizations  of the profile $U$ in the space of $2\pi/k$-periodic functions.

\begin{proposition}
	\label{lem-var}
	The periodic profile $U$ given by (\ref{per-wave}) is a critical point for $E(u) + c M(u)$ in $H^{1/2}_{\rm per}$ and a critical point for 
	$R(u) + \gamma M(u)$ in $H^1_{\rm per}$ with $\gamma = -\frac{1}{4} (k^2 + 3 c^2)$.
\end{proposition}

\begin{proof}
A critical point of $E(u) + c M(u)$ satisfies the stationary BO equation 
\begin{equation}
\label{BO-stat}
-H(u_x) - u^2 + c u = 0,
\end{equation}
which is obtained after integration of the traveling wave reduction 
of the BO equation
\begin{equation}
\label{BO-stat-der}
-cu_x + 2 u u_x + H(u_{xx}) = 0, 
\end{equation}
with the integration constant being zero. By using (\ref{Hirota}) and (\ref{Hirota-H}), the stationary BO equation (\ref{BO-stat}) is reduced to the bilinear form 
$$
(i c D_x + D_x^2) f' \cdot f = 0,
$$ 
which is satisfied by $f = 1 + e^{ik \xi + \phi}$ and $f' = 1 + e^{ik \xi - \phi}$ in agreement with (\ref{Hirota}). Hence, the integration constant of the traveling wave reduction (\ref{BO-stat-der}) is indeed zero. 

Since $H^2 = -I$ on periodic functions with zero mean, we derive for the traveling wave profile that 
\begin{align}
\quad & H(U') = c U - U^2, \label{b-1} \\ 
& \Rightarrow \;\; -U' = c H(U) - H(U^2), \notag \\
& \Rightarrow \;\;
2 H (U U') = c H(U') + U'', \notag \\
& \Rightarrow \;\;
2 H(U U') = c^2 U - c U^2 +  U''. \notag
\end{align}
This allows us to compute 
\begin{align*}
R'(U) &= -U'' + \frac{3}{2} U H(U') + \frac{3}{2} H(U U') + U^3 \\
&= -\frac{1}{4} \left[ U'' + 2 U^3 - 3c U^2 - 3 c^2 U \right].
\end{align*}
By using (\ref{per-wave}), we obtain 
$$
U'' + 2 U^3 = \frac{k^3 \sinh \phi [2 \cosh \phi - \cos(k\xi)]}{[\cos(k\xi) + \cosh \phi]^2}
$$
and 
$$
U'' + 2 U^3 - 3 c U^2 = -\frac{k^3 \sinh \phi}{\cos(k\xi) + \cosh \phi} = -k^2 U.
$$
Hence, the profile $U$ is a critical point of $R(u) + \gamma M(u)$ if and only if 
\begin{equation}
\label{b-0}
\frac{1}{4} (k^2 + 3 c^2) + \gamma = 0,
\end{equation}
which yields the unique expression for $\gamma$.
\end{proof}

\subsection{Second variation of $E(u) + cM(u)$}

By Proposition \ref{lem-var}, the periodic profile $U$ is a critical point 
of $E(u) + cM(u)$ in $H^{1/2}_{\rm per}$. The second variation of this functional at $U$ is defined by the Hessian operator $L : H^1(\R) \subset L^2(\R) \to L^2(\R)$ given by 
$$
L = E''(U) + c M''(U) = -H \partial_x - 2 U + c.
$$
Due to the duality of $H^{-1/2}(\R)$ and $H^{1/2}(\R)$ under the $L^2(\R)$ inner product, the quadratic form $\langle L q, q \rangle_{L^2(\R)}$ is well-defined for every $q \in H^{1/2}(\R)$. 

The linearized BO equation (\ref{BO-lin}) can be rewritten for the 
traveling waves with the profile $u(x,t) = U(x-ct)$ 
as a spectral stability problem for $q(x,t) = Q(x-ct) e^{t \Lambda}$, 
\begin{equation}
\partial_x L Q = \Lambda Q.
\label{spec-BO}
\end{equation}
It follows from Proposition \ref{prop-eigenfunctions-adjoint} that there exist four particular solutions of the spectral stability problem (\ref{spec-BO}) with the explicit expressions for eigenvalues 
\begin{align}
\label{eigenvalues-explicit-1}
\begin{cases}
\Lambda_I(\lambda) = i \lambda ( \lambda + c), \\
\Lambda_{II}(\lambda) = -i \lambda (\lambda + c), 
\end{cases} 
\quad \lambda \in [0,\infty), 
\end{align}
and
\begin{align}
\label{eigenvalues-explicit-2}
\begin{cases}
\Lambda_{III}(\lambda) = i (\lambda - \lambda_0) (\lambda - \mu_0), \\
\Lambda_{IV}(\lambda) = -i (\lambda - \lambda_0) (\lambda - \mu_0),
\end{cases}
\quad \lambda \in [\lambda_0,\mu_0],
\end{align}
as well as eigenfunctions 
\begin{align}
\label{eigenfunction-explicit-1}
\begin{cases}
Q_I(\lambda) = e^{i \lambda \xi} \frac{f}{f'} \left[ 1 + \frac{k}{(\lambda - \mu_0) f'} \right], \\
Q_{II}(\lambda) = e^{-i \lambda \xi} \frac{f'}{f} \left[ 1 - \frac{k}{(\lambda - \lambda_0) f} \right], 
\end{cases}
\quad \lambda \in [0,\infty), 
\end{align}
and
\begin{align}
\label{eigenfunction-explicit-2}
\begin{cases}
Q_{III}(\lambda) = e^{i (\lambda - \lambda_0) \xi} \left\{ \frac{1}{f'} \left[ 1 + \frac{k}{(\lambda - \mu_0) f'} \right] - \frac{\lambda_0}{\lambda f} \left[ 1 + \frac{k}{(\lambda - \mu_0) f} \right] \right\}, \\
Q_{IV}(\lambda) = e^{-i (\lambda - \mu_0) \xi} \left\{ \frac{1}{f} \left[ 1 - \frac{k}{(\lambda - \lambda_0) f} \right] - \frac{\mu_0}{\lambda f'} \left[ 1 - \frac{k}{(\lambda - \lambda_0) f'} \right] \right\}, 
\end{cases}
\quad \lambda \in (\lambda_0,\mu_0).
\end{align}
By using (\ref{solution-adjoint-12}) and (\ref{solution-adjoint-34}), we can write explicitly:
\begin{equation}
\label{operators-explicit}
\begin{cases}
L Q_I(\lambda) = (\lambda + c) \Psi_{I}(\lambda), \\
L Q_{II}(\lambda) =  (\lambda + c) \Psi_{II}(\lambda),
\end{cases}
\quad 
\begin{cases}
L Q_{III}(\lambda) = \frac{(\lambda - \lambda_0) (\lambda - \mu_0)}{\lambda} \Psi_{III}(\lambda), \\
L Q_{IV}(\lambda) = \frac{(\lambda - \lambda_0) (\lambda - \mu_0)}{\lambda} \Psi_{IV}(\lambda), 
\end{cases}
\end{equation}
where the adjoint eigenfunctions follow from Proposition \ref{prop-eigenfunctions}:
\begin{align}
\label{eigenfunction-adjoint-explicit-1}
\begin{cases}
\Psi_{I}(\lambda) = e^{i \lambda \xi} \frac{f}{f'} \left[ 1 + \frac{k}{(\lambda - \mu_0) f} \right], \\
\Psi_{II}(\lambda) = e^{-i \lambda \xi} \frac{f'}{f} \left[ 1 - \frac{k}{(\lambda - \lambda_0) f'} \right],
\end{cases}
\quad \lambda \in [0,\infty), 
\end{align}
and 
\begin{align}
\label{eigenfunction-adjoint-explicit-2}
\begin{cases}
\Psi_{III}(\lambda) = e^{i (\lambda - \lambda_0) \xi} \frac{1}{f'} \left[ 1 + \frac{k}{(\lambda - \mu_0) f} \right], \\
\Psi_{IV}(\lambda) = e^{-i (\lambda - \mu_0) \xi} \frac{1}{f} \left[ 1 - \frac{k}{(\lambda - \lambda_0) f'} \right],
\end{cases}
\quad \lambda \in (\lambda_0,\mu_0).
\end{align}
We also collect together the complex-conjugate relations (\ref{conjugation-1}), (\ref{conjugation-2}), (\ref{conjugation-1-q}), and (\ref{conjugation-2-q}):
\begin{align}
\label{complex-conjugate}
\begin{cases}
\bar{Q}_I(\lambda) = e^{2 \phi} \frac{\lambda - \lambda_0}{\lambda - \mu_0} Q_{II}(\lambda), \\
\bar{Q}_{II}(\lambda) = e^{-2 \phi} \frac{\lambda - \mu_0}{\lambda - \lambda_0} Q_{I}(\lambda),
\end{cases}
\quad 
 \begin{cases}
\bar{Q}_{III}(\lambda) = e^{\phi} \frac{\lambda - \lambda_0}{\lambda - \mu_0} Q_{IV}(\lambda), \\
\bar{Q}_{IV}(\lambda) = e^{-\phi} \frac{\lambda - \mu_0}{\lambda - \lambda_0} Q_{III}(\lambda),
\end{cases}
\end{align}
and
\begin{align}
\label{complex-conjugate-adjoint}
\begin{cases}
\bar{\Psi}_I(\lambda) = e^{2 \phi} \frac{\lambda - \lambda_0}{\lambda - \mu_0} \Psi_{II}(\lambda), \\
\bar{\Psi}_{II}(\lambda) = e^{-2 \phi} \frac{\lambda - \mu_0}{\lambda - \lambda_0} \Psi_{I}(\lambda), 
\end{cases}
\quad 
\begin{cases}
\bar{\Psi}_{III}(\lambda) = e^{\phi} \frac{\lambda - \lambda_0}{\lambda - \mu_0} \Psi_{IV}(\lambda), \\
\bar{\Psi}_{IV}(\lambda) = e^{-\phi} \frac{\lambda - \mu_0}{\lambda - \lambda_0} \Psi_{III}(\lambda).
\end{cases}
\end{align}

The following proposition uses these explicit expressions to diagonalize 
the quadratic form $\langle L q, q \rangle_{L^2(\R)}$ for $q \in H^{1/2}(\R)$ by using the basis of eigenfunctions of the spectral problem (\ref{spec-BO}).

\begin{proposition}
	\label{prop-diagoalization-L}
	For every $q \in H^{1/2}(\R)$, we have
	\begin{align}
	& \langle L q, q \rangle_{L^2(\R)} = 
	2\pi \int_0^{\infty} (\lambda + c) \left[ |a_{I}(\lambda)|^2 e^{2 \phi} \frac{\lambda - \lambda_0}{\lambda - \mu_0} + |a_{II}(\lambda)|^2 e^{-2 \phi} \frac{\lambda - \mu_0}{\lambda - \lambda_0} \right] d \lambda \notag \\
	& + 2\pi
	\int_{\lambda_0}^{\mu_0} \frac{(\lambda - \lambda_0) (\lambda - \mu_0)}{ (e^{\phi} - e^{-\phi}) \lambda} \left[ |a_{III}(\lambda)|^2 e^{\phi} \frac{\lambda - \lambda_0}{\lambda - \mu_0} + |a_{IV}(\lambda)|^2 e^{-\phi} \frac{\lambda - \mu_0}{\lambda - \lambda_0} \right] d \lambda,
	\label{quadratic-form-L}
	\end{align}
	where $a_I(\lambda)$, $a_{II}(\lambda)$, $a_{III}(\lambda)$, $a_{IV}(\lambda)$ are uniquely defined by $q \in H^{1/2}(\R)$.
	Consequently, the first integral in (\ref{quadratic-form-L}) is positive semi-definite and the second integral in (\ref{quadratic-form-L}) is negative semi-definite. 
\end{proposition}

\begin{proof}
By Proposition \ref{prop-complete-2}, we can represent $q \in H^{1/2}(\R)$ by using the basis of eigenfunctions in the form:
	\begin{equation}
	\label{q-expansion}
	q = \int_0^{\infty} [a_{I}(\lambda) Q_I(\lambda) + a_{II}(\lambda) Q_{II}(\lambda)] d \lambda + \int_{\lambda_0}^{\mu_0} [a_{III}(\lambda) Q_{III}(\lambda) + a_{IV}(\lambda) Q_{IV}(\lambda)] d \lambda,
	\end{equation}
with some uniquely defined coefficients $a_I(\lambda)$, $a_{II}(\lambda)$, $a_{III}(\lambda)$, $a_{IV}(\lambda)$, see (\ref{projection-1}) and 
(\ref{projection-2}) which hold if $q \in L^2(\R) \cap L^1(\R)$.
By using (\ref{operators-explicit}) and (\ref{q-expansion}), we obtain 
\begin{align*}
L q &= \int_0^{\infty} [a_{I}(\lambda) L Q_I(\lambda) + a_{II}(\lambda) L Q_{II}(\lambda)] d \lambda + \int_{\lambda_0}^{\mu_0} [a_{III}(\lambda) L Q_{III}(\lambda) + a_{IV}(\lambda) L Q_{IV(\lambda)}] d \lambda\\
&= \int_0^{\infty} (\lambda + c) [a_{I}(\lambda) \Psi_I(\lambda) + a_{II}(\lambda) \Psi_{II}(\lambda)] d \lambda \\
& \quad + \int_{\lambda_0}^{\mu_0} \frac{(\lambda - \lambda_0) (\lambda - \mu_0)}{\lambda}  [a_{III}(\lambda) \Psi_{III}(\lambda) + a_{IV}(\lambda) \Psi_{IV}(\lambda)] d \lambda.
\end{align*}
By using (\ref{complex-conjugate}), we take the complex conjugation of (\ref{q-expansion}) and obtain
\begin{align}
\bar{q} &= \int_0^{\infty} \left[\bar{a}_{I}(\lambda) e^{2 \phi} \frac{\lambda - \lambda_0}{\lambda - \mu_0} Q_{II}(\lambda) + \bar{a}_{II}(\lambda) e^{-2 \phi} \frac{\lambda - \mu_0}{\lambda - \lambda_0} Q_{I}(\lambda) \right] d \lambda 
\notag \\
& \quad  + \int_{\lambda_0}^{\mu_0} \left[ \bar{a}_{III}(\lambda) e^{\phi} \frac{\lambda - \lambda_0}{\lambda - \mu_0} Q_{IV}(\lambda) + \bar{a}_{IV}(\lambda) e^{-\phi} \frac{\lambda - \mu_0}{\lambda - \lambda_0} Q_{III}(\lambda) \right] d \lambda. \label{q-bar-expansion}
\end{align}
Substituting (\ref{q-bar-expansion}) into $\langle L q, q\rangle_{L^2(\R)}$ and using the orthogonality relations of Propositions \ref{prop-orthog-1}, \ref{prop-orthog-2}, and \ref{prop-orthog-3}, we obtain 
the diagonalization formula (\ref{quadratic-form-L}).
\end{proof}

\subsection{Second variation of $R(u) + \gamma M(u)$}

By Proposition \ref{lem-var}, the periodic profile $U$ is a critical point 
of $R(u) + \gamma M(u)$ in $H^1_{\rm per}$ 
with $\gamma = -\frac{1}{4} (k^2 + 3 c^2)$. 
The second variation of this functional at $U$ is defined by 
the Hessian operator $\tilde{L}: H^2(\R) \subset L^2(\R) \to L^2(\R)$ given by 
$$
\tilde{L} = R''(U) + \gamma M''(U) = -\partial_x^2 + \frac{3}{2} H(U') 
+ \frac{3}{2} U H \partial_x + \frac{3}{2} H \partial_x (U \cdot ) + 3 U^2 + \gamma.
$$
Due to the duality of $H^{-1}(\R)$ and $H^{1}(\R)$ under the $L^2(\R)$ inner product, the quadratic form $\langle \tilde{L} q, q \rangle_{L^2(\R)}$ is well-defined for every $q \in H^{1}(\R)$. 

The eigenfunctions (\ref{eigenfunction-explicit-1}) and (\ref{eigenfunction-explicit-2}) of $\partial_x L$ can be rewritten 
in the Floquet--Bloch form as
\begin{align}
\label{eigenfunction-Floquet-1}
\begin{cases}
Q_{I,n}(\theta) = e^{i (k n+\theta) \xi} \frac{f}{f'} \left[ 1 + \frac{k}{(k n + \theta - \mu_0) f'} \right], \\
Q_{II,n}(\theta) = e^{-i (k n + \theta) \xi} \frac{f'}{f} \left[ 1 - \frac{k}{(k n + \theta - \lambda_0) f} \right], \\
\end{cases}
\quad \lambda = k n + \theta, \quad \theta \in [0,k),
\end{align}
with integer $n \geq 0$, and
\begin{align}
\label{eigenfunction-Floquet-2}
\begin{cases}
Q_{III}(\theta) = e^{i \theta \xi} \left\{ \frac{1}{f'} \left[ 1 + \frac{k}{(\theta - k) f'} \right] - \frac{\lambda_0}{(\lambda_0 + \theta) f} \left[ 1 + \frac{k}{(\theta - k) f} \right] \right\}, \\
Q_{IV}(\theta)= e^{-i \theta \xi + ik\xi} \left\{ \frac{1}{f} \left[ 1 - \frac{k}{\theta f} \right] - \frac{\mu_0}{(\lambda_0 + \theta) f'} \left[ 1 - \frac{k}{\theta f'} \right] \right\},
\end{cases}
\lambda = \lambda_0 + \theta, \;\; \theta \in [0,k),
\end{align}
where $Q_{IV}(\theta)$ can be renormalized as $\frac{\theta}{\theta - k} Q_{IV}(\theta)$ to avoid the divergence as $\theta \to 0$. This modification does not affect the rest of our analysis.

When $\theta \in [0,k)$ is fixed, the eigenvalues of $\partial_x L$ given by (\ref{eigenvalues-explicit-1}) and (\ref{eigenvalues-explicit-2}) form simple sequences
\begin{align}
\label{eigenvalues-Floquet-1}
\begin{cases}
\Lambda_{I,n}(\theta)= i (k n+\theta) (k n + \theta + c), \\
\Lambda_{II,n}(\theta) = -i  (k n+\theta) (k n + \theta + c), 
\end{cases} 
\quad \theta \in [0,k),
\end{align}
where $n \geq 0$ is an integer, and
\begin{align}
\label{eigenvalues-Floquet-2}
\begin{cases}
\Lambda_{III}(\theta) = i \theta (\theta - k), \\
\Lambda_{IV}(\theta) = -i  \theta (\theta - k),
\end{cases}
\quad \theta \in [0,k).
\end{align}

It follows from the existence of a recursion operator between two commuting flows (\ref{BO-flow}) and (\ref{higher-BO-flow}), see \cite[Eqs. (1.12), (2.21), (2.23)]{LanWang2025}, that the linearized operators $\partial_x L$ and $\partial_x \tilde{L}$ commute. In the Floquet--Bloch form for a fixed $\theta \in [0,k)$, the eigenfunctions $\{ Q_{I,n}(\theta) \}_{n \geq 0}$ and $Q_{III}(\theta)$ in (\ref{eigenfunction-Floquet-1}) and (\ref{eigenfunction-Floquet-2}) satisfy 
$$
Q_{I,n}(\theta;\xi) = e^{i \theta \xi}\tilde{Q}_{I,n}(\theta;\xi), \quad 
Q_{III}(\theta;\xi) = e^{i \theta \xi}\tilde{Q}_{III}(\theta;\xi), \quad \theta \in [0,k),
$$
where $\tilde{Q}_{I,n}(\theta;\cdot),\tilde{Q}_{III}(\theta;\cdot) \in L^2_{\rm per}$ correspond to simple eigenvalues $\{ \Lambda_{I,n}(\theta) \}_{n \geq 0}$ and $\Lambda_{III}(\theta)$
given by (\ref{eigenvalues-Floquet-1}) and (\ref{eigenvalues-Floquet-2}).
The complex-conjugate eigenfunctions $\{ Q_{II,n}(\theta) \}_{n \geq 0}$ and $Q_{IV}(\theta)$  correspond to the simple complex-conjugate eigenvalues $\{ \Lambda_{II,n}(\theta) \}_{n \geq 0}$ and $\Lambda_{IV}(\theta)$ and satisfy
$$
Q_{II,n}(\theta;\xi) = e^{-i \theta \xi}\tilde{Q}_{II,n}(\theta;\xi), \quad 
Q_{IV}(\theta;\xi) = e^{-i \theta \xi}\tilde{Q}_{IV}(\theta;\xi), \quad \theta \in [0,k),
$$
where $\tilde{Q}_{II,n}(\theta;\cdot),\tilde{Q}_{IV}(\theta;\cdot) \in L^2_{\rm per}$. Due to simplicity of eigenvalues for $\partial_x L$, the two commuting operators $\partial_x L$ and $\partial_x \tilde{L}$ share a common basis of eigenfunctions in the Floquet--Bloch form with $\theta \in [0,k)$ with the same set of eigenfunctions (\ref{eigenfunction-Floquet-1}) and (\ref{eigenfunction-Floquet-2}). Therefore, these eigenfunctions  satisfy  the spectral problem 
\begin{equation}
\partial_x \tilde{L} Q = \tilde{\Lambda}  Q.
\label{higher-spec-BO}
\end{equation}
However, the eigenvalues $\tilde{\Lambda}$ in (\ref{higher-spec-BO}) are different from the eigenvalues $\Lambda$  in (\ref{spec-BO}). 
The following proposition gives the explicit expressions for 
eigenvalues $\tilde{\Lambda}$. 

\begin{proposition}
	\label{prop-eigenvalues}
The spectral problem (\ref{higher-spec-BO}) is satisfied with the eigenfunctions (\ref{eigenfunction-explicit-1}) and (\ref{eigenfunction-explicit-2}) for the eigenvalues:
	\begin{align}
	\label{eigenvalues-explicit-tilde-1}
	\begin{cases}
	\tilde{\Lambda}_I(\lambda) =  i \lambda (\lambda^2 + \gamma),\\
	\tilde{\Lambda}_{II}(\lambda) =  -i \lambda (\lambda^2 + \gamma),
	\end{cases}
	\quad \lambda \in [0,\infty),
\end{align}
and	
	\begin{align}
\label{eigenvalues-explicit-tilde-2}
\begin{cases}
	\tilde{\Lambda}_{III}(\lambda) = i (\lambda - \lambda_0) (\lambda - \mu_0) (\lambda - c),\\
	\tilde{\Lambda}_{IV}(\lambda) = -i (\lambda - \lambda_0) (\lambda - \mu_0) (\lambda - c),
	\end{cases}
		\quad \lambda \in [\lambda_0,\mu_0],
	\end{align}
such that 
\begin{align}
\label{operators-explicit-tilde}
\begin{cases}
\tilde{L} Q_I(\lambda) = (\lambda^2 + \gamma) \Psi_{I}(\lambda), \\
\tilde{L} Q_{II}(\lambda) =  (\lambda^2 + \gamma) \Psi_{II}(\lambda), 
\end{cases}
\quad 
\begin{cases}
\tilde{L} Q_{III}(\lambda) = \frac{(\lambda - \lambda_0) (\lambda - \mu_0) (\lambda - c)}{\lambda} \Psi_{III}(\lambda), \\
\tilde{L} Q_{IV}(\lambda) = \frac{(\lambda - \lambda_0) (\lambda - \mu_0) (\lambda - c)}{\lambda} \Psi_{IV}(\lambda).
\end{cases} 
\end{align}
\end{proposition}

\begin{proof}
	We recall that $\tilde{L}$ is a linear operator with real-valued coefficients. Since $Q_{II}(\lambda)$ and $\Psi_{II}(\lambda)$ are proportional to the complex-conjugate versions of $Q_I(\lambda)$ and $\Psi_I(\lambda)$ in (\ref{complex-conjugate}) and (\ref{complex-conjugate-adjoint}), whereas $Q_{IV}(\lambda)$ and $\Psi_{IV}(\lambda)$ are proportional to the complex-conjugate versions of $Q_{III}(\lambda)$ and $\Psi_{III}(\lambda)$, we obtain $\tilde{\Lambda}_{II}(\lambda) = -\tilde{\Lambda}_{I}(\lambda)$ and $\tilde{\Lambda}_{IV}(\lambda) = -\tilde{\Lambda}_{III}(\lambda)$. Therefore, it suffices to prove (\ref{eigenvalues-explicit-tilde-1}), (\ref{eigenvalues-explicit-tilde-2}), and (\ref{operators-explicit-tilde}) for $\tilde{\Lambda}_I(\lambda)$ and $\tilde{\Lambda}_{III}(\lambda)$ only. 
	
We compute these eigenvalues from the asymptotic behavior of the spectral problem (\ref{higher-spec-BO}) in the complex $x$-plane at infinity. We have also checked validity of the relations (\ref{operators-explicit-tilde}) with  direct computations which are included in Appendix \ref{app-A} for completeness. \\

\underline{Eigenvalue $\tilde{\Lambda}_I(\lambda)$:} We note that $Q_I(\lambda)$ is analytic in $\mathbb{C}^+$ in the complex $x$-plane and since $\lambda > 0$, the exponential factor in $Q_I(\lambda)$ is bounded and decaying as ${\rm Im}(x) \to \infty$. Hence, 
$H Q_I = i Q_I$ so that $L Q_I$ is equivalent to 
$$
L Q_I = -i Q_I' -2 U  Q_I + c Q_I.
$$
As ${\rm Im}(x) \to -\infty$, we have exponentially large terms in $f \sim e^{i k \xi + \phi}$, $f' \sim e^{i k \xi - \phi}$, so that 
\begin{equation}
\label{first-beh}
Q_I(\lambda;\xi) \sim e^{i \lambda \xi + 2 \phi}, \quad \Psi_I(\lambda;\xi) \sim e^{i \lambda \xi + 2 \phi} \quad \mbox{\rm as} \;\; {\rm Im}(x) \to -\infty.
\end{equation}
Since 
$$
U(\xi) = \frac{k \sinh \phi}{\cos(k \xi) + \cosh \phi}, \quad 
H(U')(\xi) = k^2 \frac{1 + \cos(k\xi) \cosh \phi}{[\cos(k\xi) + \cosh \phi]^2},
$$
the term $-2U Q_I$ in $LQ_{I}$ is smaller than $-iQ_I' + c Q_I$ by an exponential factor as ${\rm Im}(x) \to -\infty$. The balance of exponential terms in the limit ${\rm Im}(x) \to -\infty$ for 
\begin{equation}
\label{first-eig}
L Q_I = -i Q_I' -2 U  Q_I + c Q_I = -\frac{i \Lambda_{I}}{\lambda} \Psi_{I}
\end{equation}
yields $\Lambda_I = i \lambda (\lambda + c)$, in agreement with (\ref{eigenvalues-explicit-1}). For the operator $\tilde{L}$, we write 
$$
\tilde{L} Q_I  = -Q_I'' + \frac{3}{2} Q_I H(U')  + \frac{3}{2} U 
H(Q_I') + \frac{3}{2} H(U Q_I)' + 3 U^2 Q_I + \gamma Q_I.
$$
It follows from (\ref{first-eig}) after using the Hilbert transform that 
\begin{align}
2 H(U Q_I) &= Q_I' + c H(Q_I) - (\lambda + c) H(\Psi_I) \notag \\
&= Q_I' + i c Q_I - i (\lambda + c) \Psi_I, \label{b-2}
\end{align}
where we used the analyticity of $Q_I$ and $\Psi_I$ in $\mathbb{C}^+$ and their exponential decay to $0$ as ${\rm Im}(x) \to +\infty$. By using the asymptotic behavior (\ref{first-beh}), we confirm that 
$H(U Q_I)$ is smaller than $Q_I$ by an exponential factor as ${\rm Im}(x) \to -\infty$, similar to $Q_I H(U')$, $U H(Q_I') = iU Q_I'$, and $U^2 Q_I$. The balance of exponential terms in the limit ${\rm Im}(x) \to -\infty$ for 
\begin{equation*}
\tilde{L} Q_I = -\frac{i \tilde{\Lambda}_{I}}{\lambda} \Psi_{I}
\end{equation*}
yields $\tilde{\Lambda}_I = i \lambda (\lambda^2 + \gamma)$, in agreement with (\ref{eigenvalues-explicit-tilde-1}).\\

\underline{Eigenvalues $\tilde{\Lambda}_{III}(\lambda)$:} The eigenfunction $Q_{III}(\lambda)$ consists of the sum of two terms. The first term 
$$
e^{i (\lambda - \lambda_0) \xi} \frac{1}{f'} \left[ 1 + \frac{k}{(\lambda - \mu_0) f'} \right]
$$
is analytic in $\mathbb{C}^+$ and since $\lambda > \lambda_0$, the exponential factor is bounded and decaying as ${\rm Im}(x) \to \infty$. On the other hand, 
the second term 
$$
- e^{i (\lambda - \lambda_0) \xi} \frac{\lambda_0}{\lambda f} \left[ 1 + \frac{k}{(\lambda - \mu_0) f} \right] = -e^{i (\lambda - \mu_0) \xi - \phi} \frac{\lambda_0}{\lambda (1 + e^{-ik\xi - \phi})} \left[ 1 + \frac{k}{(\lambda - \mu_0) f} \right]
$$
is analytic in $\mathbb{C}^-$ and since $\lambda < \mu_0$, the exponential factor is bounded and decaying as ${\rm Im}(x) \to -\infty$. As a result, we obtain 
\begin{align*}
H(Q_{III}) = i e^{i (\lambda - \lambda_0) \xi} \left\{ \frac{1}{f'} \left[ 1 + \frac{k}{(\lambda - \mu_0) f'} \right] + \frac{\lambda_0}{\lambda f} \left[ 1 + \frac{k}{(\lambda - \mu_0) f} \right] \right\}.
\end{align*}
As ${\rm Im}(x) \to -\infty$, we have exponentially large terms in $f \sim e^{i k \xi + \phi}$, $f' \sim e^{i k \xi - \phi}$, so that 
\begin{equation}
\label{third-beh}
\begin{cases} 
Q_{III}(\lambda;\xi) \sim \frac{\lambda - \mu_0}{\lambda} e^{i (\lambda - \mu_0) \xi +  \phi}, \\ 
\Psi_{III}(\lambda;\xi) \sim e^{i (\lambda - \mu_0) \xi +  \phi}, 
\end{cases}
\quad \mbox{\rm as} \;\; {\rm Im}(x) \to -\infty
\end{equation}
and 
\begin{equation}
\label{third-beh-2}
H(Q_{III})(\lambda;\xi) \sim i \frac{\lambda + \mu_0}{\lambda} e^{i (\lambda - \mu_0) \xi +  \phi}, \quad \mbox{\rm as} \;\; {\rm Im}(x) \to -\infty,
\end{equation}
where we have used $\lambda_0 e^{-2\phi} = \mu_0$. Note that 
$Q_{III}$, $H(Q_{III})$, and $\Psi_{III}$ are exponentially small 
as ${\rm Im}(x) \to -\infty$. Since $-2U Q_{III}$ in $LQ_{III}$ is smaller than $-H(Q_{III}') + c Q_{III}$ by an exponential factor as ${\rm Im}(x) \to -\infty$, the balance of exponential terms in the limit ${\rm Im}(x) \to -\infty$ for 
\begin{equation}
\label{third-eig}
L Q_{III} = -H(Q_{III}') -2U Q_{III} + c Q_{III} = -\frac{i \Lambda_{III}}{\lambda} \Psi_{III}
\end{equation}
yields
\begin{align*}
\Lambda_{III} = i (\lambda - \mu_0) (\lambda + \mu_0 + c) = i (\lambda - \lambda_0) (\lambda - \mu_0),
\end{align*}
in agreement with (\ref{eigenvalues-explicit-2}), where we used $c + \lambda_0 + \mu_0 = 0$. For the operator $\tilde{L}$, we write again
$$
\tilde{L} Q_{III}  = -Q_{III}'' + \frac{3}{2} Q_{III} H(U')  + \frac{3}{2} U 
H(Q_{III}') + \frac{3}{2} H(U Q_{III})' + 3 U^2 Q_{III} + \gamma Q_{III},
$$
and note that $Q_{III} H(U')$, $U H(Q_{III}')$, and $U^2 Q_{III}$ 
are smaller than $-Q_{III}'' + \gamma Q_{III}$ by an exponential factor 
as ${\rm Im}(x) \to -\infty$. On the other hand, the term 
$\frac{3}{2} H(U Q_{III})'$ cannot be neglected this time. 
It follows from (\ref{third-eig}) after using the Hilbert transform that 
\begin{align}
2 H(U Q_{III}) = Q_{III}' + c H(Q_{III}) - \lambda^{-1} (\lambda -  \lambda_0) (\lambda - \mu_0) H(\Psi_{III}).
\label{b-01}
\end{align}
The eigenfunction $\Psi_{III}$ consists of the sum of two terms. The first term $\frac{1}{f'} e^{i (\lambda - \lambda_0) \xi}$
is analytic in $\mathbb{C}^+$ and since $\lambda > \lambda_0$, the exponential factor is bounded and decaying as ${\rm Im}(x) \to \infty$. On the other hand, 
the second term can be decomposed by using (\ref{c-3.15}) into the sum of two terms:
\begin{align*}
\frac{k}{(\lambda - \mu_0) f f'} e^{i (\lambda - \lambda_0) \xi} &= 
\frac{i}{(\lambda - \mu_0) (e^{\phi} - e^{-\phi})} e^{i (\lambda - \mu_0) \xi} 
\left( \frac{f'_x}{f'} - \frac{f_x}{f} \right) \\
&= -\frac{k}{(\lambda - \mu_0) (e^{\phi} - e^{-\phi})} 
\left[ \frac{e^{i (\lambda - \lambda_0) \xi - \phi}}{f'} - \frac{e^{i(\lambda - \mu_0) \xi}}{1 + e^{-ik\xi - \phi}} \right],
\end{align*}
where the first term is analytic in $\mathbb{C}^+$ and the second term is  analytic in $\mathbb{C}^-$. Moreover, the exponential factors are bounded and decaying as ${\rm Im}(x) \to +\infty$ for the first term and 
as ${\rm Im}(x) \to -\infty$ for the second term. As a result, we obtain 
\begin{align}
\label{b-02}
H(\Psi_{III}) = i e^{i (\lambda - \lambda_0) \xi} \left[ 
\frac{1}{f'}  - \frac{k  e^{- \phi}}{(\lambda - \mu_0)(e^{\phi} - e^{-\phi}) f'}- \frac{k e^{+ \phi}}{(\lambda - \mu_0)(e^{\phi} - e^{-\phi}) f}  \right].
\end{align}
As ${\rm Im}(x) \to -\infty$, we have exponentially large terms in $f \sim e^{i k \xi + \phi}$, $f' \sim e^{i k \xi - \phi}$, so that 
\begin{align}
H(\Psi_{III})(\lambda;\xi) &\sim i e^{i (\lambda - \mu_0) \xi +  \phi} \left[ 1 - \frac{2k e^{-\phi}}{(\lambda - \mu_0)(e^{\phi} - e^{-\phi})} \right] \notag \\
&= 
i \frac{\lambda + \mu_0}{\lambda - \mu_0} e^{i (\lambda - \mu_0) \xi +  \phi}, \quad \mbox{\rm as} \;\; {\rm Im}(x) \to -\infty,
\label{third-beh-3}
\end{align}
where we have used (\ref{lambda-mu}) for $\mu_0$. 
By using (\ref{third-beh}), (\ref{third-beh-2}), and (\ref{third-beh-3}) 
we obtain the principal term in the asymptotic expansion of
\begin{align*}
2 H(U Q_{III}) &\sim i \lambda^{-1} e^{i (\lambda - \mu_0) \xi +  \phi} 
\left[ (\lambda - \mu_0)^2 + c (\lambda + \mu_0) - (\lambda - \lambda_0) (\lambda + \mu_0) \right] \\ 
&= - 4 i \mu_0  e^{i (\lambda - \mu_0) \xi +  \phi}, \quad \mbox{\rm as} \;\; {\rm Im}(x) \to -\infty,
\end{align*}
where we have used $c + \lambda_0 + \mu_0 = 0$. The balance of exponential terms in the limit ${\rm Im}(x) \to -\infty$ for 
\begin{equation*}
\tilde{L} Q_{III} = -\frac{i \tilde{\Lambda}_{III}}{\lambda} \Psi_{III}
\end{equation*}
yields
\begin{align*}
\tilde{\Lambda}_{III} &= i (\lambda - \mu_0) \left( (\lambda - \mu_0)^2 + \gamma + 3 \mu_0 \lambda \right) \\
&= i (\lambda - \mu_0) ( \lambda^2 + \mu_0 \lambda + \mu_0^2 + \gamma) \\
&= i (\lambda - \mu_0) (\lambda - \lambda_0) (\lambda + \lambda_0 + \mu_0) \\
&= i (\lambda - \lambda_0) (\lambda - \mu_0) (\lambda - c),
\end{align*}
in agreement with (\ref{eigenvalues-explicit-tilde-2}), 
where we have used 
$$
\gamma = -\frac{1}{4} (k^2 + 3 c^2) = -(\lambda_0^2 + \lambda_0 \mu_0 + \mu_0^2),
$$
which follows from (\ref{Lax-spectrum-continuous}). 
\end{proof}

The following proposition uses the explicit expressions for 
the eigenvalues in Proposition \ref{prop-eigenvalues} to diagonalize 
the quadratic form $\langle \tilde{L} q, q \rangle_{L^2(\R)}$ for $q \in H^{1}(\R)$ by using the basis of eigenfunctions of the spectral problem (\ref{spec-BO}).

\begin{proposition}
	\label{prop-diagoalization-L-tilde}
	For every $q \in H^{1}(\R)$, we have
	\begin{align}
	& \langle \tilde{L} q, q \rangle_{L^2(\R)} = 
	2\pi \int_0^{\infty} (\lambda^2 + \gamma) \left[ |a_I(\lambda)|^2 e^{2 \phi} \frac{\lambda - \lambda_0}{\lambda - \mu_0} + |a_{II}(\lambda)|^2 e^{-2 \phi} \frac{\lambda - \mu_0}{\lambda - \lambda_0} \right] d \lambda  \notag \\
	& + 2\pi
	\int_{\lambda_0}^{\mu_0} \frac{(\lambda - \lambda_0) (\lambda - \mu_0) (\lambda - c)}{(e^{\phi} - e^{-\phi}) \lambda} \left[ |a_{III}(\lambda)|^2 e^{\phi} \frac{\lambda - \lambda_0}{\lambda - \mu_0} + |a_{IV}(\lambda)|^2 e^{-\phi} \frac{\lambda - \mu_0}{\lambda - \lambda_0} \right] d \lambda,
	\label{quadratic-form-L-tilde}
	\end{align}
	where $a_I(\lambda)$, $a_{II}(\lambda)$, $a_{III}(\lambda)$,  $a_{IV}(\lambda)$ are uniquely defined by $q \in H^{1}(\R)$.
	Consequently, the first integral in (\ref{quadratic-form-L-tilde}) is sign-indefinite and the second integral is positive semi-definite.
\end{proposition}

\begin{proof}
By using (\ref{q-expansion}) and (\ref{operators-explicit-tilde}), we obtain 
\begin{align*}
\tilde{L} q &= \int_0^{\infty} [a_I(\lambda) \tilde{L} Q_I(\lambda) + a_{II}(\lambda) \tilde{L} Q_{II}(\lambda)] d \lambda + \int_{\lambda_0}^{\mu_0} [a_{III}(\lambda) \tilde{L} Q_{III}(\lambda) + a_{IV}(\lambda) \tilde{L} Q_{IV}(\lambda)] d \lambda \\
&= \int_0^{\infty} (\lambda^2 + \gamma) [a_I(\lambda) \Psi_I(\lambda) + a_{II}(\lambda) \Psi_{II}(\lambda)] d \lambda \\
& \quad + \int_{\lambda_0}^{\mu_0} \frac{(\lambda - \lambda_0) (\lambda - \mu_0) (\lambda - c)}{\lambda}  [a_{III}(\lambda) \Psi_{III}(\lambda) + a_{IV}(\lambda) \Psi_{IV}(\lambda)] d \lambda.
\end{align*}
Substituting the decomposition (\ref{q-bar-expansion}) into $\langle \tilde{L} q, q \rangle_{L^2(\R)}$ and using  the orthogonality relations of Propositions \ref{prop-orthog-1}, \ref{prop-orthog-2}, and \ref{prop-orthog-3}, we obtain the diagonalization formula (\ref{quadratic-form-L-tilde}).
\end{proof}

\subsection{Positivity of the second variation of $S(u)$}

To get the Lyapunov functional (\ref{Lyapunov}), we consider a linear superposition of the two Hamiltonian structures (\ref{BO-flow}) and (\ref{higher-BO-flow}) together. By Proposition \ref{lem-var}, the spatial profile $U$ is a critical point of 
$$
S(u) := R(u) + \gamma M(u) + \delta \left[ E(u) + c M(u) \right],  
$$
for every parameter $\delta \in \mathbb{R}$ to be defined. 
The second variation of this functional at $U$ is defined by the Hessian operator $\tilde{L} + \delta L : H^2(\R) \subset L^2(\R) \to L^2(\R)$. 
We show that the Hessian operator is positive semi-definite if $\delta = \frac{1}{2}(3c-k)$, which then yields (\ref{Lyapunov}) since 
$$
\gamma + \frac{1}{2} (3c-k) c = \frac{1}{4} (3c^2 -2kc - k^2) = \frac{1}{4} (3c + k)(c - k).
$$
By using (\ref{quadratic-form-L}) and (\ref{quadratic-form-L-tilde}), we obtain 
\begin{align}
\langle (\tilde{L} + \delta L) q, q \rangle_{L^2(\R)} &= 
2\pi \int_0^{\infty} (\lambda^2 + \delta \lambda + \gamma + \delta c) S_+(\lambda) d \lambda \notag \\
& \quad +  2\pi
\int_{\lambda_0}^{\mu_0} \frac{(\lambda - \lambda_0) (\lambda - \mu_0) (\lambda - c + \delta)}{(e^{\phi} - e^{-\phi}) \lambda}  S_-(\lambda) d \lambda,
\label{diag-final}
\end{align}
where 
\begin{align*}
S_+(\lambda) &= |a_I(\lambda)|^2 e^{2 \phi} \frac{\lambda - \lambda_0}{\lambda - \mu_0} + |a_{II}(\lambda)|^2 e^{-2 \phi} \frac{\lambda - \mu_0}{\lambda - \lambda_0} \geq 0, \quad \lambda \in (0,\infty),  \\
S_-(\lambda) &= |a_{III}(\lambda)|^2 e^{\phi} \frac{\lambda - \lambda_0}{\lambda - \mu_0} + |a_{IV}(\lambda)|^2 e^{-\phi} \frac{\lambda - \mu_0}{\lambda - \lambda_0} \leq 0, \quad \lambda \in (\lambda_0,\mu_0).
\end{align*}
The second integral in (\ref{diag-final}) is positive semi-definite if $\lambda - c + \delta < 0$ for all $\lambda \in (\lambda_0,\mu_0)$. This determines the range of admissible values of $\delta$ by $\delta \leq c - \mu_0$. Taking the maximal value 
of this range as $\delta := c - \mu_0$, we ensure that $\lambda - c + \delta = \lambda - \mu_0 < 0$ if $\lambda \in (\lambda_0,\mu_0)$. Since $\delta = -\lambda_0 - 2 \mu_0 > 0$, the first integral in (\ref{diag-final}) is positive semi-definite
for all $\lambda \in (0,\infty)$ if $\gamma + \delta c > 0$, which is true since 
$$
\gamma + \delta c = -(\lambda_0^2 + \lambda_0 \mu_0 + \mu_0^2) + (\lambda_0 + \mu_0) (\lambda_0 + 2 \mu_0) = \mu_0^2 + 2 \lambda_0 \mu_0 > 0.
$$
Hence, the self-adjoint operator $\tilde{L} + \delta L$ is positive semi-definite in $L^2(\R)$ and $S(u)$ with $\delta = c - \mu_0 = \frac{1}{2} (3c - k)$ is the Lyapunov functional for the periodic waves of the BO equation. 
The first assertion of Theorem \ref{th-main} is proven.

\section{Degeneracy of the second variation $S''(U)$ in Theorem \ref{th-main}}
\label{sec-6}

To prove the second assertion of Theorem \ref{th-main} for $N = 1$, we consider $S''(U) = \tilde{L} + \delta L$ with $\delta =  c-\mu_0$ for co-periodic perturbations in $L^2_{\rm per}$. We extend the orthogonality 
and completeness of eigenfunctions in $L^2_{\rm per}$ in order to  
characterize ${\rm Ker}_{L^2_{\rm per}} S''(U)$ and to show that 
it is spanned by the translational mode $U'$. Since 
the spectrum of $S''(U)$ in $L^2_{\rm per}$ consists of a sequence 
of isolated eigenvalues due to the compact embedding of $H^1_{\rm per}$ into $L^2_{\rm per}$, the non-degeneracy result implies coercivity 
of the quadratic form $\langle S''(U) v, v \rangle_{L^2_{\rm per}}$ for every $v \in H^1_{\rm per} |_{(U')^{\perp}}$ needed for the standard proof 
of Theorem \ref{th-orbit} with $N = 1$. For every fixed $N \in \mathbb{N}$, we extend the proof with the Floquet--Bloch theory in $L^2_{\rm Nper}$ 
to complete the proof of the second assertion of Theorem \ref{th-main} for $N \in \mathbb{N}$.

\subsection{Orthogonality of eigenfunctions in $L^2_{\rm per}$}

By \cite[Proposition 2.2]{GK-21}, the Lax equations (\ref{Lax-pair}) 
admit the purely discrete spectrum in $L^2_{\rm per}$ located at
\begin{equation}
\label{Lax-spectrum-discrete}
\sigma_L |_{L^2_{\rm per}} = \{ \lambda_0 \} \cup \{ 0\} \cup \{ k n \}_{n \in \mathbb{N}}.
\end{equation}
The eigenfunctions of the Lax equations (\ref{Lax-pair}) for 
$\lambda = kn$ with integer $n \geq 0$ are given by (\ref{eigenfunction-1}), whereas the eigenfunction for $\lambda = \lambda_0$ is obtained from the residue term of (\ref{eigenfunction-2}) given by 
$$
\lim_{\lambda \to \lambda_0} (\lambda - \lambda_0) \varphi_{II}^+(\lambda) = -\frac{k}{f'}, \quad \lim_{\lambda \to \lambda_0} (\lambda - \lambda_0) \varphi_{II}^-(\lambda) = -\frac{k}{f}.
$$
The first two eigenfunctions of the adjoint linearized BO equation (\ref{BO-adjoint}) are given by Proposition \ref{prop-eigenfunctions}:
\begin{align}
\begin{cases}
\psi_{I}(kn;x,t) = e^{i \lambda x + i \lambda^2 t} \frac{f(x-ct)}{f'(x-ct)} \left[ 1 + \frac{k}{(\lambda - \mu_0) f(x-ct)} \right],  \\
\psi_{II}(kn;x,t) = e^{-i \lambda x - i \lambda^2 t} \frac{f'(x-ct)}{f(x-ct)} \left[ 1 - \frac{k}{(\lambda - \lambda_0) f'(x-ct)} \right],
\end{cases}
\label{solution-discrete}
\end{align}
whereas first two eigenfunctions of the linearized BO equation (\ref{BO-lin}) are given by Proposition \ref{prop-eigenfunctions-adjoint}:
\begin{align}
\begin{cases}
q_{I}(kn;x,t) = e^{i \lambda x + i \lambda^2 t} \frac{f(x-ct)}{f'(x-ct)} \left[ 1 + \frac{k}{(\lambda - \mu_0) f'(x-ct)} \right], \\
q_{II}(kn;x,t) = e^{-i \lambda x - i \lambda^2 t} \frac{f'(x-ct)}{f(x-ct)} \left[ 1 - \frac{k}{(\lambda - \lambda_0) f(x-ct)} \right], 
\end{cases}
\label{solution-adjoint} 
\end{align}
where $\lambda = kn$ with integer $n \geq 0$. For $\lambda = 0$ ($n=0$), the two eigenfunctions in (\ref{solution-discrete}) are linearly dependent because 
\begin{equation}
\label{solutions-zero}
\psi_I(0) = \frac{f}{f'} - \frac{k}{\mu_0 f'} = \frac{\lambda_0}{\mu_0}, \quad 
\psi_{II}(0) = \frac{f'}{f} + \frac{k}{\lambda_0 f} = \frac{\mu_0}{\lambda_0},
\end{equation}
whereas the two eigenfunctions in (\ref{solution-adjoint}) are linearly independent because
\begin{equation}
\label{solutions-adjoint-zero}
q_I(0) = \frac{\lambda_0}{\mu_0} \left[ 1 - \frac{ik f'_x}{\lambda_0 \mu_0 (f')^2} \right], \quad q_{II}(0) = \frac{\mu_0}{\lambda_0} \left[ 1 - \frac{ik f_x}{\lambda_0 \mu_0 f^2} \right].
\end{equation}
However, we have the complex-conjugation symmetry between the two eigenfunctions:
\begin{equation}
\label{conjugation-zero}
\overline{q_{II}(0)} = \frac{\mu_0^2}{\lambda_0^2} q_I(0),
\end{equation}
which agrees with (\ref{conjugation-1}) for $\lambda = 0$ due to the relation (\ref{lambda-mu}). As a result, we can construct two real-valued solutions by using (\ref{solutions-adjoint-zero}):
\begin{align*}
\frac{1}{2} \left[ \frac{\mu_0}{\lambda_0} q_I(0) + \frac{\lambda_0}{\mu_0} q_{II}(0) \right] &= 1 + \frac{ik}{2 \lambda_0 \mu_0} \partial_x \left[ \frac{1}{f'} + \frac{1}{f} \right] \\
&= 1 + \frac{ik}{2 \lambda_0 \mu_0} \partial_x \frac{e^{-ik\xi} + \cosh \phi}{\cos(k\xi) + \cosh \phi} \\
&= 1 + \frac{k}{2\lambda_0 \mu_0} \partial_x \frac{\sin(k\xi)}{\cos(k\xi) + \cosh \phi} \\
&= 1 + \frac{k}{2\lambda_0 \mu_0} \partial_{\phi} U
\end{align*}
and 
\begin{align*}
\frac{1}{2i} \left[ \frac{\mu_0}{\lambda_0} q_I(0) - \frac{\lambda_0}{\mu_0} q_{II}(0) \right] &= \frac{k}{2 \lambda_0 \mu_0} \partial_x \left[ \frac{1}{f'} - \frac{1}{f} \right] \\
&= \frac{k}{2 \lambda_0 \mu_0} \partial_x \frac{\sinh \phi}{\cos(k\xi) + \cosh \phi} \\
&= \frac{1}{2\lambda_0 \mu_0} \partial_{x} U,
\end{align*}
where we have used (\ref{per-wave}), from which we get
\begin{equation}
\label{der-phi}
\partial_{\phi} U = \frac{k (1 + \cos(k\xi) \cosh \phi)}{[\cos(k\xi) + \cosh \phi]^2} = \partial_x \frac{\sin(k\xi)}{\cos(k\xi) + \cosh \phi}.
\end{equation}
Since $\partial_x U$ appears in the basis of eigenfunctions from solutions associated with the eigenvalue $\lambda = \lambda_0$, see (\ref{solutions-adjoint-extra}) below, we include only one solution 
to the basis of eigenfunctions from the two solutions (\ref{solutions-adjoint-zero}):
\begin{equation}
\label{q-0}
q_0(x,t) = 1 + \frac{k}{2\lambda_0 \mu_0} \partial_{\phi} U(x-ct).
\end{equation}
Similarly, we include only one solution 
to the basis of eigenfunctions from the two solutions (\ref{solutions-zero}):
\begin{equation}
\label{psi-0}
\psi_0(x,t) = 1.
\end{equation}
The following proposition establishes orthogonality between the eigenfunctions in $L^2_{\rm per}$, where we use the standard inner product,
$$
f,g \in L^2_{\rm per} : \quad \langle f, g \rangle_{L^2_{\rm per}} = 
\oint f(x) \bar{g}(x) dx.
$$

\begin{proposition}
	\label{prop-orthog-discrete}
	The set $\{ \bar{\psi}_I, \bar{\psi}_{II}\}$ is orthogonal to the set $\{ q_I, q_{II} \}$ in the sense:
	\begin{equation}
	\label{orth-2}
\begin{cases}
	\langle q_I(kn), \bar{\psi}_I(km) \rangle_{L^2_{\rm per}} =  
	\langle q_{II}(k n), \bar{\psi}_{II}(k m) \rangle_{L^2_{\rm per}} = 0, \\	
	\langle q_I(k n), \bar{\psi}_{II}(k m) \rangle_{L^2_{\rm per}} =  
	\langle q_{II}(k n), \bar{\psi}_{I}(k m) \rangle_{L^2_{\rm per}} = \frac{2\pi}{k} \delta_{n,m}, 
\end{cases} \quad n,m \in \mathbb{N}, 
	\end{equation}
	where $\delta_{n,m}$ is Kronecker's symbol. Furthermore, 
	\begin{equation}
	\label{orth-2-a}
\begin{cases}
\langle q_I(kn), \psi_0 \rangle_{L^2_{\rm per}} =  
\langle q_{II}(k n), \psi_0 \rangle_{L^2_{\rm per}} = 0, \\	
\langle q_0, \bar{\psi}_I(k n) \rangle_{L^2_{\rm per}} =  
\langle q_0, \bar{\psi}_{II}(kn) \rangle_{L^2_{\rm per}} = 0,
\end{cases} \quad n \in \mathbb{N}, 		
	\end{equation}
	and 
\begin{equation}
\label{orth-4}
\langle q_0, \psi_0 \rangle_{L^2_{\rm per}} = \frac{2\pi}{k}.
\end{equation}
\end{proposition}

\begin{proof}
	We follow the proof of Proposition \ref{prop-orthog-1}. 
	
	To prove the first line in (\ref{orth-2}), we use the representation 
	(\ref{representation-1}) for $\lambda = kn$ and $\lambda' = km$, where the function is analytic in $\mathbb{C}^+$ and the exponential factor $e^{i (\lambda + \lambda') x}$ decays to $0$ exponentially fast as ${\rm Im}(x) \to +\infty$. Integrating $q_I(\lambda) \psi_I(\lambda')$ along the boundary of the following rectangle
	$$
	\left\{ z \in \mathbb{C} : \quad |{\rm Re}(z)| \leq \frac{\pi}{k}, \;\; 
	0 \leq {\rm Im}(z) \leq R \right\}, \quad R > 0,
	$$
	in the positive orientation and taking the limit $R \to \infty$ yields 
	\begin{align*}
	\langle q_I(kn), \bar{\psi}_I(km) \rangle_{L^2_{\rm per}} = \lim_{R\to \infty} \oint  q_I(kn;x+iR,t) \psi_I(km;x+iR,t) dx = 0, 
	\end{align*}	
	due to the periodicity of the integrand for $n,m \in \mathbb{N}$. 
	The proof for $\langle q_{II}(\lambda), \bar{\psi}_{II}(\lambda') \rangle_{L^2_{\rm per}} = 0$ is similar but the contour of integration is closed in $\mathbb{C}^-$. 
	
	To prove the second line in (\ref{orth-2}), we use the representation 
	(\ref{representation-2}). Since 
$$
\oint e^{i k (n-m) x + i k^2 (n^2 - m^2)t} dx = \frac{2\pi}{k} \delta_{n,m}
$$ 
we get $\langle q_I(kn), \bar{\psi}_{II}(km) \rangle_{L^2_{\rm per}} = \frac{2\pi}{k} \delta_{n,m}$ if and only if the integral of all remaining terms after expanding the product vanishes:
	\begin{align*}
0 &= \frac{ik}{(kn - \mu_0) (km - \lambda_0)} \oint \frac{d}{dx} \left( \frac{e^{i k (n-m) x + i k^2 (n^2-m^2) t} }{f'} \right) dx,
	\end{align*}	
which is true. The result for $\langle q_{II}(kn), \bar{\psi}_{I}(km) \rangle_{L^2_{\rm per}} = \frac{2\pi}{k} \delta_{n,m}$
holds by the complex conjugate symmetry (\ref{conjugation-1}) and (\ref{conjugation-1-q}).

The orthogonality relations (\ref{orth-2-a}) are obtained by setting either $n \neq 0$ and $m = 0$ or $n = 0$ and $m \neq 0$ in (\ref{orth-2}). The orthogonality relations (\ref{orth-4}) are obtained directly from  (\ref{q-0}) and (\ref{psi-0}), where we use $\langle \partial_{\phi} U, 1 \rangle_{L^2_{\rm per}} = 0$, which follows from the fact that 
\begin{align*}
\langle  U, 1 \rangle_{L^2_{\rm per}} &= k \sinh \phi \oint \frac{dx}{\cos(kx) + \cosh \phi} \\
&= \sinh \phi \int_0^{2\pi} \frac{dx}{\cos x + \cosh \phi} = 2 \pi
\end{align*}
is independent of $\phi > 0$. 
\end{proof}

It remains to augment the set of eigenfunctions (\ref{solution-discrete}) and (\ref{solution-adjoint}) for $\lambda = kn$ with integer $n \geq 0$ 
by the set of eigenfunctions for $\lambda = \lambda_0$, according to the Lax spectrum (\ref{Lax-spectrum-discrete}). It follows from (\ref{solution-3}) and (\ref{solution-4}) that
$$
\psi_{III}(\lambda_0) = \frac{f-1}{f f'} = -\frac{\lambda_0}{k^2} U, \quad 
\psi_{IV}(\mu_0) = \frac{f'-1}{f f'} = -\frac{\mu_0}{k^2} U,
$$
where $U$ is the periodic profile (\ref{per-wave}) and we have used (\ref{lambda-mu}). Thus, both solutions are constant multiples of $U$. This is natural due to the translational symmetry, which implies that the linearized BO equation (\ref{BO-lin}) with $u(x,t) = U(x-ct)$ admits the solution $\partial_x U$ and the adjoint linearized BO equation (\ref{BO-adjoint}) admits the solution $U$.

To get the linearly independent solutions, we compute derivatives of $\psi_{III}(\lambda)$ at $\lambda = \lambda_0$ and $\psi_{IV}(\lambda)$ at $\lambda = \mu_0$:
\begin{align*}
\psi_{III}'(\lambda_0) = i (x + 2 \lambda_0 t) \frac{f-1}{f f'} - \frac{1}{k f f'}, \quad
\psi_{IV}'(\mu_0) = -i (x + 2 \mu_0 t) \frac{f'-1}{f f'} + \frac{1}{k f f'}.
\end{align*}
To get the $\frac{2\pi}{k}$-periodic solution, we obtain the linear combination of the two derivatives:
\begin{align*}
\mu_0 \psi_{III}'(\lambda_0) + \lambda_0 \psi_{IV}'(\mu_0) &= \frac{2 i \lambda_0 \mu_0}{k}  t U - \frac{1}{ff'} \\
&= \frac{2 i \lambda_0 \mu_0 t \sinh \phi}{\cos(k \xi) + \cosh \phi} - \frac{1}{2} \frac{e^{-ik\xi}}{\cos(k \xi) + \cosh \phi}. 
\end{align*}
The real part of this function is given by a linear combination of $1$ and $U$ and the imaginary part of this function is a linear combination of $\partial_x^{-1} \partial_{\phi} U$ and $t U$, where the uniquely defined 
real-valued odd function $\partial_x^{-1} \partial_{\phi} U$ is obtained from (\ref{der-phi}) as 
$$
\partial_x^{-1} \partial_{\phi} U(x) := \frac{\sin(k x)}{\cos(k x) + \cosh \phi}.
$$
The adjoint linearized BO equation (\ref{BO-adjoint}) in $L^2_{\rm per}$ admits two linearly independent eigenfunctions
\begin{equation}
\label{solutions-extra}
\psi_1(x,t) = U(x-ct), \quad \psi_2(x,t) = \partial_x^{-1} \partial_{\phi} U(x-ct) + \frac{k}{\sinh^2 \phi} t U(x-ct).
\end{equation}
The linearized BO equation (\ref{BO-lin}) admits two linearly independent eigenfunctions 
\begin{equation}
\label{solutions-adjoint-extra}
q_1(x,t) = \partial_x U(x-ct), \quad q_2(x,t) = \partial_{\phi} U(x-ct) + \frac{k}{\sinh^2 \phi} t \partial_x U(x-ct),
\end{equation}
where we have used (\ref{lambda-mu}) to obtain 
\begin{equation}
\label{tech-1}
\frac{4 \lambda_0 \mu_0}{k} = \frac{k}{\sinh^2 \phi} = -\frac{dc}{d\phi}.
\end{equation}
We have the following orthogonality result between the eigenfunctions.

\begin{proposition}
	\label{prop-orth-additional}
	Solutions (\ref{solutions-extra}) and (\ref{solutions-adjoint-extra}) satisfy the following orthogonality conditions:
	\begin{equation}
	\label{orth-5}
\begin{cases}
	\langle q_1, \psi_1 \rangle_{L^2_{\rm per}} = \langle q_2, \psi_2 \rangle_{L^2_{\rm per}} = 0, \\	
	\langle q_1, \psi_2 \rangle_{L^2_{\rm per}} =  
	-\langle q_2, \psi_1 \rangle_{L^2_{\rm per}} = \frac{\pi k}{\sinh^2 \phi}.
\end{cases}
	\end{equation}
	Furthermore, we have 
		\begin{equation}
	\label{orth-5-a}
\begin{cases}
	\langle q_I(kn), \psi_1 \rangle_{L^2_{\rm per}} =  
	\langle q_{II}(k n), \psi_1 \rangle_{L^2_{\rm per}} = 0, \\	
	\langle q_I(kn), \psi_2 \rangle_{L^2_{\rm per}} =  
\langle q_{II}(k n), \psi_2 \rangle_{L^2_{\rm per}} = 0, \\ 
	\langle q_1, \bar{\psi}_I(kn) \rangle_{L^2_{\rm per}} =  
\langle q_1, \bar{\psi}_{II}(kn) \rangle_{L^2_{\rm per}} = 0, \\	
\langle q_2, \bar{\psi}_I(kn) \rangle_{L^2_{\rm per}} =  
\langle q_2, \bar{\psi}_{II}(kn) \rangle_{L^2_{\rm per}} = 0, 
\end{cases}	
\quad n \in \mathbb{N},
	\end{equation}
	and 
	\begin{equation}
	\label{orth-5-b}
\begin{cases}
	\langle q_1, \psi_0 \rangle_{L^2_{\rm per}} =  
\langle q_2, \psi_0 \rangle_{L^2_{\rm per}} = 0, \\	
	\langle q_0, \psi_1 \rangle_{L^2_{\rm per}} =  
\langle q_0, \psi_2 \rangle_{L^2_{\rm per}} = 0.
\end{cases}
	\end{equation}
\end{proposition}

\begin{proof}
	The first line in (\ref{orth-5}) follows from the opposite parity between $U$ and $\partial_x U$ as well as $\partial_x^{-1} \partial_{\phi} U$ and $\partial_{\phi} U$ integrated over the periodic domain together with 
	$$
	\langle \partial_x U, \partial_x^{-1} \partial_{\phi} U \rangle_{L^2_{\rm per}} + 	\langle U, \partial_{\phi} U \rangle_{L^2_{\rm per}} = 0.
	$$	
	The second line in (\ref{orth-5}) is obtained after integration by parts in the form:
$$
\langle q_1, \psi_2 \rangle_{L^2_{\rm per}} =  
-\langle q_2, \psi_1 \rangle_{L^2_{\rm per}} = -\frac{1}{2} \frac{d}{d\phi} \| U \|^2_{L^2_{\rm per}}.
$$	
The numerical coefficient $\frac{d}{d\phi} \| U \|^2_{L^2_{\rm per}}$ is well-known in the theory of stability of periodic waves in the BO equation and can be computed explicitly as follows. Since 
	$$
	\frac{d}{dx} \frac{\sin(kx)}{\cos(kx) + \cosh \phi} = k \frac{1+ \cosh \phi \cos(kx)}{(\cos(kx) + \cosh \phi)^2},
	$$
	we obtain 
	\begin{align*}
	\| U \|^2_{L^2_{\rm per}} &= k^2 \sinh^2 \phi \oint \frac{dx}{(\cos(kx) + \cosh \phi)^2} \\
	&= -k^2 \sinh^2 \phi \cosh \phi \oint \frac{\cos(kx) dx}{(\cos(kx) + \cosh \phi)^2} \\
	&=-k^2 \sinh^2 \phi \cosh \phi \oint \frac{dx}{\cos(kx) + \cosh \phi} 
	+ k^2 \sinh^2 \phi \cosh^2 \phi \oint \frac{dx}{(\cos(kx) + \cosh \phi)^2}.
	\end{align*}
Eliminating the integral of $(\cos(kx) + \cosh \phi)^{-2}$ yields
	\begin{align*}
	\| U \|^2_{L^2_{\rm per}} &= k^2 \cosh \phi \oint \frac{dx}{\cos(kx) + \cosh \phi} \\
	&= k \cosh \phi \int_0^{2\pi} \frac{dx}{\cos x + \cosh \phi} \\
	&= 2 \pi k \coth \phi,
	\end{align*}
so that 
\begin{equation}
\label{tech-2}
	-\frac{1}{2} \frac{d}{d\phi} \| U \|^2_{L^2_{\rm per}} = \frac{\pi k}{\sinh^2 \phi} = - \pi \frac{dc}{d\phi},
\end{equation}
	which completes the proof of the second line in (\ref{orth-5}).
	
	The orthogonality conditions (\ref{orth-5-a}) follow from the orthogonality conditions (\ref{orth-3-b}) and (\ref{orth-3-c}) in Proposition \ref{prop-orthog-2} for $\lambda = kn$ with $n \in \mathbb{N}$. To verify the orthogonality conditions (\ref{orth-5-b}), we use (\ref{q-0}), (\ref{psi-0}), (\ref{solutions-extra}), and (\ref{solutions-adjoint-extra}) directly. 
	For the first line in (\ref{orth-5-b}), we obtain 
\begin{align*}
\begin{cases}
	\langle q_1, \psi_0 \rangle_{L^2_{\rm per}} =  \langle \partial_x U, 1 \rangle_{L^2_{\rm per}} = 0, \\
	\langle q_2, \psi_0 \rangle_{L^2_{\rm per}} = \langle \partial_{\phi} U, 1 \rangle_{L^2_{\rm per}} = 0,
	\end{cases}		
\end{align*}
where we have used that $\langle U, 1 \rangle_{L^2_{\rm per}} = 2\pi$ is independent of $\phi > 0$. For the second line in (\ref{orth-5-b}), we get 
\begin{align*}
\langle q_0, \psi_1 \rangle_{L^2_{\rm per}} & = \langle 1, U \rangle_{L^2_{\rm per}} + \frac{k}{2\lambda_0 \mu_0} \langle \partial_{\phi} U, U \rangle_{L^2_{\rm per}} \\
& = 2\pi + \frac{k}{4\lambda_0 \mu_0} \frac{d}{d\phi} \| U \|^2_{L^2_{\rm per}} = 0,
\end{align*}
where we have used (\ref{tech-1}) and (\ref{tech-2}), and 
\begin{align*}
\langle q_0, \psi_2 \rangle_{L^2_{\rm per}} & = \langle 1, \partial_x^{-1} \partial_{\phi} U \rangle_{L^2_{\rm per}} + \frac{k}{2\lambda_0 \mu_0} \langle \partial_{\phi} U, \partial_x^{-1} \partial_{\phi} U  \rangle_{L^2_{\rm per}} = 0,
\end{align*}
where we have used the opposite parity of $1$ and $\partial_{\phi} U$ relative to $\partial_x^{-1} \partial_{\phi} U$.
\end{proof}

\subsection{Completeness of eigenfunctions in $L^2_{\rm per}$}

To simplify the proof of completeness, we set $t = 0$ again and suppress the time variable from the list of arguments. We claim that solutions (\ref{solution-adjoint}) and (\ref{solutions-adjoint-extra}) provide a basis in $L^2_{\rm per}$ in the sense that every element $q \in L^2_{\rm per}$ can be uniquely represented by the superposition of elements of the basis:
\begin{equation}
\label{basis-1}
q(x) = \sum_{n \in \mathbb{N}} a_n q_I(kn;x) + b_n q_{II}(kn;x) + c_0 q_0(x) + c_1 \partial_x U(x) + c_2  \partial_{\phi} U(x),
\end{equation}
where the coefficients in the decomposition (\ref{basis-1}) are computed by using the orthogonality relations of Propositions \ref{prop-orthog-discrete} and \ref{prop-orth-additional}:
\begin{align*}
a_n = \frac{k}{2\pi} \langle q, \bar{\psi}_{II}(kn) \rangle_{L^2_{\rm per}}, \quad
b_n = \frac{k}{2\pi} \langle q, \bar{\psi}_{I}(kn) \rangle_{L^2_{\rm per}}, \quad n \in \mathbb{N},
\end{align*}
\begin{align*}
c_0 = \frac{k}{2\pi} \langle q, \psi_0 \rangle_{L^2_{\rm per}},
\end{align*}
and 
\begin{align*}
c_1 = \frac{\sinh^2 \phi}{\pi k}  \langle q, \psi_2 \rangle_{L^2_{\rm per}}, \quad
c_2 = -\frac{\sinh^2 \phi}{\pi k} \langle q, \psi_1 \rangle_{L^2_{\rm per}}.
\end{align*}
Substituting all projections into (\ref{basis-1}) yields the identity
\begin{align*}
q(x) &= \frac{k}{2\pi} \oint q(y) \left( \sum_{n \in \mathbb{N}} q_{I}(kn;x) \psi_{II}(kn;y) + q_{II}(kn;x) \psi_{I}(kn;y) \right. \\
& \qquad \left. + \frac{1}{2} \left[ q_{I}(0;x) \psi_{II}(0;y) + q_{II}(0;x) \psi_{I}(0;y)\right]   + \frac{2}{k^2} \sinh^2 \phi \left[ q_1(x) \psi_2(y) - q_2(x) \psi_1(y)\right] \right) dy,
\end{align*}
where we have used
$$
q_0(x) \psi_0 = \frac{1}{2} \left[ \frac{\mu_0}{\lambda_0} q_I(0;x) + \frac{\lambda_0}{\mu_0} q_{II}(0;x) \right] = \frac{1}{2} \left[ q_{I}(0;x) \psi_{II}(0;y) + q_{II}(0;x) \psi_{I}(0;y)\right],
$$
which follows from (\ref{solutions-zero}) and (\ref{solutions-adjoint-zero}). The identity holds for every $f \in L^2_{\rm per}$ if and only if 
the eigenfunctions satisfy the completeness relation 
\begin{align}
& \sum_{n \in \mathbb{N}} q_I(kn;x) \psi_{II}(kn;y) +  q_{II}(kn;x) \psi_I(kn;y) + \frac{1}{2} (q_I(0;x) \psi_{II}(0;y) + q_{II}(0;x) \psi_I(0;y)) \notag \\
& \quad  + \frac{2}{k^2} \sinh^2 \phi (q_1(x) \psi_2(y) - q_2(x) \psi_1(y)) = \frac{2\pi}{k} \delta(x-y),
\label{complete}
\end{align}
where the Dirac delta distribution $\delta(x-y)$ in the periodic domain is understood in the sense of the generalized identity:
\begin{equation}
\label{Dirac-periodic}
\sum_{n \in \mathbb{Z}} e^{i k n (x-y)} = \frac{2\pi}{k} \delta(x-y).
\end{equation}
The following proposition confirms the validity of the completeness relation (\ref{complete}). 

\begin{proposition}
	\label{prop-complete}
	The completeness relation (\ref{complete}) is satisfied. 
\end{proposition}

\begin{proof}
	We follow the proof of Proposition \ref{prop-complete-2} and compute each term in (\ref{complete}). The first and second terms are given by the same expressions (\ref{term-first}) and (\ref{term-second}) with $\lambda = kn$ for integer $n \geq 0$. Terms in the last brackets in (\ref{term-first}) and (\ref{term-second}) produce different outcomes depending on the summation in $n$ from $1$ to $\infty$ or from $0$ to $\infty$. By using $\mu_0 = \lambda_0 + k$, we obtain 
	\begin{align}
	\sum_{n=1}^{\infty} \left[ \frac{e^{i k (n+1)x - i k n y}}{k (n+1) - \mu_0} - \frac{e^{i k n x - i k (n-1) y}}{k n - \mu_0} \right] &= \frac{e^{ikx}}{\lambda_0}, \label{sum-1} \\
	\sum_{n=0}^{\infty} \left[ \frac{e^{i k (n+1)x - i k n y}}{k (n+1) - \mu_0} - \frac{e^{i k n x - i k (n-1) y}}{k n - \mu_0} \right] &= \frac{e^{iky}}{\mu_0}, 
	\label{sum-2} \\
	\sum_{n=1}^{\infty} \left[ \frac{e^{-i k n x + i k (n+1) y}}{k (n+1) - \mu_0} - \frac{e^{-i k (n-1) x + i k n y}}{k n - \mu_0} \right] &= \frac{e^{iky}}{\lambda_0}, \label{sum-3} \\
	\sum_{n=0}^{\infty} \left[ \frac{e^{-i k n x + i k (n+1) y}}{k (n+1) - \mu_0} -
	\frac{e^{-i k (n-1) x + i k n y}}{k n - \mu_0} \right] &= \frac{e^{ikx}}{\mu_0}.
	\label{sum-4} 
	\end{align}
	Splitting symmetrically the series in (\ref{complete}) and using (\ref{sum-1}), (\ref{sum-2}), (\ref{sum-3}), and (\ref{sum-4}), we obtain 
	\begin{align*}
	& \sum_{n \in \mathbb{N}} q_I(kn;x) \psi_{II}(kn;y) +  q_{II}(kn;x) \psi_I(kn;y) + \frac{1}{2} (q_I(0;x) \psi_{II}(0;y) + q_{II}(0;x) \psi_I(0;y)) \\
	&= \frac{1}{2} \sum_{n=0}^{\infty} q_I(kn;x) \psi_{II}(kn;y) +  q_{II}(kn;x) \psi_I(kn;y) \\
	& \qquad + \frac{1}{2} \sum_{n=1}^{\infty} q_I(kn;x) \psi_{II}(kn;y) +  q_{II}(kn;x) \psi_I(kn;y) \\
	&= \sum_{n \in \mathbb{Z}} e^{i k n (x-y)} + \frac{i}{2k} (e^{\phi} - e^{-\phi})^2 \frac{\partial}{\partial x} \left[ \frac{e^{ikx - \phi}}{f'(x) f(y)} + \frac{e^{iky + \phi}}{f'(x) f(y)} - \frac{e^{iky-\phi}}{f(x) f'(y)} - \frac{e^{ikx + \phi}}{f(x) f'(y)} \right], 
	\end{align*}
	where we have used (\ref{lambda-mu}) again. Since the first term gives $\frac{2\pi}{k} \delta(x-y)$ due to (\ref{Dirac-periodic}), the completeness relation (\ref{complete}) is satisfied if and only if the remaining terms
$$
\frac{2i}{k} \sinh^2 \phi \frac{\partial}{\partial x} \left[ \frac{e^{ikx - \phi}}{f'(x) f(y)} + \frac{e^{iky + \phi}}{f'(x) f(y)} - \frac{e^{iky-\phi}}{f(x) f'(y)} - \frac{e^{ikx + \phi}}{f(x) f'(y)} \right]
$$
cancel out with the last two terms in the left-hand-side of (\ref{complete}):
$$
\frac{2}{k^2} \sinh^2 \phi (q_1(x) \psi_2(y) - q_2(x) \psi_1(y)).
$$
We note that 
	\begin{align*}
	\frac{e^{iky + \phi}}{f'(x) f(y)} - \frac{e^{iky-\phi}}{f(x) f'(y)} 
	&= \frac{e^{i k y} (e^{\phi}-e^{-\phi}) (1 + e^{ikx} (e^{\phi} + e^{-\phi}) + e^{ik(x+y)})}{f(x) f'(x) f(y) f'(y)}, \\ 
	\frac{e^{ikx -\phi}}{f'(x) f(y)} - \frac{e^{ikx+\phi}}{f(x) f'(y)} 
	&= -\frac{e^{i k x} (e^{\phi}-e^{-\phi}) (1 + e^{iky} (e^{\phi} + e^{-\phi}) + e^{ik(x+y)})}{f(x) f'(x) f(y) f'(y)}, 
	\end{align*}
	hence, we get 
	\begin{align}
& \frac{2i}{k} \sinh^2 \phi \frac{\partial}{\partial x} \left[ \frac{e^{ikx - \phi}}{f'(x) f(y)} + \frac{e^{iky + \phi}}{f'(x) f(y)} - \frac{e^{iky-\phi}}{f(x) f'(y)} - \frac{e^{ikx + \phi}}{f(x) f'(y)} \right] \notag \\
	& \quad = \frac{4i}{k} \sinh^3 \phi 
	\frac{\partial}{\partial x} \frac{(e^{i k y} - e^{ikx}) (1 + e^{ik(x+y)})}{f(x) f'(x) f(y) f'(y)}. \label{help-1}
	\end{align}
	On the other hand, we expand 
	\begin{align}
	& \frac{2}{k^2} \sinh^2 \phi (q_1(x) \psi_2(y) - q_2(x) \psi_1(y)) \notag \\
	&= 
\frac{2}{k^2} \sinh^2 \phi \frac{\partial}{\partial x} \left[ U(x) \partial_y^{-1} \partial_{\phi} U(y) - U(y) \partial_x^{-1} \partial_{\phi} U(x) \right]  \notag \\
	&= \frac{2}{k^2} \sinh^2 \phi \frac{\partial}{\partial x} \left[ 
	\frac{(f'_x f - f_x f') (e^{2iky}-1)}{f(x) f'(x) f(y) f'(y)} - \frac{(f'_y f - f_y f') (e^{2ikx}-1)}{f(x) f'(x) f(y) f'(y)}  \right] \notag \\
	&= -\frac{4i}{k} \sinh^3 \phi \frac{\partial}{\partial x} 
	\frac{e^{ikx} (e^{2iky}-1) - e^{iky} (e^{2ikx}-1)}{f(x) f'(x) f(y) f'(y)} \notag \\
	&= -\frac{4i}{k} \sinh^3 \phi \frac{\partial}{\partial x} 
	\frac{(e^{iky}-e^{ikx}) (1 + e^{ik(x+y)})}{f(x) f'(x) f(y) f'(y)}.
	\label{help-2}
	\end{align}
The sum of (\ref{help-1}) and (\ref{help-2}) is zero so that the completeness relation (\ref{complete}) is proven. 
\end{proof}

\subsection{Diagonalization of $L$}

For the set of eigenfunctions (\ref{solution-adjoint}), we introduce the same expressions $Q_{I,n}$ and $Q_{II,n}$ as in (\ref{eigenfunction-explicit-1}) but for $\lambda = k n$ with integer $n \geq 0$. These expressions correspond to the Floquet--Bloch form (\ref{eigenfunction-Floquet-1}) with $\theta = 0$.
The eigenfunctions satisfy the spectral problem (\ref{spec-BO}) in the periodic domain with the eigenvalues 
\begin{align*}
\begin{cases}
\Lambda_{I,n} = i kn (kn + c), \\
\Lambda_{II,n} = -ikn (kn + c), 
\end{cases} \quad n \geq 0.
\end{align*}
For the set of eigenfunctions (\ref{solutions-adjoint-extra}), we have 
\begin{equation}
\label{3-4-eigen}
L \partial_x U = 0, \quad  L \partial_{\phi} U = -\frac{dc}{d\phi} U,
\end{equation}
due to the differentiation of the stationary BO equation $-H(U') - U^2 + c U =0$ in $x$ and $\phi$.

We can write explicitly similarly to (\ref{operators-explicit}) that 
\begin{equation}
\label{1-2-eigen}
L Q_{I,n} = (kn + c) \Psi_{I,n}, \quad L Q_{II,n} = (kn + c) \Psi_{II,n}, 
\end{equation}
which holds for any integer $n \geq 1$. The relations (\ref{1-2-eigen}) are extended to eigenfunctions (\ref{q-0}) and (\ref{psi-0}) for $n = 0$ since 
$$
L q_0 = c - 2 U + \frac{k}{2\lambda_0 \mu_0} L \partial_{\phi} U = 
c - 2 U - \frac{k}{2\lambda_0 \mu_0} \frac{dc}{d\phi} U = c \psi_0,
$$
where we have used (\ref{tech-1}) and (\ref{3-4-eigen}).

By using the decomposition (\ref{basis-1}) and the completeness result of Proposition \ref{prop-complete}, we obtain for every $q \in H^{1/2}_{\rm per}$, similarly to Proposition \ref{prop-diagoalization-L}, that 
	\begin{align}
\langle L q, q \rangle_{L^2_{\rm per}} &= 
\frac{2\pi}{k} \sum_{n\in \mathbb{N}} (kn + c) \left[ |a_n|^2 e^{2 \phi} \frac{kn - \lambda_0}{kn - \mu_0} + |b_n|^2 e^{-2 \phi} \frac{kn - \mu_0}{kn - \lambda_0} \right] \notag \\
& \quad + \frac{2\pi c}{k} |c_0|^2 - \frac{\pi k^2}{\sinh^4 \phi} |c_2|^2, 
\label{quadratic-form-L-discrete}
\end{align}
which follows from 
\begin{align*}
L q &= \sum_{n\in \mathbb{N}} [a_n L Q_{I,n} + b_n L Q_{II,n}] + c_0 L q_0 
+ c_1 L \partial_x U + c_2 L \partial_{\phi} U \\
&= \sum_{n\in \mathbb{N}} (kn + c) [a_n \Psi_{I,n} + b_n \Psi_{II,n}] + c_0 c \psi_0 - c_2 \frac{dc}{d\phi} U
\end{align*}
and
\begin{align*}
\bar{q} &= \sum_{n\in \mathbb{N}} \left[\bar{a}_n e^{2 \phi} \frac{kn - \lambda_0}{kn - \mu_0} Q_{II,n} + \bar{b}_n e^{-2 \phi} \frac{kn - \mu_0}{kn- \lambda_0} Q_{I,n} \right] + \bar{c}_0 q_0 + \bar{c}_1 \partial_x U + \bar{c}_2 \partial_{\phi} U,
\end{align*}
due to the orthogonality result of Propositions \ref{prop-orthog-discrete} and \ref{prop-orth-additional}. Similarly to the outcome of Proposition \ref{prop-diagoalization-L}, all terms but the last one in (\ref{quadratic-form-L-discrete}) are positive semi-definite, whereas the last term in (\ref{quadratic-form-L-discrete}) is negative semi-definite. 
The coefficient $c_1$ does not appear in (\ref{quadratic-form-L-discrete}).

\subsection{Diagonalization of $\tilde{L}$}

The eigenfunctions $Q_{I,n}$ and $Q_{II,n}$ satisfy the spectral problem (\ref{higher-spec-BO}) in the periodic domain for the eigenvalues 
\begin{align*}
\begin{cases}
\tilde{\Lambda}_{I,n} = i kn (k^2 n^2 + \gamma), \\
\tilde{\Lambda}_{II,n} = -ikn (k^2 n^2 + \gamma), 
\end{cases} \quad n \geq 0.
\end{align*}
The eigenfunctions (\ref{solutions-adjoint-extra}) satisfy
\begin{equation}
\label{3-4-eigen-tilde}
\tilde{L} \partial_x U = 0, \quad  \tilde{L} \partial_{\phi} U = -\frac{d\gamma}{d\phi} U,
\end{equation}
due to the differentiation of the stationary higher-order BO equation 
$$
-U'' + \frac{3}{2} U H(U') + \frac{3}{2} H(U U') + U^3 + \gamma U = 0
$$ 
in $x$ and $\phi$.

We can write explicitly similarly to (\ref{operators-explicit-tilde}) that 
\begin{equation}
\label{1-2-eigen-tilde}
\tilde{L} Q_{I,n} = (k^2 n^2 + \gamma) \Psi_{I,n}, \quad \tilde{L} Q_{II,n} = (k^2 n^2 + \gamma) \Psi_{II,n},
\end{equation}
which holds for any integer $n \geq 1$. The relations (\ref{1-2-eigen-tilde}) are extended to eigenfunctions (\ref{q-0}) and (\ref{psi-0}) for $n = 0$ since 
$$
\tilde{L} q_0 = \gamma + 3 U^2 + 3 H(U') + \frac{k}{2\lambda_0 \mu_0} \tilde{L}\partial_{\phi} U = 
\gamma + 3 cU - \frac{k}{2\lambda_0 \mu_0} \frac{d\gamma}{d\phi} U = \gamma \psi_0,
$$
where we have used the stationary BO equation $-H(U') - U^2 + c U =0$, as well as relations (\ref{tech-1}), (\ref{3-4-eigen-tilde}), and $\frac{d\gamma}{d\phi} = -\frac{3}{2} c \frac{dc}{d\phi}$.

By using the decomposition (\ref{basis-1}) and the completeness result of Proposition \ref{prop-complete}, we obtain for every $q \in H^1_{\rm per}$, similarly to Proposition \ref{prop-diagoalization-L-tilde}, that 
\begin{align}
\langle \tilde{L} q, q \rangle_{L^2_{\rm per}} &= 
\frac{2\pi}{k} \sum_{n\in \mathbb{N}} (k^2 n^2 + \gamma) \left[ |a_n|^2 e^{2 \phi} \frac{kn - \lambda_0}{kn - \mu_0} + |b_n|^2 e^{-2 \phi} \frac{kn - \mu_0}{kn - \lambda_0} \right] \notag \\
& \quad + \frac{2\pi \gamma}{k} |c_0|^2 + \frac{3 \pi k^2 c}{2 \sinh^4 \phi} |c_2|^2,
\label{quadratic-form-L-discrete-tilde}
\end{align}
where we have used (\ref{tech-1}), (\ref{orth-5}), and $\frac{d\gamma}{d\phi} = -\frac{3}{2} c \frac{dc}{d\phi}$. 
Similarly to the outcome of Proposition \ref{prop-diagoalization-L-tilde}, the last term in (\ref{quadratic-form-L-discrete-tilde}) is positive semi-definite, whereas the other terms in (\ref{quadratic-form-L-discrete-tilde}) are sign-indefinite. The coefficient $c_1$ does not appear in (\ref{quadratic-form-L-discrete-tilde}).

\subsection{Diagonalization of $\tilde{L} + \delta L$ in $L^2_{\rm per}$}

We define $\delta = c - \mu_0 = -\lambda_0 - 2 \mu_0 > 0$ and obtain 
the following diagonalization of $\tilde{L} + \delta L$ for every $q \in H^1_{\rm per}$ from (\ref{quadratic-form-L-discrete}) and (\ref{quadratic-form-L-discrete-tilde}):
\begin{align}
\langle (\tilde{L} + \delta L) q, q \rangle_{L^2_{\rm per}} &= 
\frac{2\pi}{k} \sum_{n\in \mathbb{N}} (k^2 n^2 + \delta kn + \gamma + \delta c) S_n \notag \\
& \quad + \frac{2\pi}{k} (\gamma + \delta c) |c_0|^2 + \frac{\pi k^2 }{\sinh^4 \phi} \left( \frac{3}{2} c - \delta \right) |c_2|^2,
\label{diagonalization-final}
\end{align}
where 
\begin{align*}
S_n &= |a_n|^2 e^{2 \phi} \frac{kn - \lambda_0}{kn - \mu_0} + |b_n|^2 e^{-2 \phi} \frac{kn - \mu_0}{kn - \lambda_0} \geq 0. 
\end{align*}
Since $\delta > 0$, $\gamma + \delta c = \mu_0^2 + 2 \lambda_0 \mu_0 > 0$, 
and 
$$
\frac{3}{2} c - \delta = \frac{c}{2} + \mu_0 = \frac{k}{2} > 0,
$$
all terms in (\ref{diagonalization-final}) are positive semi-definite and the only term missing is the projection $c_1$ of $q \in H^1_{\rm per}$ to $\partial_x U$. Hence, ${\rm Ker}_{L^2_{\rm per}} (\tilde{L} + \delta L) = {\rm span}(U')$, where $U'$ is the translational mode of the linearized BO equation (\ref{BO-lin}) with $u(x,t) = U(x-ct)$. Coercivity 
of the quadratic form $\langle (\tilde{L} + \delta L) q, q \rangle_{L^2_{\rm per}}$ for every $q \in H^1_{\rm per} |_{(U')^{\perp}}$ 
follows from positivity of all other coefficients in  (\ref{diagonalization-final}) with the explicit lower bound for some $C > 0$:
\begin{align*}
\langle (\tilde{L} + \delta L) q, q \rangle_{L^2_{\rm per}} &\geq 
\frac{2\pi}{k} \frac{\mu_0^2}{\lambda_0^2} \sum_{n\in \mathbb{N}} (k^2 n^2 + \gamma + \delta c) (|a_n|^2 + |b_n|^2)  \notag \\
& \quad + \frac{2\pi}{k} (\gamma + \delta c) |c_0|^2 + \frac{\pi k^2 }{\sinh^4 \phi} \left( \frac{3}{2} c - \delta \right) |c_2|^2 \\
&\geq C \| q - c_1 \partial_x U \|_{H^1_{\rm per}}^2,
\end{align*}
where we have used that 
$$
\frac{kn-\lambda_0}{kn - \mu_0} \geq 1, \quad \frac{kn-\mu_0}{kn - \lambda_0} \geq \frac{\mu_0}{\lambda_0}, \quad e^{2\phi} = \frac{\lambda_0}{\mu_0} > 1.
$$
The proof of Theorem \ref{th-main} is complete for $N = 1$.

\subsection{Diagonalization of $\tilde{L} + \delta L$ in $L^2_{\rm N per}$ with fixed $N \in \mathbb{N}$} 

Working in $L^2_{\rm N per}$ for any fixed $N \in \mathbb{N}$, we can use the Floquet--Bloch form of eigenfunctions (\ref{eigenfunction-Floquet-1}) and (\ref{eigenfunction-Floquet-2}) for the following $N$ values of parameter $\theta \in [0,k)$:
$$
\theta_m := \frac{m k}{N}, \quad m \in \mathcal{M}_N := \{ 0,1,\dots,N-1\} .
$$
Eigenfunctions $\{ Q_{I,n}(\theta_m)\}_{n \geq 0, m \in \mathcal{M}_N}$ and $\{ Q_{II,n}(\theta_m) \}_{n \geq 0, m \in \mathcal{M}_N}$ are defined by (\ref{eigenfunction-Floquet-1}) for $\theta = \theta_m$ with the only exception for $(n,m) = (0,0)$, for which the proper eigenfunction is $q_0$ in (\ref{q-0}). Eigenfunctions $\{ Q_{III}(\theta_m) \}_{m \in \mathcal{M}_N}$ and $\{ Q_{IV}(\theta_m) \}_{m \in \mathcal{M}_N}$ are defined by (\ref{eigenfunction-Floquet-2}) for $\theta = \theta_m$ with the only exception for $\theta_0 = 0$, for which the proper eigenfunctions 
are $q_1$ and $q_2$ in (\ref{solutions-adjoint-extra}). Similar preparations are done for the adjoint eigenfunctions.  

The orthogonality relations of Propositions \ref{prop-orthog-discrete} and \ref{prop-orth-additional} are extended for these eigenfunctions and their adjoint eigenfunctions in $L^2_{\rm Nper}$. The completeness relation of Proposition \ref{prop-complete} is also extended in $L^2_{\rm Nper}$ with a similar proof. Since all additional eigenfunctions in $L^2_{\rm N per}$ correspond to non-zero eigenvalues $\Lambda$ and $\tilde{\Lambda}$, 
the diagonalization formula for $\langle (\tilde{L} + \delta L) q, q \rangle_{L^2_{\rm N per}}$ has strictly positive coefficients for all terms with the exception of the projection to $\partial_x U \in L^2_{\rm N per}$, for which the coefficient is zero.  Hence, ${\rm Ker}_{L^2_{\rm N per}} (\tilde{L} + \delta L) = {\rm span}(U')$ and the coercivity bound holds for some $C_N > 0$:
\begin{align*}
\langle (\tilde{L} + \delta L) q, q \rangle_{L^2_{\rm N per}} &\geq 
C_N \| q - c_1 \partial_x U \|_{H^1_{\rm N per}}^2.
\end{align*}
The proof of Theorem \ref{th-main} is complete for every fixed $N \in \mathbb{N}$. We note that $C_N \to 0$ as $N \to \infty$ due to the lack 
of the spectral gap in the second term in (\ref{diag-final}).

\bigskip

{\bf Acknowledgements.} The work of J. Chen was supported by the National Natural Science Foundation of China (No. 12371248), and the Project ``333'' of Jiangsu Province. The work of C. Grant was a part of an undergraduate thesis at McMaster University (Math 4P06). The work of D. E. Pelinovsky was supported by the NSERC Discovery Grant.

\appendix

\section{Proof of relations (\ref{operators-explicit-tilde}) for $Q_I(\lambda)$ and $Q_{III}(\lambda)$}
\label{app-A}

\underline{Relation (\ref{operators-explicit-tilde}) for $Q_I(\lambda)$.}
We recall from \eqref{eigenfunction-explicit-1} and \eqref{eigenfunction-adjoint-explicit-1} that $Q_I$ and $\Psi_I$ are analytic in $\mathbb{C}^+$ in the complex $x$-plane. Then, we have $H (Q_I') = i Q_I'$ and $H \Psi_I = i \Psi_I$. Applying the Hilbert operator to both sides of \eqref{b-2}, we obtain  
\begin{equation} \label{b-3}
Q_I' = 2 i U Q_I - i c Q_I + i (\lambda + c) \Psi_I,
\end{equation}
which immediately results in 
\begin{equation} \label{b-4}
Q_I'' = 2 i U' Q_I + 2 i U Q_I' - i c Q_I' - \lambda (\lambda + c) Q_I.
\end{equation}
Resorting to \eqref{c-2}, we arrive at 
\begin{equation}\label{b-5}
U' = i k^2 \left(\frac{1}{f'} - \frac{1}{f}\right) \left(\frac{1}{f'} + \frac{1}{f} - 1\right).
\end{equation}
By combining \eqref{Lax-spectrum-continuous}, \eqref{lambda-mu}, \eqref{Hirota}, \eqref{c-2}, \eqref{b-1}, \eqref{b-0}, \eqref{b-2},  \eqref{b-3}, \eqref{b-4}, and \eqref{b-5}, we compute 
\begin{align}
\tilde{L} Q_I  = & -Q_I'' + \frac{3}{2} Q_I H(U')  + \frac{3}{2} U 
H(Q_I') + \frac{3}{2} H(U Q_I)' + 3 U^2 Q_I + \gamma Q_I \nonumber \\
= & - \frac14 Q_I'' + \frac32 i U Q_I' + \frac34 i c Q_I' + \frac32 c U Q_I + \frac32 U^2 Q_I + \frac34 \lambda (\lambda + c) Q_I 
+ \gamma Q_I \nonumber \\
= & Q_I \left[ \lambda (\lambda + c) + \lambda_0 \mu_0  + \frac12 c U - \frac12 U^2 - \frac{1}{2} i U' \right] - \Psi_{I} (\lambda + c) (c + U). \label{b-6}
\end{align}
Substituting the expressions of \eqref{Lax-spectrum-continuous}, \eqref{lambda-mu}, \eqref{c-2}, \eqref{eigenfunction-explicit-1}, \eqref{eigenfunction-adjoint-explicit-1}, and \eqref{b-5} into \eqref{b-6}, we arrive at 
\begin{align*}
\tilde{L} Q_I  = & \frac{e^{i \lambda \xi}}{(\lambda - \mu_0) f'}\left[ \lambda^3 f - \lambda^2 \lambda_0 f' - \lambda f (k^2 + 3 \lambda_0 \mu_0)  - k (k^2 + 3 \lambda_0 \mu_0) + \mu_0 f (k^2 + 3 \lambda_0 \mu_0) \right] \\
= & \frac{e^{i \lambda \xi}}{(\lambda - \mu_0) f'} \left[ \lambda^2 (\lambda f - \lambda_0 f') - (k^2 + 3 \lambda_0 \mu_0) (\lambda f - \mu_0 f + k) \right] \\
= &  \frac{e^{i \lambda \xi}}{(\lambda - \mu_0) f'} \left[ \lambda^2 - (k^2 + 3 \lambda_0 \mu_0) \right] (\lambda f - \lambda_0 f') \\
=&  \frac{e^{i \lambda \xi}(\lambda^2 + \gamma)}{(\lambda - \mu_0)f'} [f(\lambda - \mu_0) + k] = (\lambda^2 + \gamma) \Psi_{I},
\end{align*}
which coincides with the first relation in (\ref{operators-explicit-tilde}),
where we have used the identity 
\begin{equation} \label{b-13}
k (f - 1) + \lambda_0 f = \lambda_0 f', \qquad \lambda_0 f' = \mu_0 f - k.
\end{equation}

\underline{Relation (\ref{operators-explicit-tilde}) for $Q_{III}(\lambda)$.}
Similar to the derivation of \eqref{b-6}, using \eqref{b-1}, \eqref{third-eig}, and \eqref{b-01}, we obtain
\begin{align}
\tilde{L} Q_{III}  = & -Q_{III}'' + \frac{3}{2} Q_{III} H(U')  + \frac{3}{2} U 
H(Q_{III}') + \frac{3}{2} H(U Q_{III})' + 3 U^2 Q_{III} + \gamma Q_{III} \nonumber \\
= & - \frac14 Q_{III}'' + \frac32 c U Q_{III} - \frac32 U^2 Q_{III} - \frac14  k^2 Q_{III} \nonumber  \\ & - \frac34 \lambda^{-1} 
(\lambda - \lambda_0) (\lambda - \mu_0) \left[2 U \Psi_{III} + c \Psi_{III} + H(\Psi_{III}') \right], \label{b-8} 
\end{align}
Resorting to \eqref{lambda-mu}, \eqref{Hirota}, and \eqref{b-02}, we have
\begin{align}
H(\Psi_{III}') = & \frac{d }{d \xi } \left[ i e^{i (\lambda - \lambda_0) \xi} \left( \frac{\lambda }{(\lambda -\mu_0) f'} + \frac{\lambda_0}{(\lambda -\mu_0) f}
\right) \right] \nonumber \\
=& - e^{i (\lambda - \lambda_0) \xi} \left( \frac{\lambda}{f'} + \frac{\lambda_0}{f}\right) - \frac{k e^{i (\lambda - \lambda_0) \xi} }{\lambda - \mu_0} \left[
\frac{\lambda}{(f')^2} + \frac{\lambda_0}{f^2} \right]. \label{b-10}
\end{align}
By using \eqref{lambda-mu}, \eqref{Hirota}, and  \eqref{eigenfunction-explicit-2}, we arrive at 
\begin{align}
Q_{III}' &= i (\lambda - \mu_0) Q_{III} + \frac{i k e^{i (\lambda - \lambda_0) \xi} }{\lambda - \mu_0} \left[
\frac{\lambda - \mu_0 - k}{(f')^2} + \frac{2 k }{(f')^3} - \frac{\lambda_0}{\lambda} \left(
\frac{\lambda - \mu_0 - k}{f^2} + \frac{2 k }{f^3}  \right) \right], 
\label{b-11} \\ 
Q_{III}'' &=   - \frac{k (\lambda - \lambda_0) e^{i (\lambda - \lambda_0) \xi} }{\lambda - \mu_0} \left[
\frac{\lambda - \mu_0 - k}{(f')^2} + \frac{2 k }{(f')^3} - \frac{\lambda_0}{\lambda} \left(
\frac{\lambda - \mu_0 - k}{f^2} + \frac{2 k }{f^3}  \right) \right]  \nonumber \\
& + i (\lambda - \mu_0) Q_{III}' - \frac{2 k^2 e^{i (\lambda - \lambda_0) \xi} }{\lambda - \mu_0} \left[ \frac{ \lambda - \mu_0 - 4 k }{(f')^3}
- \frac{\lambda - \mu_0 - k}{(f')^2} + \frac{3 k}{(f')^4} \right. \nonumber \\ & \left. - \frac{\lambda_0}{\lambda} \left(
\frac{ \lambda - \mu_0 - 4 k }{f^3}
- \frac{\lambda - \mu_0 - k}{f^2} + \frac{3 k}{f^4}\right) \right]. \label{b-12}
\end{align}
Once again, plugging \eqref{Lax-spectrum-continuous}, \eqref{lambda-mu}, \eqref{c-2}, \eqref{eigenfunction-explicit-2}, \eqref{eigenfunction-adjoint-explicit-2}, \eqref{b-13}, \eqref{b-10}, \eqref{b-11}, and \eqref{b-12} into \eqref{b-8}, a tedious but direct algebraic computation gives rise to
\begin{align*}
\tilde{L} Q_{III}  = & \frac{e^{i (\lambda - \lambda_0) \xi} }{ \lambda (\lambda - \mu_0 ) f^3 (f')^3} \left[  \lambda^4 f^3 (f')^2 -  \lambda^3 \lambda_0 f^2 (f')^3 
-  \lambda^2 f^3 (f')^2 (k^2 + 3 \lambda_0 \mu_0) \right.  \nonumber \\ & \left. 
+  \lambda \lambda_0 f^2 (f')^2 [(\lambda_0 + \mu_0) k + f' (4 \lambda_0 \mu_0 + k^2 + \lambda_0^2)] -  \lambda_0^2 \mu_0 (2 \lambda_0 + k) f^2 (f')^3 \right]  \nonumber \\
= & \frac{e^{i (\lambda - \lambda_0) \xi} }{\lambda f f' (\lambda - \mu_0)} \left\{ k\left[ \lambda^3 - \lambda (k^2 + 3 \lambda_0 \mu_0) + \lambda_0 \mu_0 (\lambda_0 + \mu_0)  \right]   \right. \nonumber \\ & 
\left. + f \left[ \lambda^4 - \mu_0 \lambda^3 - \lambda^2 (k^2 + 3 \lambda_0 \mu_0)  + \lambda \mu_0 (k^2 + \lambda_0^2 + 4 \lambda_0 \mu_0) -\lambda_0 \mu_0^2 (\lambda_0 + \mu_0)	\right] \right\} \nonumber \\
= & e^{i (\lambda - \lambda_0) \xi} \frac{(\lambda - \lambda_0) (\lambda - \mu_0) (\lambda - c)}{\lambda f f' (\lambda - \mu_0)} \left[
f (\lambda - \mu_0) + k \right] \nonumber \\
= & \frac{(\lambda - \lambda_0) (\lambda - \mu_0) (\lambda - c)}{\lambda} \Psi_{III},
\end{align*}
which coincides with the third relation in (\ref{operators-explicit-tilde}).

\end{document}